\documentclass[a4paper]{article}
\usepackage[a4paper, total={6in, 9in}]{geometry}

\usepackage{amsmath,amsthm}
\usepackage{amssymb}

\usepackage{enumitem}

\usepackage{caption}
\usepackage{array} 
\newcolumntype{M}[1]{>{\centering\arraybackslash}m{#1}} 
\usepackage{stackengine} 

\makeatletter
\newcommand{\customlabel}[2]{%
	\protected@write \@auxout {}{\string \newlabel {#1}{{#2}{}}}}
\makeatother
\newenvironment{proofbold}[1]{\paragraph{Proof of {#1}.}}{\hfill$\square$}

\newtheorem{theorem}{Theorem}[section]

\newtheorem{lemma}[theorem]{Lemma}
\newtheorem{proposition}[theorem]{Proposition}

\theoremstyle{definition}

\newtheorem{remark}[theorem]{Remark}

\newtheorem{innercustomthm}{Theorem}
\newenvironment{ctheorem}[1]
  {\renewcommand\theinnercustomthm{#1}\innercustomthm}
  {\endinnercustomthm}

\numberwithin{equation}{section}

\newtheorem*{theorem*}{Theorem}

\usepackage{sagetex}

\usepackage{dsfont}

\DeclareMathAlphabet{\mathpzc}{OT1}{pzc}{m}{it}

\newenvironment{mysage}{\sagesilent}{\endsagesilent}

\begin{document}

\begin{sagesilent}
################
### Rounding ####
###############

def roundup(x,d):

    return float(ceil(x*10^d)/10^d)
def rounddown(x,d):
    return float(floor(x*10^d)/10^d)

##################################
### Managing Number Representation ###
##################################

def Interval(x): ### Converts a real x into a neighborhood of x, expressed by a closed interval. Allows us to handle interval arithmetic. ###
    X = RBF(RIF(x)) 
    I = RealSet([X.lower(),X.upper()])[0]
    return I

def Upper(x,d): ### Gives a ''close'' upper bound with d decimal digits for a real x. ###
    I = Interval(x)
    u = float(ceil(I.upper()*10^d) /10^d)
    return u

def Lower(x,d): ### Gives a ''close'' lower bound with d decimal digits for a real x. ###
    I = Interval(x)
    l = float(floor(I.lower()*10^d)/10^d)
    return l

def Trunc(x,d): ### Given a real x, it gives an expression of x with at most d EXACT decimals; handy for managing displaying of numbers without losing accuracy (upon inserting ''...''). It is a truncation that manages to take care of short decimals as well as numbers of the form 2.9999999999999 (that are taken to be 3). It handles other numbers normally. ###

    if d>10:
        return false
    I = Interval(x)
    s=0
    while (I.lower()*10^s).trunc()==(I.upper()*10^s).trunc() and s<=10:
        s=s+1
    if (x-RIF(x))!=0: ### Ex. x=29.99999999999999999 => RIF(x)=30. This operation returns 29.99999 (d=5), rather than 30 ###
        return Lower(x,d)
    if d>s and (x-RIF(x))==0: ### Ex. x=1.24 and d=5. This operation returns 1.24 rather than 1.23999 given by above orders. ###
        return float(x)
    return float(floor(I.lower()*10^d)/10^d) 

def Numb(x,y,d):
    z=RIF((x+y)/2)
    return Trunc(z,d)

gamma = 0.5772156649015328606065120900824024310421

def A(a,d):
    constantpaper = 1
    if a==1:
        if 0<d<=1:
            return RIF(max(gamma,1/d/exp(gamma*d+1)))
        return False
    if a!=1:
        if d+1==a:
            return RIF(max(constantpaper,1/abs(a-1), zeta(a)-1/(a-1)))
        if 0<d<a<d+1:
            l1=RIF((d-a+1)/abs(zeta(a))/abs(a-1))
            return RIF(max(constantpaper, (l1^(d-a+1)/d^d)^(1/(a-1)), zeta(a)-1/(a-1)))             
        if d==a:
            return RIF(constantpaper)
        return False

def E(a,v):
    if v!=1 and v!=2:
        return False
    e1=0.43*(1+abs(a-1)/(a-1/2))
    e2=abs(zeta(a)/zeta(2*a)-6/(a-1)^2/pi^2) 
    e3=abs(a-1)/(a-1/2) * ( 3*zeta(2*a) / ( (a-1/2)*pi^2*abs(zeta(a)*(a-1)) ) )^(2/(a-1)) 
    if v==2:
        e1=0.12*(1+abs(a-1)/(a-1/2))
        e2=(sqrt(2)-1)/sqrt(2) * abs(2^a/(2^a+1) * zeta(a)/zeta(2*a)-2/(2+1) * 6/(a-1)^2/pi^2)
        e3=(sqrt(2)-1)/sqrt(2) * abs(a-1)/(a-1/2) * ( 3*(2^a+1)*zeta(2*a) / ( (a-1/2)*2^(a-1)*(2+1)*pi^2*abs(zeta(a)*(a-1)) ) )^(2/(a-1))  
    E= RIF(max(e1,e2,e3))
    return E
	
\end{sagesilent}

\begin{mysage} 

################################################################################# MAIN PARAMETERS ###

digits = 3 #Precision seems to be affected by more than 7 digits

dlong = 6 # (Trunc(,) is not defined for higher number of digits)

choice = 2/3

PW = 10^(12)
theta = RIF(1-1/12/log(10))

program = 5*10^8 # expected boundary

provisoire = 10^7 # current boundary

bound = 8

exp_prog = 7
exp_prog_bound = 25/2

DET = 10^exp_prog
TRICK = 10^exp_prog_bound

BOU = 10^bound

################################################################################### tildetildem ###

I_sum1_l=0.755366607315099
I_sum1_u=0.755366626776258	### Precision: 10^9, Time: 8555.87245297s, using SAGE ###
C_sum1=Numb(I_sum1_l+gamma, I_sum1_u+gamma,8)

k2_1half = RIF ( 3/(3-3^(1-1/2)) )
k2_theta = RIF ( 3/(3-3^(1-theta)) )

Delta_1half_l = zeta(3/2)*1.36843276585094 	### Precision: 4*10^9, using C++ ###
Delta_1half_u = zeta(3/2)*1.36843284087041
Delta_1half = Delta_1half_u

Delta_theta_l = zeta(1+theta)*1.00724550163589
Delta_theta_u = zeta(1+theta)*1.00724557626645		
Delta_theta = Delta_theta_u

################################################################################# Lemma_{sum2} ###

I_sum2_l=1.13992197915589		### Precision: 10^9, Time: 8555.87245297s, using SAGE ###
I_sum2_u=1.13992201807822		
C_sum2=Numb(I_sum2_l+gamma, I_sum2_u+gamma,digits)

################################################################################# Lemma_{sum1} ###

constant_ram = 4.4
constant_ram_v2 = 0.493

I_sum1_l=0.755366607315099
I_sum1_u=0.755366626776258	
C_sum1=Numb(I_sum1_l+gamma, I_sum1_u+gamma,digits)

################################################################################# Lemma_{sumvar1log} ###

delta_sumvar1log = 1/3

Prod_sumvar1log_l = 1.6035247003888
Prod_sumvar1log_u = 1.60354379786847

Sum_sumvar1log_l = -0.19526982241707
Sum_sumvar1log_u = -0.194491292426689

deltaProd_sumvar1log_l = 33.792832437589
deltaProd_sumvar1log_u = 39.522139388739

errProd_sumvar1log  = A(1,delta_sumvar1log) * deltaProd_sumvar1log_u 

constant_A2 = 1
Tau3_v2 = 1 - constant_A2 / (2-1+constant_A2)
constant_sg2 =  log(2)  * constant_A2 / (constant_A2+2-1)
Tau4_v2 = 1 + (2*(2-1)-constant_A2*(2+2^delta_sumvar1log) ) / ( (2-1)*2^(1-delta_sumvar1log) + constant_A2*(2+2^delta_sumvar1log) - 2 +1)

analytic_sumvar1log_v1 = Prod_sumvar1log_u * ( 1/2 + (Sum_sumvar1log_u+gamma) / log(provisoire) ) + errProd_sumvar1log  / delta_sumvar1log / log(provisoire)^2
analytic_sumvar1log_v2 = Prod_sumvar1log_u * Tau3_v2 * ( 1/2 + (Sum_sumvar1log_u+gamma+constant_sg2) / log(provisoire) ) + errProd_sumvar1log  * Tau4_v2 / delta_sumvar1log / log(provisoire)^2

program_sumvar1log_v1 = 1.09937069175235
program_sumvar1log_v2 = 0.694356698566237 

constant_eta_v1 = max( analytic_sumvar1log_v1, program_sumvar1log_v1 )
constant_eta_v2 = max( analytic_sumvar1log_v2, program_sumvar1log_v2 )

################################################################################# Lemma_{S1:total} ###

Upsilon1_v1 =  2 * constant_ram * constant_ram_v2 * constant_eta_v2 + 2 * constant_ram * constant_eta_v2
Upsilon1_v2 =  2 * constant_ram * constant_ram_v2 * constant_eta_v2

################################################################################# Lemma_{sumvarp} ###

provisoire1 = 10^6  # CORRECT HERE

Prod_sumvarp_l =  1.94359643387259  ### Precision: 3*10^9, Time: <360s, using C++ ---> int_double
Prod_sumvarp_u = 1.94359649909918

uu2 = RIF(1 - 2/(2^2-2+1) )
vv2 = RIF(1+(2^2-4*2+2)/((sqrt(2)-1)*(2-1)^2+2*2-1) )
 
errProd_sumvarp_l = 60.5407619086305
errProd_sumvarp_u = 60.5425689714784

value_sumvarp_f2 = 1 / (2-1)^2
walue__sumvarp_1 = ( sqrt(2)-1) / ( sqrt(2)-1+abs(2^2*value_sumvarp_f2-1) ) * ( E(2,1) + E(2,2) * abs(2^2*value_sumvarp_f2-1) / (sqrt(2)-1) ) 
walue__sumvarp_2 = E(2,2)

i_v1_l = walue__sumvarp_1 * errProd_sumvarp_l
i_v1_u = walue__sumvarp_1 * errProd_sumvarp_u

i_v2_l = walue__sumvarp_2 * errProd_sumvarp_l
i_v2_u =  walue__sumvarp_2 * errProd_sumvarp_u 

a_constant_phi_v1 = Prod_sumvarp_u + i_v1_u / sqrt(provisoire1)
a_constant_phi_v2 = uu2 * Prod_sumvarp_u + vv2 * i_v2_u / sqrt(provisoire1)

constant1_threshold_sumvarp = 2.37215689830106
constant2_threshold_sumvarp = 0.725107915681235

constant1_phi_v1 = max(a_constant_phi_v1,constant1_threshold_sumvarp)
constant1_phi_v2 = max(a_constant_phi_v2,constant2_threshold_sumvarp) 

################################################################################# Lemma_{Ss1} ###

provisoire2 = 10^8 # CORRECT HERE (uses phi)

a_constant_phi2_v1 = Prod_sumvarp_u + 5 * i_v1_u / sqrt(provisoire2)
a_constant_phi2_v2 = uu2 * Prod_sumvarp_u + 5 * vv2 * i_v2_u / sqrt(provisoire2)

constant1_threshold_Ss1 = 1.94386748623864
constant2_threshold_Ss1 = 0.647934571702803

constant2_phi_v1 = max(a_constant_phi2_v1,constant1_threshold_Ss1)
constant2_phi_v2 = max(a_constant_phi2_v2,constant2_threshold_Ss1) 

################################################################################# Lemma_{Parity} ###

Cc1 = 1.044
Cc2 = 0.232

ParityProd_l = 1.51478012597114
ParityProd_u = 1.51478017037289
 
errParityProd_l = 11.0518744225813
errParityProd_u = 11.0522043537281

value_Parity_f2 = 1
walue__Parity_1 = ( sqrt(2)-1) / ( sqrt(2)-1+abs(2*value_Parity_f2-1) ) * ( Cc1 + Cc2 * abs(2*value_Parity_f2-1) / (sqrt(2)-1) ) 
walue__Parity_2 = Cc2

h_v1_l = walue__Parity_1 * errParityProd_l
h_v1_u = walue__Parity_1 * errParityProd_u

h_v2_l = walue__Parity_2 * errParityProd_l
h_v2_u =  walue__Parity_2 * errParityProd_u

################################################################################# Lemma_{S2:total} ###

Prod_Parity_l = 0.660161800282638  ### Precision: 3*10^9, Time: <360s, using C++ ---> int_double
Prod_Parity_u = 0.660161816820513

Upsilon2_v1 = Prod_Parity_u * sqrt(2) * h_v2_u
Upsilon2_v2 = 4 * Prod_Parity_u * h_v2_u

Upsilon3_v1 = constant1_phi_v1 * constant2_phi_v1
Upsilon3_v2 = constant1_phi_v2 * constant2_phi_v2 * 2^2

################################################################################# Lemma_{sum.1/2} ###

delta_sumHalf = 1/3

Prod_sumHalf_l = 15.0333977306198  ### Precision: 3*10^9, Time: <360s, using C++ ---> int_double
Prod_sumHalf_u = 15.0337644348976

deltaProd1_sumHalf_l = 10.7125980293757
deltaProd1_sumHalf_u = 10.7200446392398

deltaProd2_sumHalf_l = 150.759502870216
deltaProd2_sumHalf_u = 165.591131996911

deltaProd3_sumHalf_l = deltaProd1_sumHalf_l * deltaProd2_sumHalf_l
deltaProd3_sumHalf_u = deltaProd1_sumHalf_u * deltaProd2_sumHalf_u

Sum_sumHalf_l = -1.73634511834745+gamma
Sum_sumHalf_u = -1.73581610256903+gamma

Error_sumHalf = A(1,delta_sumHalf) * deltaProd3_sumHalf_u

Error_sum2Half = Error_sumHalf * (1 +  1/(1-delta_sumHalf) )

################################################################################# Lemma_{sum.1/2threshold} ###

F2 = RIF( 1- 1/(2-2*sqrt(2)+2) )
G2 = RIF( log(2) / (2-2*sqrt(2)+2) )
H2 = RIF( 1+(2-4*sqrt(2)-2^delta_sumHalf+2) / ( (sqrt(2)-1)^2*2^(1-delta_sumHalf)+2*sqrt(2)+2^delta_sumHalf-1 ) )

constant1_threshold_1half = 14.087466570648  
constant2_threshold_1half =  2.1362281243806

a_constant1_chi_v1 = Prod_sumHalf_u + Error_sumHalf / program^delta_sumHalf / log(program)
a_constant1_chi_v2 = F2 * Prod_sumHalf_u + H2 * Error_sumHalf / program^delta_sumHalf / log(program)

constant1_chi_v1 = max( constant1_threshold_1half, a_constant1_chi_v1)
constant1_chi_v2 = max( constant2_threshold_1half, a_constant1_chi_v2)

################################################################################# Lemma-label{sum2.1/2}  ###

constant1_threshold2_1half = 16.7682417501771
constant2_threshold2_1half = 2.50640699370728 

a_constant2_chi_v1 = Prod_sumHalf_u + Error_sum2Half / program^delta_sumHalf 
a_constant2_chi_v2 = F2 * Prod_sumHalf_u + H2 * Error_sum2Half / program^delta_sumHalf 

constant2_chi_v1 = max( constant1_threshold2_1half, a_constant2_chi_v1)
constant2_chi_v2 = max( constant2_threshold2_1half, a_constant2_chi_v2)

################################################################################# Lemma-label{Ss1Log}   ###

f2_2_Ss1Log = RIF( 1 - 1/(2^(2-2*theta)*(2^theta-1)^2+1) )
f2_3_Ss1Log = RIF( 1 - 1/(2^(3/2-2*theta)*(2^theta-1)^2+1) )

Prod2_Ss1Log_l = 2.93392267813336  ### Precision: 3*10^9, Time: <360s, using C++ ---> int_double
Prod2_Ss1Log_u = 2.93392292726709

Prod3_Ss1Log_l = 5.30171336966289 ### Precision: 3*10^9, Time: <360s, using C++ ---> int_double
Prod3_Ss1Log_u = 5.30175957628429

constant_tau_v1 = Prod2_Ss1Log_u + ( 4/e/log(PW) + 16/e^2/log(PW)^2 ) * Prod3_Ss1Log_u 
constant_tau_v2 = f2_2_Ss1Log * Prod2_Ss1Log_u + ( 4/e/log(PW) + 16/e^2/log(PW)^2 ) * f2_3_Ss1Log * Prod3_Ss1Log_u 

################################################################################# Lemma-label{Ss2Log}   ###

Prod_Ss2Log_l = 2.08813954913666  ### Precision: 3*10^9, Time: <360s, using C++ ---> int_double
Prod_Ss2Log_u = 2.0881396785312

Sum_Ss2Log_l = 0.282661139651601+gamma  ### Precision: 3*10^9, Time: <360s, using C++ ---> int_double
Sum_Ss2Log_u = 0.282661147670366+gamma

xx2 = RIF( 1-(2-1) / ( 2^(2-2*theta)*(2^theta-1)^2+2*2^(1-theta)-2^(1-2*theta)-1 ) )
ss2 = RIF( log(2)/( 2^(1-2*theta)*(2^theta-1)^2+1 ) )
yy2 = RIF( 1 + ( 2^(2*theta)-4*2^theta+2 )/( (sqrt(2)-1)*(2^theta-1)^2 + 2*2^theta-1 ) )

errProd_Ss2Log_l = 74.7209571772
errProd_Ss2Log_u = 74.725281409463

value_Ss2Log_f2 = 1 / 2^(1-2*theta) / (2^theta-1)^2
walue_Ss2Log_1 = ( sqrt(2)-1) / ( sqrt(2)-1+abs(2*value_Ss2Log_f2-1) ) * ( Cc1 + Cc2 * abs(2*value_Ss2Log_f2-1) / (sqrt(2)-1) ) 
walue_Ss2Log_2 = Cc2

j_v1_l = walue_Ss2Log_1 * errProd_Ss2Log_l
j_v1_u = walue_Ss2Log_1 * errProd_Ss2Log_u

j_v2_l = walue_Ss2Log_2 * errProd_Ss2Log_l
j_v2_u =  walue_Ss2Log_2 * errProd_Ss2Log_u

kk = 1/2

ccc = 10
ccc2 =  16
ccc3 = 70

lll = log(PW^(1-choice) / ccc)
lll2 = log(PW^(1-choice) / ccc2)
lll3 = log(PW^(1-choice) / ccc3)

sqqq = log(PW^(1/kk-choice) / ccc )
sqqq2 = log(PW^(1/kk-choice) / ccc2 )
sqqq3 = log(PW^(1/kk-choice) / ccc3 )

WWW = 1/ sqrt(ccc) / sqrt(PW^choice) + 4 / ( kk * sqqq ) + 4 / ( ccc^(kk/2) * PW^(kk/2 * choice) * lll)
WWW2 = 1/ sqrt(ccc2) / sqrt(PW^choice) + 4 / ( kk * sqqq2 ) + 4 / ( ccc2^(kk/2) * PW^(kk/2 * choice) * lll2)
WWW3 = 1/ sqrt(ccc3) / sqrt(PW^choice) + 4 / ( kk * sqqq3 ) + 4 / ( ccc3^(kk/2) * PW^(kk/2 * choice) * lll3)

YYY_v1 = j_v1_u * WWW / lll / ( 1-choice-log(ccc)/log(PW) )
YYY2_v1 = j_v1_u * WWW2 / lll2 / ( 1-choice-log(ccc2)/log(PW) )
YYY3_v1 = j_v1_u * WWW3 / lll3 / ( 1-choice-log(ccc3)/log(PW) )

YYY_v2 = yy2 * j_v2_u * WWW / lll / ( 1-choice-log(ccc)/log(PW) )
YYY2_v2 = yy2 * j_v2_u * WWW2 / lll2 / ( 1-choice-log(ccc2)/log(PW) )
YYY3_v2 = yy2 * j_v2_u * WWW3 / lll3 / ( 1-choice-log(ccc3)/log(PW) )

ttt = log(ccc) / log(PW) + choice
ttt2 = log(ccc2) / log(PW) + choice
ttt3 = log(ccc3) / log(PW) + choice

XXX_v1 = Prod_Ss2Log_u * ( 1 / ttt + Sum_Ss2Log_u / log(PW) )
XXX2_v1 = Prod_Ss2Log_u * ( 1 / ttt2 + Sum_Ss2Log_u / log(PW) )
XXX3_v1 = Prod_Ss2Log_u * ( 1 / ttt3 + Sum_Ss2Log_u / log(PW) )

XXX_v2 = xx2 * Prod_Ss2Log_u * ( 1 / ttt + (Sum_Ss2Log_u + ss2 ) / log(PW) )
XXX2_v2 = xx2 * Prod_Ss2Log_u * ( 1 / ttt2 + (Sum_Ss2Log_u + ss2 ) / log(PW) )
XXX3_v2 = xx2 * Prod_Ss2Log_u * ( 1 / ttt3 + (Sum_Ss2Log_u + ss2 ) / log(PW) )

constant_xi_v1 = XXX_v1 + YYY_v1
constant2_xi_v1 = XXX2_v1 + YYY2_v1
constant3_xi_v1 = XXX3_v1 + YYY3_v1

constant_xi_v2 = XXX_v2 + YYY_v2
constant2_xi_v2 = XXX2_v2 + YYY2_v2
constant3_xi_v2 = XXX3_v2 + YYY3_v2

################################################################################# Theorem-label{S3:total}   ###

boundomega_v1_1 = constant1_chi_v1 * constant2_chi_v1
boundomega_v2_1 = constant1_chi_v2 * constant2_chi_v2 * 2 / (sqrt(2)-1)^2 

boundomega_v1_2 = 1 / 389 * constant_tau_v1 * constant_xi_v1 
boundomega_v2_2 = 1 / 389 * constant_tau_v2 * constant_xi_v2 * 2^(2*theta) / (2^theta-1)^2

boundomega2_v1_2 = 1 / 389 * constant_tau_v1 * constant2_xi_v1 
boundomega2_v2_2 = 1 / 389 * constant_tau_v2 * constant2_xi_v2 * 2^(2*theta) / (2^theta-1)^2

boundomega3_v1_2 = 1 / 389 * constant_tau_v1 * constant3_xi_v1 
boundomega3_v2_2 = 1 / 389 * constant_tau_v2 * constant3_xi_v2 * 2^(2*theta) / (2^theta-1)^2

optimal_omega_v1 =  1 / (constant1_chi_v1 * constant2_chi_v1 ) 
optimal_omega_v2 = (sqrt(2)-1)^2 /  (constant1_chi_v2 * constant2_chi_v2 * 2 *2  ) # we are considering v^1

value_c_v1 = optimal_omega_v1 # IT MAY CHANGE
value_c_v2 = optimal_omega_v2

Upsilon4_v1 = (1+value_c_v1) * constant1_chi_v1 * constant2_chi_v1 
Upsilon4_v2 = (1+value_c_v2) * constant1_chi_v2 * constant2_chi_v2 * 2 / (sqrt(2)-1)^2 

Upsilon5_v1 = (1+1/value_c_v1) / 389^2 * constant_tau_v1 * constant_xi_v1 
Upsilon5_v2 = (1+1/value_c_v2) / 389^2 * constant_tau_v2 * constant_xi_v2 * 2^(2*theta) / (2^theta-1)^2

Upsilon5_v1_2 = (1+1/value_c_v1) / 389^2 * constant_tau_v1 * constant2_xi_v1 
Upsilon5_v2_2 = (1+1/value_c_v2) / 389^2 * constant_tau_v2 * constant2_xi_v2 * 2^(2*theta) / (2^theta-1)^2

Upsilon5_v1_3 = (1+1/value_c_v1) / 389^2 * constant_tau_v1 * constant3_xi_v1 
Upsilon5_v2_3 = (1+1/value_c_v2) / 389^2 * constant_tau_v2 * constant3_xi_v2 * 2^(2*theta) / (2^theta-1)^2

################################################################################# Lemma-label{sq1.half}   ###

delta_sqHalf = 1/3 

kk2 = 1 - ( 2+1 ) / ( 2^2-2^(3/2)+sqrt(2)+1 )
sum_ff2 = log(2) / ( 2-sqrt(2)+1 )
ll2 = 1 + ( 2^(1-delta_sqHalf) -2*2^(1/2-delta_sqHalf)-1 ) / ( 2^(2-2*delta_sqHalf)-2^(3/2-2*delta_sqHalf)+2^(1/2-delta_sqHalf))

Prod_sqHalf_l =  zeta(3/2)/zeta(3)
Prod_sqHalf_u = zeta(3/2)/zeta(3)

delta_sqHalf_l = 51.3170892955494
delta_sqHalf_u = 59.1744268501228

err_sqHalf = A(1,delta_sqHalf ) * delta_sqHalf_u 

Sum_sqHalf_l = -0.210028802803383+gamma
Sum_sqHalf_u = -0.20792551024676+gamma

analytic_sqHalf_v1 = zeta(3/2)/zeta(3) * ( 1/2 + Sum_sqHalf_u / log(provisoire) ) + err_sqHalf / log(provisoire)^2
analytic_sqHalf_v2 = kk2 * zeta(3/2)/zeta(3) * ( 1/2 + (Sum_sqHalf_u+sum_ff2) / log(provisoire) ) + ll2 * err_sqHalf / log(provisoire)^2

program_sqHalf_v1 = 1.4256628496167
program_sqHalf_v2 = 0.711297411532317

constant_psi_v1 = max( analytic_sqHalf_v1, program_sqHalf_v1 )
constant_psi_v2 = max( analytic_sqHalf_v2, program_sqHalf_v2 )

################################################################################# label{source}   ###

coeff_r_v1 = RIF(2* Cc2/ Cc1 / sqrt(2)/(sqrt(2)-1))

Tau1_v1 =  RIF(2 * Cc1 * (coeff_r_v1 + 1) * (constant_psi_v1 * constant_psi_v2 + constant_psi_v1^2))
Tau1_v2 =  RIF(2 * Cc2 * sqrt(2) / (sqrt(2)-1) * constant_psi_v2^2)

################################################################################# label{integral}   ###

Tau2_v1 = 3.83717 #it should be improved for X>=100
Tau3_v1 = 4.89606 # =
Tau4_v1 = 0.000033536

Tau2_v2 = 4.99703 #it should be improved for X>=100
Tau3_v2 = 9.57182 # = 
Tau4_v2 = 0.0000615022

################################################################################# CONCLUSION   ###

lower = 20
Psi_v1 = Tau2_v1 + Tau3_v1 / log(lower)
Psi_v2 = Tau2_v2 + Tau3_v2 / log(lower)

Merge1_v1 = RIF( Tau1_v1 * (1-choice)^4 / sqrt(ccc)  + 2 * Upsilon1_v1 * (1-choice)^2 / sqrt(ccc) / log(DET)^2 + 6 * (1-choice) * ccc * Tau2_v1 / pi^2 / log(DET)^3 + (1-choice) * Upsilon4_v1 / log(DET)^3)
Merge1_v2 = RIF (Tau1_v2 * (1-choice)^4 / sqrt(ccc)  + 2 * Upsilon1_v2 * (1-choice)^2 / sqrt(ccc) / log(DET)^2 + 6 * (1-choice) * ccc * Tau2_v2 / pi^2 / log(DET)^3 * 2/(2+1) + (1-choice) * Upsilon4_v2 / log(DET)^3 )

Merge2_v1 = RIF( 1 / log(DET)^4 * ( 6 / pi^2 * ccc * (Tau2_v1 + Tau3_v1 ) + Upsilon2_v1 / sqrt(ccc) + ccc * Upsilon3_v1 + Upsilon4_v1 ) )
Merge2_v2 = RIF( 1 / log(DET)^4 * ( 6 / pi^2 * ccc * (Tau2_v2 + Tau3_v2 ) * 2/(2+1) + Upsilon2_v2 / sqrt(ccc) + ccc * Upsilon3_v2 + Upsilon4_v2 ) )

Alt_Merge1_v1 = RIF( Tau1_v1 * (1-choice)^4 / sqrt(ccc)  + 2 * Upsilon1_v1 * (1-choice)^2 / sqrt(ccc) / log(DET)^2 + 6 * (1-choice) * ccc * Psi_v1 / pi^2 / log(DET)^3 + (1-choice) * Upsilon4_v1 / log(DET)^3)
Alt_Merge1_v2 = RIF (Tau1_v2 * (1-choice)^4 / sqrt(ccc)  + 2 * Upsilon1_v2 * (1-choice)^2 / sqrt(ccc) / log(DET)^2 + 6 * (1-choice) * ccc * Psi_v2 / pi^2 / log(DET)^3 * 2/(2+1) + (1-choice) * Upsilon4_v2 / log(DET)^3 )

Alt_Merge2_v1 = RIF( 1 / log(DET)^4 * ( 6 / pi^2 * ccc * Psi_v1 + Upsilon2_v1 / sqrt(ccc) + ccc * Upsilon3_v1 + Upsilon4_v1 ) )
Alt_Merge2_v2 = RIF( 1 / log(DET)^4 * ( 6 / pi^2 * ccc * Psi_v2 * 2/(2+1) + Upsilon2_v2 / sqrt(ccc) + ccc * Upsilon3_v2 + Upsilon4_v2 ) )

MergeLog_v1 = RIF( 12 * Tau4_v1 / (1-choice-log(ccc)/log(DET) ) / pi^2 + Upsilon5_v1 )
MergeLog_v2 = RIF( 12 * Tau4_v2 / (1-choice-log(ccc)/log(DET) ) / pi^2 * 2/(2+1) + Upsilon5_v2 )

NUM1_v1 = RIF( Tau1_v1 * (1-choice)^4 / sqrt(ccc3) * log(TRICK/ccc3^(1/(1-choice)))^4 / TRICK^(choice/2) + 2 * Upsilon1_v1 * (1-choice)^2 / sqrt(ccc3) * log(TRICK/ccc3^(1/(1-choice)))^2 / TRICK^(choice/2) + 6 * (1-choice) * ccc3 * Psi_v1 / pi^2 * log(TRICK/ccc3^(1/(1-choice))) / TRICK^(1-choice) + (1-choice) * Upsilon4_v1 * log(TRICK) / TRICK^(1-choice) )
NUM1_v2 = RIF( Tau1_v2 * (1-choice)^4 / sqrt(ccc2) * log(TRICK/ccc2^(1/(1-choice)))^4 / TRICK^(choice/2) + 2 * Upsilon1_v2 * (1-choice)^2 / sqrt(ccc2) * log(TRICK/ccc2^(1/(1-choice)))^2 / TRICK^(choice/2) + 6 * (1-choice) * ccc2 * Psi_v2 / pi^2 * log(TRICK/ccc2^(1/(1-choice))) / TRICK^(1-choice) + (1-choice) * Upsilon4_v2 * log(TRICK) / TRICK^(1-choice) )

NUM2_v1 = RIF( 6 / pi^2 * ccc3 * Psi_v1 / TRICK^(1-choice) + Upsilon2_v1 / sqrt(ccc3) / TRICK^(choice/2)+ ccc3 * Upsilon3_v1 / TRICK^(1-choice)+ Upsilon4_v1 / TRICK^(1-choice)  )
NUM2_v2 = RIF( 6 / pi^2 * ccc2 * Psi_v2 * 2/(2+1) / TRICK^(1-choice) + Upsilon2_v2 / sqrt(ccc2) / TRICK^(choice/2)+ ccc2 * Upsilon3_v2 / TRICK^(1-choice)+ Upsilon4_v2 / TRICK^(1-choice)  )

NUMLog_v1 = RIF( 12 * Tau4_v1 / (1-choice-log(ccc3)/log(TRICK) ) / pi^2 / log(TRICK) + Upsilon5_v1_3 / log(TRICK) ) 
NUMLog_v2 = RIF( 12 * Tau4_v2 / (1-choice-log(ccc2)/log(TRICK) ) / pi^2 * 2/(2+1)  / log(TRICK) + Upsilon5_v2_2 / log(TRICK) )

numerical_v1 = NUM1_v1 + NUM2_v1 + NUMLog_v1
numerical_v2 = NUM1_v2 + NUM2_v2 + NUMLog_v2

KV1 = numerical_v1 * log(TRICK)
KV2 = numerical_v2 * log(TRICK)

################################################################################# CONSTANT   ###

barrier_threshold_1 = 0.607398570962174
barrier_threshold_2 = 1.4731118309395

B_constant_integral1_l = -0.0495100113498
B_constant_integral1_u = -0.049510010626

B_constant_integral2_l = 2.63481269161
B_constant_integral2_u =  2.63481271383

B_error_integral_v1 =  RIF( Tau2_v1 * log(BOU)/ BOU+ ( Tau2_v1 + Tau3_v1 )/ BOU  + 2 * Tau4_v1/ log(PW) )
B_error_integral_v2 =  RIF( Tau2_v2 * log(BOU)/ BOU + ( Tau2_v2 + Tau3_v2 )/ BOU + 2 * Tau4_v2/ log(PW) )

B_constant_s1_l = RIF(gamma - 6 / pi^2 * B_constant_integral1_u - B_error_integral_v1)
B_constant_s1_u = RIF(gamma - 6 / pi^2 * B_constant_integral1_l + B_error_integral_v1)

B_constant_s2_l = RIF(2 * (gamma + log(2) ) - 4 / pi^2 * B_constant_integral2_u - B_error_integral_v2)
B_constant_s2_u = RIF(2 * (gamma + log(2) ) - 4 / pi^2 * B_constant_integral2_l + B_error_integral_v2)

Psi2_v1 = Tau2_v1 + Tau3_v1 / log(BOU)
Psi2_v2 = Tau2_v2 + Tau3_v2 / log(BOU)

C_constant_integral1_l = -0.0495100113498
C_constant_integral1_u = -0.049510010626

C_constant_integral2_l = 2.63481269161
C_constant_integral2_u =  2.63481271383

C_error_integral_v1 =  RIF( Psi2_v1 * ( log(BOU)/ BOU - log(PW)/PW + 1/ BOU - 1/PW )  + 2 * Tau4_v1/ log(PW) )
C_error_integral_v2 =  RIF( Psi2_v2 * ( log(BOU)/ BOU - log(PW)/PW + 1/ BOU - 1/PW ) + 2 * Tau4_v2/ log(PW) )

C_constant_s1_l = RIF(gamma - 6 / pi^2 * C_constant_integral1_u - C_error_integral_v1)
C_constant_s1_u = RIF(gamma - 6 / pi^2 * C_constant_integral1_l + C_error_integral_v1)

C_constant_s2_l = RIF(2 * (gamma + log(2) ) - 4 / pi^2 * C_constant_integral2_u - C_error_integral_v2)
C_constant_s2_u = RIF(2 * (gamma + log(2) ) - 4 / pi^2 * C_constant_integral2_l + C_error_integral_v2)

sci1 = Upper(C_error_integral_v1,11) #we seek to avoid scientific notation
sci11 = '

sci2 = Upper(C_error_integral_v2,11)
sci22 = '

################################################################################# CALCULATIONS   ###

ds = 5
ds2 = 4

# The value of CCC lies in (0.002,733]

CCC_1 = 1/2
Brun_value_1 = -0.8264
Xiv1_value_1 = 12.48749

CCC_2 = 6
Brun_value_2 = -0.0069
Xiv1_value_2 = 1.42087

CCC_3 = 7
Brun_value_3 = 0.0036
Xiv1_value_3 = 1.25199

CCC_41 = 10
Brun_value_41 = 0.0201
Xiv1_value_41 = 0.94685

CCC_42 = 14
Brun_value_42 = 0.0267
Xiv1_value_5 = 0.745

CCC_4 = 15
Brun_value_4 = 0.027
Xiv1_value_4 = 0.71185

CCC_5 = 16
Brun_value_5 = 0.0271
Xiv1_value_5 = 0.68307

CCC_6 = 17
Brun_value_6 = 0.0269
Xiv1_value_6 = 0.65787

CCC_7 = 38
Brun_value_7 = 0.0004
Xiv1_value_7 = 0.45734

CCC_8 = 39
Brun_value_8 = -0.0013
Xiv1_value_8 = 0.45419

CCC_81 = 50
Brun_value_81 = -0.0199
Xiv1_value_81 = 0.43225

CCC_9 = 60
Brun_value_9 = -0.0372
Xiv1_value_9 = 0.42499

CCC_10 = 70
Brun_value_10 = -0.0544 
Xiv1_value_10 = 0.42425 

CCC_101 = 80
Brun_value_101 = -0.0715
Xiv1_value_101 = 0.42747

CCC_11 = 98
Brun_value_11 = -0.1018
Xiv1_value_11 = 0.43908

CCC_12 = floor( 10^(exp_prog_bound*(1-choice))/e^4)
Brun_value_12 = -0.3566
Xiv1_value_12 = 0

################################################################################# BRUN-TITCHMARSH   ###

ebrun = 2*exp_prog_bound

iot = 2*(1-4/pi^2)

br = numerical_v2

Brun = RIF( 4 * ( -C_constant_s2_l / 2 + br / 2 + 4 * ( 4 / pi^2 + iot / (10^ebrun)^(1/4) )^2 ) )

constantprime = RIF( ( Brun + log( 10^ebrun )^2 / sqrt(10^ebrun) ) / 2 )

\end{mysage}


\baselineskip=17pt


\title{On a logarithmic sum related to the Selberg sieve}

\author{Sebastian Zuniga Alterman\\
Institut de Math\'{e}matiques de Jussieu\\ 
Universit\'{e} Paris Diderot P7\\ 
B\^atiment Sophie Germain, 8 Place Aur\'elie Nemours \\ 
75013 Paris, France\\
E-mail: sebastian.zuniga-alterman@imj-prg.fr}

\date{}

\maketitle


\renewcommand{\thefootnote}{}



\renewcommand{\thefootnote}{\arabic{footnote}}
\setcounter{footnote}{0}


\begin{abstract}
We study the sum $\Sigma_q(U)=\sum_{\substack{d,e\leq U\\(de,q)=1}}\frac{\mu(d)\mu(e)}{[d,e]}\log\left(\frac{U}{d}\right)\log\left(\frac{U}{e}\right)$, $U>1$, so that a continuous, monotonic and explicit version of Selberg's sieve can be stated. 

Thanks to Barban--Vehov (1968), Motohashi (1974) and Graham (1978), it has been long known, but never explicitly, that $\Sigma_1(U)$ is asymptotic to $\log(U)$. In this article, we discover not
only that $\Sigma_q(U)\sim\frac{q}{\varphi(q)}\log(U)$ for all $q\in\mathbb{Z}_{>0}$, but also we find a closed-form expression for its secondary order term of $\Sigma_q(U)$, a constant $\mathfrak{s}_q$, which we are able to estimate explicitly when $q=v\in\{1,2\}$. We thus have $\Sigma_v(U)= \frac{v}{\varphi(v)}\log(U)-\mathfrak{s}_v+O_v^*\left(\frac{K_v}{\log(U)}\right)$,  for some explicit constant $K_v > 0$, where $\mathfrak{s}_1=0.60731\ldots$ and $\mathfrak{s}_2=1.4728\ldots$. 

As an application, we show how our result gives an explicit version of the Brun--Titchmarsh theorem within a range.
\end{abstract}
\normalsize

\section{Notation and basic definitions}

Throughout the present work the variable $p$ denotes a prime number, $q$ denotes an arbitrary positive integer and the function $X>0\mapsto\log^+(X)$ corresponds to $\max\{\log(X),0\}$. We also use the {\em $O^*$ notation}: we write $f(X)=O^*(h(X))$, as $X\to a$ to indicate that $|f(X)|\leq h(X)$ in a neighborhood of $a$, where, in absence of precision, $a$ corresponds to $\infty$.
Finally, we consider the {\em Euler $\varphi_s$ and Kappa $\kappa_s$} functions: let $s$ be any complex number, we define $\varphi_s:\mathbb{Z}_{>0}\to\mathbb{C}$ as $q\mapsto q^s\prod_{p|q}\left(1-\frac{1}{p^s}\right)$ and $\kappa_s:\mathbb{Z}_{>0}\to\mathbb{C}$ as $q\mapsto q^s\prod_{p|q}\left(1+\frac{1}{p^s}\right)$.

\section{Introduction}\label{Intro}

Let $U_1>1$ and $U_0>0$ such that $U_1>U_0$. Consider the Barban--Vehov weights $d\in\mathbb{Z}_{>0}\mapsto\mathbf{L}_d=\log^+\left(\frac{U_1}{d}\right)-\log^+\left(\frac{U_0}{d}\right)$; $d\mapsto\mathbf{L}_d$ is continuous on $(0,\infty)$ and satisfies $\frac{\mathbf{L}_1}{\log(U_1)}=1$ and $\frac{\mathbf{L}_d}{\log(U_1)}=0$ for all $d\geq U_1$. Therefore  $\left\{\frac{\mathbf{L}_d}{\log(U_1)}\right\}_{d=1}^\infty$ is a sequence of parameters as in Selberg's sieve that are also continuous and monotonic. In particular, if $U_1>2$, then the sum
\begin{align}\label{Selberg} 
\sum_{\substack{n\leq X\\(n,v)=1}}\left(\sum_{d|n}\frac{\mu(d)\mathbf{L}_d}{\log(U_1)}\right)^2
\end{align}
sifts the prime numbers in the interval $[U_1,X]$, where $v\in\{1,2\}$.   

Our main motivation is to give a continuous and monotonic version of the Selberg sieve, which, by stating it explicitly, has important consequences, as Theorem \ref{CON} below. Here, the condition $U_1^2\leq X$ is mandatory.
However, it is significant that Graham, in \cite{Gra78}, carries out a non-explicit asymptotic analysis for the sum \eqref{Selberg} that considers not only the case $U_1^2\leq X$ but also $X\leq U_1^2$ (see \cite[\S 4]{Gra78}) 

Sums like the one in \eqref{Selberg} have been studied non-explicitly, with $\frac{\mathbf{L}_d}{\log(U)}$ replaced by $d\mapsto\mathds{1}_{\{d\leq U\}}(d)$, by Dress, Iwaniec and Tenenbaum in \cite{DIT83} and recently by de la Bret\`eche, Dress and Tenenbaum in \cite{BDT19}. In this case, the analogous main term coefficient $\sum_{\substack{d,e\leq U}}\frac{\mu(d)\mu(e)}{[d,e]}$ converges to a positive constant whose rigorous estimation has been given by Helfgott in \cite[Prop. 6.30]{Hel19}, showing that $\sum_{\substack{d,e}}\frac{\mu(d)\mu(e)}{[d,e]}=0.440729+O^*(0.0000213).$

In this article, we study the asymptotic expression of \eqref{Selberg} with a particular choice of Barban--Vehov weights, namely, $U_1=U>1$ and $U_0=1$, so that $\mathbf{L}_d$ becomes the one parameter logarithmic weight $d\mapsto\log^+\left(\frac{U}{d}\right)$. With this choice, in our main result, we conclude in an explicit manner not only that the main term coefficient 
$\sum_{\substack{d,e\\(de,q)}}\frac{\mu(d)\mu(e)}{[d,e]}\mathbf{L}_d\mathbf{L}_e$ of \eqref{Selberg}  is asymptotic to $\frac{q}{\varphi(q)}\log(U)$ for all $q\in\mathbb{Z}_{>0}$, as shown in \cite{Gra78} without coprimality conditions, but we are also able to obtain its second order term, which is a constant value for all $q\in\mathbb{Z}_{>0}$. It is stated in \S\ref{Conclusion} and reads as follows.

\begin{ctheorem}{1}\label{MAI}
Let $U>1$. Then for all $q\in\mathbb{Z}_{>0}$, one can determine explicit constants $\mathfrak{s}_q$ and  $K_q>0$ such that
\begin{equation*}
\sum_{\substack{d,e\\(de,q)=1}}\frac{\mu(d)\mu(e)}{[d,e]}\mathbf{L}_d\mathbf{L}_e=\frac{q}{\varphi(q)}\log(U)-\mathfrak{s}_q+O_q^*\left(\frac{K_q}{\log(U)}\right).
\end{equation*}  
In particular,  
\begin{equation}\label{hope}
\mathfrak{s}_v=\begin{cases}
0.60731\ldots&\text{ if }v=1,\\
1.4728\ldots&\text{ if }v=2, 
\end{cases}
\end{equation}
and, if $U\geq 10^{\sage{exp_prog_bound}}$, we can select $K_1=\sage{Upper(KV1,digits+1)}$, $K_2=\sage{Upper(KV2,digits+1)}$.
\end{ctheorem}
The error term magnitude in the theorem above has been obtained by resting on modern estimations given by Balazard \cite{BA12}, Bord\`elles \cite{BO15}, El-Marraki \cite{EM95}, Helfgott \cite[\S 6]{Hel19} and Ramar\'e \cite{RA13}, \cite{RA15}, among others; nonetheless, when using non-explicit tools, one should expect an error of magnitude $e^{-w\sqrt{\log(U)}}$ for some constant $w>0$.  Furthermore, in order to derive Theorem \ref{MAI}, one must differ from a classical prime-number-theorem-like approach (as in \cite{Gra78}, for example), as it gives ineffective or well explicit estimations that involve huge numbers, thus too inconvenient or impractical to be used.  For this reason, in order to succeed, we introduce in \S\ref{bounds} some averaging functions involving the M\"obius function, from which we can derive \emph{explicit} estimations that are crucial for the results given in sections \S\ref{Logarithm} and \S\ref{Conclusion}.    

Our work is general enough to consider any coprimality condition, and not the specific ones $q=v\in\{1,2\}$ that we have chosen to work with to derive explicit constants. Hence, we are able to derive Theorem \ref{CONCLUSION} for any $q\in\mathbb{Z}_{>0}$, provided that the specific constant $K_q$ is found. In that case, given its closed-form and with the help of analogous results to propositions \ref{tail} and \ref{TT}, any constant $\mathfrak{s}_q$ with $q>2$ may be rigorously estimated. 

Subsequently, closely following \cite[\S 3.2]{MV07}, Theorem \ref{MAI} can be applied to derive an explicit version of the Brun--Titchmarsh theorem within a range, by providing a uniform upper bound for the number of primes in arithmetic progressions. That is, if we consider the quantity $\pi(X;q,a)=\{p\leq X, p\equiv a\ (mod\ q)\}$, we have the following.
\begin{ctheorem}{2}[\textbf{Brun--Titchmarsh inequality}]\label{CON} Let $a,q\in\mathbb{Z}_{>0}$ such that $(a,q)=1$. Let $Y$ be a real number such that $Y\geq 10^{\sage{ebrun}}q$. Then for all $X\geq 0$,
\begin{equation*}   
\pi(X+Y;q,a)-\pi(X;q,a)\leq\frac{2Y}{\varphi(q)\log\left(\frac{Y}{q}\right)}\left(1-\frac{\sage{-Upper(constantprime,digits+1)}}{\log\left(\frac{Y}{q}\right)}\right).
\end{equation*}  
\end{ctheorem}

\section{Bounds on functions involving the M\"obius function}\label{bounds}

In \S\ref{Logarithm}, we will need to reduce our estimations to some well-known and simpler functions. For any $X>0$, define $m_q(X)=\sum_{\substack{n\leq X\\(n,q)=1}}\frac{\mu(n)}{n}$ and consider
 \begin{align}
 \check{m}_q(X)=\sum_{\substack{n\leq X\\(n,q)=1}}\frac{\mu(n)}{n}\log\left(\frac{X}{n}\right),\quad\check{\check{m}}_q(X)&=\sum_{\substack{n\leq X\\(n,q)=1}}\frac{\mu(n)}{n}\log^2\left(\frac{X}{n}\right),\label{checks}\\
 \tilde{m}_q(X)=\sum_{\substack{n\leq X\\(n,q)=1}}\frac{\mu(n)}{\kappa(n)}\log\left(\frac{X}{n}\right),\quad\tilde{\tilde{m}}_q(X)&=\sum_{\substack{n\leq X\\(n,q)=1}}\frac{\mu(n)}{\kappa(n)}\log^2\left(\frac{X}{n}\right).\label{tildes}
\end{align}   
It is straightforward to see by summation by parts that $\check{m}_q$ and $\check{\check{m}}_q$, as well as $\tilde{m}_q$ and $\tilde{\tilde{m}}_q$, are related by the following identity.
 \begin{lemma}\label{integral}
 Let $X\geq 1$. We have 
\begin{equation*}
 \int_1^X\check{m}_q(s)\frac{ds}{s}=\frac{1}{2}\check{\check{m}}_q(X),\quad\int_1^X\tilde{m}_q(s)\frac{ds}{s}=\frac{1}{2}\tilde{\tilde{m}}_q(X).
\end{equation*}
 \end{lemma}

\noindent Moreover, we have the following explicit estimations
\begin{align}\label{checkbounds}  
|\check{m}(X)-1|&\leq\frac{1}{\sqrt{X}},&\text{ if }& 0< X\leq 10^{12}&\text{ \cite[Lemma 5.9]{Hel19}},&&\nonumber\\  
&\leq\frac{1}{389\log(X)},&\text{ if }& X\geq 3155&\text{ \cite[Thm. 1.5]{RA15}},&&\\
|\check{\check{m}}(X)-2\log(X)+2\gamma|&\leq\frac{4e^{\frac{\gamma}{2}-1}}{\sqrt{X}},&\text{ if }& 0< X\leq 10^{12}&\text{ \cite[Lemma 5.9]{Hel19}},&&\nonumber\\ 
&\leq\frac{1}{103\log(X)},&\text{ if }& X\geq 9&\text{ \cite[Thm. 1.8]{RA15}},&&\label{checkcheckbounds}
\end{align}
where the first bound in each case has been obtained with the help of computer calculations using an implementation of interval arithmetic.   

Consider now $q,d\in\mathbb{Z}_{>0}$. We write $d|q^\infty$, meaning that $d$ is in the set $\{d',\ p|d'\implies p|q\}$. If $n\in\mathbb{Z}_{>0}$, then $(q^\infty,n)$ is the greatest divisor $d'$ of $n$ such that $d'|q^\infty$. Therefore $(q^\infty,n)=1$ if and only if $(n,q)=1$; otherwise, $(q^\infty,n)$ and $(n,q)$ may differ. 
With this definition, and by using the following identity, established for example throughout \cite[Lemma 2]{Gra78} and in \cite[Eq. (5.72)]{Hel19}, one can study a sum with coprimality conditions from the same sum without such conditions.  

 \begin{lemma}\label{Idd}
 We have the identity $ \sum_{d|q^{\infty},d|n}\mu\left(\frac{n}{d}\right)=\mu(n)\mathds{1}_{\{(n,q)=1\}}(n)$.
Hence, for any function $h:\mathbb{Z}_{>0}\to\mathbb{C}$, we have the formal identity
 \begin{equation}\label{sumi}
 \sum_{\substack{n\\(n,q)=1}}\frac{\mu(n)}{n}h(n)=\sum_{d|q^\infty}\frac{1}{d}\sum_{n}\frac{\mu(n)}{n}h(dn).
 \end{equation}
 \end{lemma} 

 Lemma \ref{Idd} is meaningful since, as Ramar\'e points out in \cite[\S 1]{RA13}, on using merely a M\"obius inversion $\sum_{d|q,d|n}\mu\left(d\right)=\mathds{1}_{\{(n,q)=1\}}(n)$ in \eqref{sumi}, one would have been taken back to a sum having, again, coprimality conditions. For example, with the help of Lemma \ref{Idd}, we have

 \begin{lemma}\label{check} 
 Let $X>0$ and $\theta=1-\frac{1}{\log(10^{12})}$. Then 
 \begin{equation*} 
\left|\check{m}_q(X)-\frac{q}{\varphi(q)}\right|\leq\frac{\sqrt{q}}{\varphi_{\frac{1}{2}}(q)}\frac{1}{\sqrt{X}}+\frac{q^\theta}{\varphi_{\theta}(q)}\frac{\mathds{1}_{\{X\geq 10^{12}\}}(X)}{389\log(X)}.
 \end{equation*}
 \end{lemma}

The proof of the Lemma \ref{check} is given by Helfgott in \cite[Prop 5.15]{Hel19}; although it is given in the range $X\geq 1$, it is not difficult to derive it for any $X>0$ as long as the term of order $\frac{1}{\log(X)}$ is considered for sufficiently large values of $X$. This proof consists on finding a way to put together the bounds \eqref{checkbounds} so that resulting estimation comes from a direct application of identity \eqref{sumi}; it is indeed very convenient to have general inequalities that prevent us from splitting a summation into many ranges that are not in general simple to handle on their own. Inspired by this remark, we derive our first result, that will help us to further understand the estimations given in lemmas \ref{SI:estimation} and \ref{S1:total}.

 \begin{lemma}\label{checkk}  
 Let $X>0$ and $\theta=1-\frac{1}{\log(10^{12})}$. Define $f_q:X\geq 1\mapsto \log(X)-\gamma-\sum_{p|q}\frac{\log(p)}{p-1}$. Then
 \begin{align*}
 \left|\check{\check{m}}_q(X)-\frac{2q}{\varphi(q)}f_q(X)\right|\leq\frac{\sqrt{q}}{\varphi_{\frac{1}{2}}(q)}\frac{4e^{\frac{\gamma}{2}-1}}{\sqrt{X}}+\frac{q^\theta}{\varphi_{\theta} (q)}\frac{\mathds{1}_{\{X\geq 10^{12}\}}(X)}{103\log(X)}. 
 \end{align*} 
 \end{lemma}

  \begin{proof} 
  By taking $h(n)=(\max\left\{\log\left(\frac{X}{n}\right),0\right\})^2$ in Lemma \ref{Idd}, we derive
  \begin{equation}\label{checkcheck} 
  \check{\check{m}}_q(X)=\sum_{d|q^\infty}\frac{1}{d}\sum_{n\leq\frac{X}{d}}\frac{\mu(n)}{n}\log^2\left(\frac{X}  {dn}\right)=\sum_{d|q^\infty}\frac{1}{d}\check{\check{m}}\left(\frac{X}{d}\right).
  \end{equation}
  Consider the Dirichlet series $\sum_{d|q^\infty}\frac{1}{d^s}$; it converges to $\frac{q^s}{\varphi_s(q)}$ for all $s\in\mathbb{C}$ such that $\Re(s)>0$. Subsequently, we can differentiate it to obtain
  \begin{equation}\label{derivate}
  -\sum_{d|q^\infty}\frac{\log(d)}{d}=\left(\sum_{d|q^\infty}\frac{1}{d^s}\right)'_{s=1}=\left(\frac{q^s}{\varphi_s(q)}\right)'_{s=1}=\left(\frac{q^s}{\varphi_s(q)}\sum_{p|q}\frac{-\log(p)}  {p^s-1}\right)_{s=1}.
  \end{equation}
  On the other hand, by combining the bounds given in \eqref{checkcheckbounds}, we have that
  \begin{align}\label{error1} 
  \left|\sum_{d|q^\infty}\frac{1}{d}\left(\check{\check{m}}\left(\frac{X}{d}\right)-2\log\left(\frac{X}{d}\right)+2\gamma\right)\right|\leq\frac{1}{\sqrt{X}}\sum_{d|q^\infty}\frac{4e^{\frac{\gamma}{2}-1}}{\sqrt{d}}+\sum_{d|  q^\infty}\frac{\mathds{1}_{\{\frac{X}{d}\geq 10^{12}\}}(d)}{103\ d\log\left(\frac{X}{d}\right)}\phantom{x}&\\ 
  \leq\frac{\sqrt{q}}{\varphi_{\frac{1}{2}}(q)}\frac{4e^{\frac{\gamma}{2}-1}}{\sqrt{X}}+\frac{\mathds{1}_{\{X\geq 10^{12}\}}(X)}{103\log(X)}\sum_{d|q^\infty}\frac{1}{d^\theta}=\frac{\sqrt{q}}{\varphi_{\frac{1}{2}}(q)}\frac{4e^{\frac{\gamma}{2}-1}}{\sqrt{X}}+\frac{q^\theta}{\varphi_{\theta}(q)}\frac{\mathds{1}_{\{X\geq 10^{12}\}}(X)}{103\log(X)},\phantom{x}&\nonumber 
  \end{align}
  where, recalling the definition of $\theta$, we have used that the function $d\mapsto\frac{1}{d^{1-\theta}\log\left(\frac{X}{d}\right)}$ is decreasing for $1\leq d\leq\frac{X}{10^{12}}$. 
  We conclude the result  from \eqref{derivate}, by identifying the following identity
  \begin{align}\label{cola}
  \sum_{d|q^\infty}\frac{1}{d}\left(\log\left(\frac{X}{d}\right)-\gamma\right)&=\frac{q}{\varphi(q)}\left(\log(X)-\gamma-\sum_{p|q}\frac{\log(p)}{p-1}\right),
  \end{align}
 that can be replaced in the leftmost expression of inequality \eqref{error1}.
 \end{proof}
 
Notice that \eqref{error1}  is valid since  the bounds \eqref{checkcheckbounds} hold regardless of whether or not $\frac{X}{d}\geq 1$; this is the reason why we do not incorporate a range condition on the variable $d$ in the outer sum of the expression \eqref{checkcheck} and also why the values of $X$ and $d$ remain independent.

We provide now the main term of the function $X\mapsto\tilde{\tilde{m}}_q(X)$. It might come as a surprise that $\mathfrak{a}_q$, defined below, is closely related to the function that optimizes the Selberg sieve (refer to \cite[\S 3.2]{MV07}), introduced in Lemma \ref{sum1}. 

 \begin{lemma}\label{tildetildem} 
 Let $X>0$ and $\theta=1-\frac{1}{\log(10^{12})}$. Then 
  \begin{equation*}
 \left|\tilde{\tilde{m}}_q(X)-\frac{2\zeta(2)\kappa(q)}{q}\left(\log(X)-\mathfrak{a}_q\right)\right|\leq\mathpzc{p}_{\frac{1}{2}}(q)\ \frac{4e^{\frac{\gamma}{2}-1}\mathit{P}_{\frac{1}{2}}}{\sqrt{X}}+\mathpzc{p}_{\theta}(q)\ \frac{\mathit{P}_{\theta}\mathds{1}_{\{X\geq 10^{12}\}}(X)}{103\log(X)},
  \end{equation*}
 \begin{flalign*} 
\text{where}&& \mathfrak{a}_q=\sum_{p}\frac{\log(p)}{p(p-1)}+\gamma+\sum_{p|q}\frac{\log(p)}{p},\quad \sum_{p}\frac{\log(p)}{p(p-1)}+\gamma=   \sage{C_sum1}    \ldots,\phantom{xxxxxx}&\\
 &&\mathpzc{p}_{\alpha}(q)=\prod_{p|q}\left(1+\frac{p^{1-\alpha}}{p+1-p^{1-\alpha}}\right),\quad\mathit{P}_{\frac{1}{2}}=   \sage{Upper(Delta_1half,digits)}   ,\quad\mathit{P}_{\theta}=   \sage{Upper(Delta_theta,digits)}.\phantom{xxxxxxxxxx} & 
 \end{flalign*}   
 \end{lemma}

  \begin{proof}
  Observe that for any square-free $n$, $\frac{n}{\kappa(n)}=\sum_{d|n}\frac{\mu(d)}{\kappa(d)}$. Therefore 
  \begin{align}
  \tilde{\tilde{m}}_q(X)=\sum_{\substack{n\leq X\\(n,q)=1}}\frac{\mu(n)}{n}\log^2\left(\frac{X}{n}\right)\sum_{d|n}\frac{\mu(d)}{\kappa(d)}&=\sum_{\substack{d\\(d,q)=1}}\frac{\mu^2(d)}{d\kappa(d)}\check{\check{m}}_{dq}\left(\frac{X}{d}\right).\label{analogous}
  \end{align}
  By \eqref{analogous}, we are now able to derive the main term of $\tilde{\tilde{m}}_q\left(\frac{X}{d}\right)$ in a similar manner to the obtainment of expression \eqref{error1}. Indeed, by using that $X$ and $d$ are independent variables, and with the help of Lemma \ref{checkk}, we have 
  \begin{align*} 
  &\sum_{\substack{d\\(d,q)=1}}\frac{\mu^2(d)}{d\kappa(d)}\left|\check{\check{m}}_{dq}\left(\frac{X}{d}\right)-\frac{2dq}{\varphi(dq)}\left(\log\left(\frac{X}{d}\right)-\gamma-\sum_{p|dq}\frac{\log(p)}{p-1}\right)\right|\\ 
 \leq&\ \frac{\sqrt{q}}{\varphi_{\frac{1}{2}}(q)}\frac{4e^{\frac{\gamma}{2}-1}}{\sqrt{X}}\sum_{\substack{d\\(d,q)=1}}\frac{\mu^2(d)}{\varphi_{\frac{1}{2}}(d)\kappa(d)}\\ 
&\phantom{xxxxxxxxxxxxxxxxx}+\frac{q^\theta}{\varphi_{\theta}(q)}   \frac{\mathds{1}_{\{X\geq 10^{12}\}}(X)}{103}\sum_{\substack{d\leq\frac{X}{10^{12}}\\(d,q)=1}} \frac{\mu^2(d)}{d^{1-\theta}\varphi_\theta(d)\kappa(d)\log\left(\frac{X}{d}\right)}\\
  \leq&\ \mathpzc{p}_{\frac{1}{2}}(q)\ \frac{4e^{\frac{\gamma}{2}-1} \mathit{P}_{\frac{1}{2}}}{\sqrt{X}}+\mathpzc{p}_{\theta}(q)\ \frac{\mathit{P}_{\theta}\mathds{1}_{\{X\geq 10^{12}\}}(X)}{103\log(X)},
  \end{align*} 
  where we have used that $d\mapsto\frac{1}{d^{1-\theta}\log\left(\frac{X}{d}\right)}$ is decreasing for $1\leq d\leq\frac{X}{10^{12}}$; thereupon, we have completed the above summations on the variable $d$ to derive a convergent sum, obtaining, for $\alpha\in\left\{\frac{1}{2},\theta\right\}$,
  \begin{align*}
  &\mathpzc{p}_{\alpha}(q)=\frac{q^\alpha}{\varphi_{\alpha}(q)}\prod_{p|q}\left(1+\frac{1}{(p^\alpha-1)(p+1)}\right)^{-1}=\prod_{p|q}\frac{p+1}{p+1-p^{1-\alpha}},\\
  &\prod_{p}\left(1+\frac{1}{(p^\alpha-1)(p+1)}\right)\in
\begin{cases}
[\sage{Lower(Delta_1half_l,digits)},\mathit{P}_{\alpha}],\quad\text{if }\alpha=\frac{1}{2},&\\ 
[\sage{Lower(Delta_theta_l,digits)},\mathit{P}_{\alpha}],\quad\text{if }\alpha=\theta.& 
\end{cases}
  \end{align*}
  On the other hand, let $F_q(s)=\sum_{\substack{d\\(d,q)=1}}\frac{\mu^2(d)}{d^s\varphi(d)\kappa(d)}$; it is a well-defined function for all $s\in\mathbb{C}$ such that $\Re(s)>-1$ and we have that
  \begin{align}\label{it}
  \sum_{\substack{d\\(d,q)=1}}\frac{\mu^2(d)}{\varphi(d)\kappa(d)}\sum_{p|d}\frac{\log(p)}{p-1}=\sum_{\substack{p\nmid q}}\frac{\log(p)}{(p-1)(p^2-1)}\sum_{\substack{e\\(e,pq)=1}}\frac{\mu^2(e)}{\varphi(e)\kappa(e)}\phantom{xxxxxxxx}&\nonumber\\
  =\sum_{\substack{p\nmid q}}\frac{\log(p)}{(p-1)(p^2-1)}\prod_{p'\nmid pq}\left(1+\frac{1}{p'^2-1}\right)=F_q(0)\sum_{\substack{p\nmid q}}\frac{\log(p)}{p^2(p-1)},&
  \end{align}
  where $F_q(0)=\zeta(2)\frac{\varphi(q)\kappa(q)}{q^2}$.
  Therefore, from \eqref{it}, we derive
  \begin{align*}
  &\frac{2q}{\varphi(q)}\sum_{\substack{d\\(d,q)=1}}\frac{\mu^2(d)}{\varphi(d)\kappa(d)}\left(\log\left(\frac{X}{d}\right)-\gamma-\sum_{p|dq}\frac{\log(p)}{p-1}\right)\\
  =&\ \frac{2q}{\varphi(q)}F_q(0)\left(\log(X)+\frac{F_q'(0)}{F_q(0)}-\gamma-\sum_{p|q}\frac{\log(p)}{p-1}-\sum_{\substack{p\nmid q}}\frac{\log(p)}{p^2(p-1)}\right)\\
  =&\ 2\zeta(2)\frac{\kappa(q)}{q}\left(\log(X)-\sum_{p}\frac{\log(p)}{p(p-1)}-\gamma-\sum_{p|q}\frac{\log(p)}{p}\right),
  \end{align*}
  and we conclude the result by observing that
  \begin{align}\label{Summation0}
  \sum_{p}\frac{\log(p)}{p(p-1)}+\gamma\in[\sage{Lower(I_sum1_l+gamma,14)},\sage{Upper(I_sum1_u+gamma,14)} ].
  \end{align} 
  \end{proof}

 A similar treatment to \eqref{analogous} allows us to derive the main term of $\tilde{m}_q(X)$ and thus Lemma \ref{tildem}. We do not bother to perform the involved calculations, since it is a result that has already been proved, by different means, in \cite[Prop. 6.8]{Hel19}; it reads as follows.

 \begin{lemma}\label{tildem} 
 Let $X>0$ and $\theta=1-\frac{1}{\log(10^{12})}$. Then 
 \begin{equation*}
 \left|\tilde{m}_q(X)-\frac{\zeta(2)\kappa(q)}{q}\right|\leq\mathpzc{p}_{\frac{1}{2}}(q)\ \frac{\mathit{P}_{\frac{1}{2}}}{\sqrt{X}}+\mathpzc{p}_{\theta}(q)\ \frac{\mathit{P}_{\theta}\mathds{1}_{\{X\geq 10^{12}\}}(X)} {389\log(X)},
  \end{equation*}
 where $\mathit{P}_{\frac{1}{2}}$, $\mathit{P}_{\theta}$, $\mathpzc{p}_{\frac{1}{2}}(q)$ and $\mathpzc{p}_{\theta}(q)$ are defined in Lemma \ref{tildetildem}.
 \end{lemma}

Consider now $h_q:s\geq 1\mapsto\sum_{\substack{d\\(d,q)=1}}\frac{\mu(d)}{\kappa(d)^2}\left(\tilde{m}_{dq}\left(\frac{s}{d}\right)-\frac{\pi^2}{6}\frac{\kappa(dq)}{dq}\right)^2$. By expanding the square, we readily have (see \cite[Lemma 6.10]{Hel19})
\begin{align}\label{h_q}  
\sum_{\substack{d\\(d,q)=1}}\frac{\mu(d)}{\kappa(d)^2}\ \tilde{m}_{dq}^2\left(\frac{s}{d}\right)=h_q(s)+\frac{\pi^2\kappa(q)}{3q}\check{m}_q(s)-\frac{\pi^2\kappa(q)}{6\varphi(q)},
\end{align}
so that, as the above left-hand side is a finite sum, the function $h_q$ is well-defined.

In order to estimate $h_q$, we will need the following result. It is a general tool that helps in deriving the correct order of arithmetic averages that are weighted by suitable negative powers of logarithms. It will also be useful to analyze propositions \ref{Ss1Log} and \ref{Ss2Log}.  

 \begin{proposition}\label{logint} Let $Z,X,m,n$ be real numbers such that $m\geq 1$ and $1\leq Z<X$. Then, $\mathbf{a)}$ if $m=1$,
  \begin{align*} 
  \int_{1}^{Z}\frac{du}{u^m\log^n\left(\frac{X}{u}\right)}=
\begin{cases}
\log\left(\frac{\log(X)}{\log\left(\frac{X}{Z}\right)}\right)\quad&\text{ if }n=1,\\ 
 \frac{1}{n-1}\left(\frac{1}{\log^{n-1}\left(\frac{X}{Z}\right)}-\frac{1}{\log^{n-1}(X)}\right)\quad&\text{ if }n\neq 1,
\end{cases} 
 \end{align*}
 $\mathbf{b)}$ if $m>1$ and $n>0$, then
 \begin{align*}
 \int_{1}^{Z}\frac{du}{u^m\log^n\left(\frac{X}{u}\right)}\leq\frac{1}{m-1}\left(\frac{1}{\log^n\left(\frac{X}{\sqrt{Z}}\right)}+\frac{1}{\log^n\left(\frac{X}{Z}\right)\sqrt{Z^{m-1}}}\right).
 \end{align*}
 \end{proposition} 
  \begin{proof}
  If $m=n=1$, we have $\int_{1}^{Z}\frac{du}{u^m\log^n\left(\frac{X}{u}\right)}=\left[-\log\left(\log\left(\frac{X}{u}\right)\right)\right|_1^{Z}$; if $m=1$, and $n\neq 1$, we have $\int_{1}^{Z}\frac{du}{u\log^n\left(\frac{X}{u}\right)}=\left[\frac{1}{n-1}\log^{-(n-1)}\left(\frac{X}{u}\right)\right|_1^{Z}$, whence $\mathbf{a)}$.

  With respect to $\mathbf{b)}$, if $n>0$, the function $u\mapsto\log^{-n}\left(\frac{X}{u}\right)$ is increasing for $1\leq u<X$ and since $m>1$, for any $0<k<1$, in particular for $k=\sage{kk}$, we conclude that
  \begin{align*}
  \int_{1}^{Z}\frac{du}{u^m\log^n\left(\frac{X}{u}\right)}&\leq\frac{1}{\log^n\left(\frac{X}{Z^k}\right)}\int_{1}^{Z^{k}}\frac{du}{u^m}+\frac{1}{\log^n\left(\frac{X}{Z}\right)}\int_{Z^{k}}^{Z}\frac{du}{u^m}\\
  &\leq\frac{1}{m-1}\left(\frac{1}{\log^n\left(\frac{X}{Z^{k}}\right)}+\frac{1}{\log^n\left(\frac{X}{Z}\right)Z^{k(m-1)}}\right).
  \end{align*} 
 Whence the result. 
  \end{proof}

\begin{proposition}\label{tail} Let $q\in\mathbb{Z}_{>0}$. The integral $\int_1^{\infty}\frac{h_q(s)}{s}ds$ converges and defines a constant depending on $q$.
Moreover, for any $X>0$, we have the following tail order estimation
\begin{equation*}  
\int_X^{\infty}\frac{h_q(s)}{s}ds=O_q\left(\frac{1}{\log(X)}\right).
\end{equation*} 
\end{proposition}    

\begin{proof} Given \cite[Thm. 3.3]{1ZA20}, we can derive a more theoretical proof. By Lemma \ref{tildem}, for any $d\in\mathbb{Z}_{>0}$ with $(d,q)=1$, the main term of $\tilde{m}_{dq}$  is $\frac{\pi^2}{6}\frac{\kappa(dq)}{dq}$ and we have the bound
\begin{align*}
|h_q(s)|\leq \sum_{\substack{d\leq s\\(d,q)=1}}\frac{\mu^2(d)}{\kappa(d)^2}\left(A_q^d(s)^2+2A_q^d(s)B_q^d(s)+B_q^d(s)^2\right),
\end{align*}
\begin{flalign*} 
\text{where}&&A_q^d(s)= \mathit{P}_{\frac{1}{2}}\mathpzc{p}_{\frac{1}{2}}(dq)\frac{\sqrt{d}}{\sqrt{s}},\quad B_q^d(s)=\mathit{P}_{\theta}\mathpzc{p}_{\theta}(dq)\frac{\mathds{1}_{\{\frac{s}{d}\geq 10^{12}\}}(d)} {389\log\left(\frac{s}{d}\right)}.\phantom{xxxxxxx}
\end{flalign*} 
Observe that 
\begin{align}
\sum_{\substack{d\leq s\\(d,q)=1}}\frac{\mu^2(d)}{\kappa(d)^2}A_q^d(s)^2&=\frac{\mathit{P}^{2}_{\frac{1}{2}}\mathpzc{p}_{\frac{1}{2}}(q)^2}{s}\sum_{\substack{d\leq s\\(d,q)=1}}\frac{\mu^2(d)\mathpzc{p}_{\frac{1}{2}}(d)^2d}{\kappa(d)^2}\leq\frac{a_q^{(1)}\log(s)+a_q^{(2)}}{s}\label{intuno},
\end{align}
and, by using Proposition \ref{logint}, 
\begin{align}
\sum_{\substack{d\leq s\\(d,q)=1}}\frac{\mu^2(d)}{\kappa(d)^2}B_q^d(s)^2&=\frac{\mathit{P}^{2}_{\theta}\mathpzc{p}_{\theta}(q)^2}{389^2}\sum_{\substack{d\leq\frac{s}{10^{12}}\\(d,q)=1}}\frac{\mu^2(d)\mathpzc{p}_{\theta}(d)^2}{\kappa(d)^2\log^2\left(\frac{s}{d}\right)}\leq\frac{b_q}{\log^2(s)}\label{integrable},
\end{align}
for some positive values $a_q^{(1)}$, $a_q^{(2)}$ and $b_{q}$ depending solely on $q$. From both estimations above, the sum $\sum_{\substack{d\\(d,q)=1}}\frac{\mu^2(d)}{\kappa(d)^2}A_q^d(s)B_q^d(s)$ can be bounded by Cauchy--Schwarz inequality, giving a term of order $O_q\left(\frac{1}{\sqrt{s}\log(s)}\right)$. 

Finally, as $\int\frac{\log(s)+1}{s^2}ds=-\frac{\log(s)+2}{s}$ and $\int\frac{ds}{s\log^2(s)}=\frac{-1}{\log(s)}$, we derive from \eqref{intuno}, \eqref{integrable} and Cauchy-Schwarz inequality for integrals that the integral $\int_1^{\infty}\frac{h_q(s)}{s}ds$ converges, and further that $\int_X^{\infty}\frac{h_q(s)}{s}ds=O_q\left(\frac{1}{\log(X)}\right)$.
\end{proof}

In \S\ref{A}, we will need explicit estimations for $h_v$, $v\in\{1,2\}$. As predicted by the estimations \eqref{intuno}, \eqref{integrable} and thanks to \cite[Prop. 6.14]{Hel19} and \cite[Prop. 6.17]{Hel19}, we have

\begin{proposition}\label{TT} For any $s\geq 1$, $|h_v(s)|$, $v\in\{1,2\}$, is at most
\begin{align*} 
\mathbf{i)}\ \frac{T_v^{(2)}\log(s)+T_v^{(3)}}{s},\quad\text{ if }1\leq s\leq 10^{12},\text{ or }\quad\mathbf{ii)}\ \frac{T_v^{(4)}}{\log^2(s)},\quad\text{ if }s\geq 10^{12},  
\end{align*}
where
\begin{align*}
T_1^{(2)} &= \sage{Trunc(Tau2_v1,dlong)},\quad   &   T_1^{(3)} &= \sage{Trunc(Tau3_v1,dlong)},\quad   &   T_1^{(4)}&= 0.000033536,\\
T_2^{(2)}&= \sage{Trunc(Tau2_v2,dlong)},\quad   &   T_2^{(3)} &=\sage{Trunc(Tau3_v2,dlong)},\quad   &    T_2^{(4)} &= 0.0000615022.
\end{align*} 
\end{proposition}

\section{A logarithmic sum involving the M\"obius function}\label{Logarithm} 

In order to start our analysis, let $U\geq 10^7$. Consider a parameter $1<Z<U$ such that $\frac{U}{Z}\geq 20$ and $Z\geq 4\times 10^5$, and write
\begin{align}\label{split} 
\sum_{\substack{d,e\\(de,v)=1}}\frac{\mu(d)\mu(e)}{[d,e]}\mathbf{L}_d\mathbf{L}_e=\sum_{\substack{\ell\leq U\\(\ell,v)=1}}\frac{\mu^2(\ell)}{\ell}\sum_{\substack{r_1,r_2\\ (r_1,r_2)=1\\(r_1r_2,\ell v)=1}}\frac{\mu(r_1)\mu(r_2)}{r_1r_2}\mathbf{L}_{\ell r_1}\mathbf{L}_{\ell r_2}=\mathit{S}_{\mathbf{I}}+\mathit{S}_{\mathbf{II}},
\end{align}
where 
\begin{align} 
\mathit{S}_{\mathbf{I}}=\mathit{S}_{\mathbf{I}}(U,Z,v)&=\sum_{\substack{Z<\ell\leq U\\(\ell,v)=1}}\frac{\mu^2(\ell)}{\ell}\sum_{\substack{r_1,r_2\\ (r_1,r_2)=1\\(r_1r_2,\ell v)=1}}\frac{\mu(r_1)\mu(r_2)}{r_1r_2}\mathbf{L}_{\ell r_1}\mathbf{L}_{\ell r_2},\label{def_S_I}\\
\mathit{S}_{\mathbf{II}}=\mathit{S}_{\mathbf{II}}(U,Z,v)&=\sum_{\substack{\ell\leq Z\\(\ell,v)=1}}\frac{\mu^2(\ell)}{\ell}\sum_{\substack{r_1,r_2\\ (r_1,r_2)=1\\(r_1r_2,\ell v)=1}}\frac{\mu(r_1)\mu(r_2)}{r_1r_2}\mathbf{L}_{\ell r_1}\mathbf{L}_{\ell r_2}\label{def_S_II}.
\end{align} 
The reason why the above sums have been introduced is related to the obtainment of actual error terms, which otherwise fail to arise. Indeed, in order to deal with lower order terms, two different approaches are required; one for $\mathit{S}_{\mathbf{I}}^{(v)}$ and another for $\mathit{S}_{\mathbf{II}}^{(v)}$, neither of them being satisfactory when applied to both sums at once.

We first show an estimation for $\mathit{S}_{\mathbf{I}}$, where the convergence of the integral below is assured by Proposition \ref{tail}. Its proof is given in \S\ref{A}.
\begin{lemma}\label{SI:estimation}         
Let $U\geq 10^{\sage{exp_prog}}$ and $v\in\{1,2\}$. If $Z$ is a real number such that $\frac{U}{Z}\geq 20$, then
\begin{align*}
\mathit{S}_{\mathbf{I}}=&\frac{6v}{\pi^2\kappa(v)}\int_1^{\infty}\frac{h_v(s)}{s}ds+\check{\check{m}}_v\left(\frac{U}{Z}\right)-\frac{v}{\varphi(v)}\log\left(\frac{U}{Z}\right)\nonumber\\
&+O^*\left(\frac{T_v^{(1)}\log^4\left(\frac{U}{Z}\right)}{\sqrt{Z}}+\frac{6v}{\pi^2\kappa(v)}\left(\Psi_v\left(\frac{Z\log\left(\frac{U}{Z}\right)}{U}+\frac{Z}{U}\right)+\frac{2T_v^{(4)}}{\log\left(\frac{U}{Z}\right)}\right)\right),
\end{align*} 
\begin{flalign*} 
\text{where}&&\phantom{x..}
T^{(1)}_v= 
\begin{cases} 
\sage{Upper(Tau1_v1,digits)}&\text{ if }v=1,\\
\sage{Upper(Tau1_v2,digits)}&\text{ if }v=2,
\end{cases}\qquad
\Psi_v=\left(T_v^{(2)}+\frac{T_v^{(3)}}{\log(20)}\right),\quad\text{ and }\phantom{xxxxxxxx}
\end{flalign*} 
$T_v^{(2)}$, $T_v^{(3)}$ and $T_v^{(4)}$ are defined in Proposition \ref{TT}.
\end{lemma} 

On the other hand, by M\"obius inversion, and recalling \eqref{checks}, we can write   
\begin{align}\label{SI}    
\mathit{S}_{\mathbf{II}}&=\sum_{\substack{\ell\leq Z\\(\ell,v)=1}}\frac{\mu^2(\ell)}{\ell}\sum_{\substack{d\leq\frac{U}{\ell}\\(d,\ell v)=1}}\frac{\mu(d)}{d^2}\ \check{m}_{\ell dv}\left(\frac{U}{\ell d}\right)^2.
\end{align}  
As the main term of $X\mapsto\check{m}_{\ell dv}(X)$ is $\frac{\ell d v}{\varphi(\ell v d)}$, we use Lemma \ref{checkk} to derive three more summations. Namely, we write $
\mathit{S}_{\mathbf{II}}=2\mathit{S}_{\mathbf{II}}^{(1)}-\mathit{S}_{\mathbf{II}}^{(2)}+\mathit{S}_{\mathbf{II}}^{(3)}$, where
\begin{align} 
\mathit{S}_{\mathbf{II}}^{(1)}=&\frac{v}{\varphi(v)}\sum_{\substack{\ell\leq Z\\(\ell,v)=1}}\frac{\mu^2(\ell)}{\varphi(\ell)}\sum_{\substack{d\leq\frac{U}{\ell}\\(d,\ell v)=1}}\frac{\mu(d)}{d\varphi(d)}\check{m}_{\ell dv}\left(\frac{U}{\ell d}\right).\label{S1:eq} \\
\mathit{S}_{\mathbf{II}}^{(2)}=&\sum_{\substack{\ell\leq Z\\(\ell,v)=1}}\frac{\mu^2(\ell)}{\ell}\sum_{\substack{d\leq\frac{U}{\ell}\\(d,\ell v)=1}}\frac{\mu(d)}{d^2}\frac{(\ell dv)^2}{\varphi(\ell d v)^2},\label{S2:eq}\\
\mathit{S}_{\mathbf{II}}^{(3)}=&\sum_{\substack{\ell\leq Z\\(\ell,v)=1}}\frac{\mu^2(\ell)}{\ell}\sum_{\substack{d\leq\frac{U}{\ell}\\(d,\ell v)=1}}\frac{\mu(d)}{d^2}\left(\check{m}_{\ell dv}\left(\frac{U}{\ell d}\right)-\frac{\ell d v}{\varphi(\ell v d)}\right)^2\label{S3:eq}.
\end{align} 
Each one of the above three summations will be estimated separately by the following lemmas that will be proven in \S\ref{S21}, \S\ref{S22} and \S\ref{S23}, respectively. 

\begin{lemma}\label{S1:total}     
Let $U\geq 10^{\sage{exp_prog}}$, $0<Z<U$ such that $\frac{U}{Z}\geq 20$ and $v\in\{1,2\}$. We have that
\begin{align*} 
\mathit{S}_{\mathbf{II}}^{(1)}=\frac{v}{\varphi(v)}\log(U)-\frac{1}{2}\ \check{\check{m}}_v\left(\frac{U}{Z}\right)+O^*\left(\frac{\Upsilon^{(1)}_v\log^2\left(\frac{U}{Z}\right)}{\sqrt{Z}}\right),
\end{align*}  
\begin{flalign*}
&\text{where }&&\Upsilon^{(1)}_v=
\begin{cases}
\sage{Upper(Upsilon1_v1,digits)} &\text{ if }v=1,\\ 
\sage{Upper(Upsilon1_v2,digits)} &\text{ if }v=2.
\end{cases}\phantom{xxxxxxxxxxxxxxxxxxxx}
\end{flalign*}
\end{lemma}
Additionally, by Lemma \ref{checkk}, we can replace $\check{\check{m}}_v\left(\frac{U}{Z}\right)$ in lemmas \ref{SI:estimation} and \ref{S1:total} by $\frac{2v}{\varphi(v)}\left(\log\left(\frac{U}{Z}\right)-\gamma-\sum_{p|v}\frac{\log(p)}{p-1}\right)$ $+O^*\left(\frac{\sqrt{q}}{\varphi_{\frac{1}{2}}(q)}\frac{4e^{\frac{\gamma}{2}-1}\sqrt{Z}}{\sqrt{U}}+\frac{q^\theta}{\varphi_{\theta} (q)}\frac{\mathds{1}_{\{U\geq 10^{12}\}}(U)}{103\log\left(\frac{U}{Z}\right)}\right)$.  

\begin{lemma}\label{S2:total}  
Let $1<Z<U$ such that $\frac{U}{Z}\geq 20$ and $Z\geq 4\times 10^5$. Then 
\begin{align*}
\mathit{S}_{\mathbf{II}}^{(2)}=\frac{v}{\varphi(v)}\left(\log(Z)+\gamma+\sum_{p|v}\frac{\log(p)}{p-1}\right)+O^*\left(\frac{\Upsilon^{(2)}_v}{\sqrt{Z}}+\frac{\Upsilon^{(3)}_v\ Z}{U}\right),
\end{align*}
\begin{flalign*}  
\text{where}&&\Upsilon^{(2)}_v&=\begin{cases} 
\sage{Upper(Upsilon2_v1,digits)}&\text{ if }v=1,\\
\sage{Upper(Upsilon2_v2,digits)}&\text{ if }v=2,
\end{cases}\quad
\Upsilon^{(3)}_v=\begin{cases} 
\sage{Upper(Upsilon3_v1,digits)}&\text{ if }v=1,\\
\sage{Upper(Upsilon3_v2,digits)}&\text{ if }v=2.
\end{cases} &
\end{flalign*} 
\end{lemma}  

With respect to $Z$, it will be clear in \S\ref{S23} and \S\ref{calcul} why we will end up selecting $Z=\mathrm{c}U^{\sage{choice}}$, with $\mathrm{c}\in\{\sage{ccc},\sage{ccc2},\sage{ccc3}\}$. With this choice, we will deduce the following.

\begin{lemma}\label{S3:total}    
Let $v\in\{1,2\}$. We have the following estimation
\begin{equation*}
|\mathit{S}_{\mathbf{II}}^{(3)}|\leq \frac{\Upsilon^{(4)}_v\log(U)}{\sage{1/(1-choice)}U^{\sage{1-choice}}}+\frac{\Upsilon^{(4)}_v}{U^{\sage{1-choice}}}+\frac{\mathds{1}_{\{U\geq 10^{12}\}}(U)\Upsilon^{(5)}_v}{\log(U)},
\end{equation*} 
\begin{flalign*}
&\text{where }\quad &&\Upsilon^{(4)}_v&&= 
\begin{cases} 
\sage{Upper(Upsilon4_v1,digits)}&\text{ if }v=1,\\ 
\sage{Upper(Upsilon4_v2,digits)}&\text{ if }v=2,
\end{cases}\ \text{ and,}\\
&\text{if }U\geq 10^{\sage{exp_prog}}\text{ and }Z=\sage{ccc}U^{\sage{choice}}\ ,\quad &&\Upsilon^{(5)}_v&&= 
\begin{cases} 
\sage{Upper(Upsilon5_v1,digits+1)}&\text{ if }v=1,\\ 
\sage{Upper(Upsilon5_v2,digits+1)}&\text{ if }v=2, 
\end{cases}\\
&\text{if }U\geq 10^{\sage{exp_prog_bound}}\text{ and }Z=\sage{ccc2}U^{\sage{choice}}\ ,\quad &&\Upsilon^{(5)}_v&&= 
\begin{cases} 
\sage{Upper(Upsilon5_v1_2,digits+1)}&\text{ if }v=1,\\ 
\sage{Upper(Upsilon5_v2_2,digits+1)}&\text{ if }v=2,
\end{cases} &\\
&\text{if }U\geq 10^{\sage{exp_prog_bound}}\text{ and }Z=\sage{ccc3}U^{\sage{choice}}\ ,\quad &&\Upsilon^{(5)}_v&&= 
\begin{cases} 
\sage{Upper(Upsilon5_v1_3,digits+1)}&\text{ if }v=1,\\ 
\sage{Upper(Upsilon5_v2_3,digits+1)}&\text{ if }v=2.
\end{cases} &
\end{flalign*}
\end{lemma} 

In \S\ref{Conclusion}, we will combine lemmas \ref{SI}, \ref{S1:total}, \ref{S2:total} and \ref{S3:total} to derive Theorem \ref{CONCLUSION}.

\subsection{The sum $\mathit{S}_{\mathbf{I}}$}\label{A}

We introduce some lemmas that will help us to estimate $\mathit{S}_{\mathbf{I}}$. The reader should keep in mind that although many of our results consider only $q=v\in\{1,2\}$, they can be stated for any $q\in\mathbb{Z}_{>0}$. 
We start quoting \cite[Lemma 4.4]{1ZA20}. 
 \begin{lemma}\label{sum2}
Let $X>0$. The following estimation holds. 
\begin{align}\label{sum2:eq} 
\sum_{\substack{\ell\leq X\\(\ell,q)}}\frac{\mu^2(\ell)}{\ell}&=\frac{q}{\kappa(q)}\frac{6}{\pi^2}\left(\log(X)+\mathfrak{b}_q\right)+O^*\left(\frac{\sqrt{q}}{\varphi_{\frac{1}{2}}(q)}\frac{\sage{Trunc(Cc1,digits)}\prod_{2|q}\sage{Upper(Cc2/Cc1,digits)}}{\sqrt{X}}\right), 
\end{align}
\begin{flalign*}  
\text{where}&&\mathfrak{b}_q =\sum_{p}\frac{2\log(p)}{p^2-1}+\gamma+\sum_{p|q}\frac{\log(p)}{p+1}, \quad
\sum_{p}\frac{2\log(p)}{p^2-1}+\gamma=\sage{C_sum2}\ldots.\phantom{xxxxxx}
\end{flalign*} 
 \end{lemma} 

\begin{proposition}\label{sq1.half}    
Let $X\geq 10$ and $v\in\{1,2\}$. Then
\begin{align*}
\frac{1}{\log^2(X)}\times\sum_{\substack{\ell\leq X\\(\ell,v)=1}}\frac{\mu^2(\ell)}{\sqrt{\ell}\ \varphi_{\frac{1}{2}}(\ell)}\log\left(\frac{X}{\ell}\right)\leq\psi_v=\begin{cases}
\sage{Upper(constant_psi_v1,digits)},&\text{ if }v=1,\\ 
\sage{Upper(constant_psi_v2,digits)},&\text{ if }v=2.\\
\end{cases} 
\end{align*} 
\end{proposition}
 \begin{proof} Let $q\in\mathbb{Z}_{>0}$. By applying \cite[Thm. 3.3]{1ZA20} with $f(p)=\frac{1}{\sqrt{p}\ \varphi_{\frac{1}{2}}(p)}=\frac{1}{\sqrt{p}(\sqrt{p}-1)}$, $\alpha=1$, $\beta=\frac{3}{2}$ and $0\leq\delta=\sage{delta_sqHalf}<\frac{1}{2}$, we obtain that
\begin{align}\label{eq:sqhalf}
\sum_{\substack{\ell\leq X\\(\ell,q)=1}}\frac{\mu^2(\ell)}{\sqrt{\ell}\ \varphi_{\frac{1}{2}}(\ell)}=\mathrm{k}_q\mathbf{F}\left(\log(X)+\mathfrak{f}_q\right)+O^*\left(\frac{\mathrm{l}_q\ \mathbf{f}}{X^{\sage{delta_sqHalf }}}\right),
\end{align} 
where
 \begin{align*}
\mathrm{k}_q&=\prod_{p|q}\left(1-\frac{p+1}{p^2-p^{\frac{3}{2}}+\sqrt{p}+1}\right),\quad\mathrm{l}_q=\prod_{p|q}\left(1+\frac{p^{\sage{1-delta_sqHalf }}-2p^{\sage{1/2-delta_sqHalf }}-1}{p^{\sage{2-2*delta_sqHalf }}-p^{\sage{3/2-2*delta_sqHalf}}+p^{\sage{1-delta_sqHalf}}+1}\right),\phantom{xxxxx}\\
\mathfrak{f}_q&=-\sum_{p}\frac{(\sqrt{p}-2)\log(p)}{(p-\sqrt{p}+1)(p-1)}+\gamma+\sum_{p|q}\frac{\log(p)}{p-\sqrt{p}+1},\phantom{xxx}\\
 \mathbf{F}&=\frac{\zeta\left(\frac{3}{2}\right)}{\zeta\left(3\right)} \in [\sage{Lower(Prod_sqHalf_l,digits)},\sage{Upper(Prod_sqHalf_u,digits)}],\quad -\sum_{p}\frac{(\sqrt{p}-2)\log(p)}{(p-\sqrt{p}+1)(p-1)}+\gamma\in [\sage{Lower(Sum_sqHalf_l,digits)},\sage{Upper(Sum_sqHalf_u,digits)}],\\ 
  \mathbf{f}&=\Delta_1^{\sage{delta_sqHalf}}\prod_{p}\left(1+\frac{p^{\sage{1/2-delta_sqHalf }}+1}{p^{\sage{3/2-2*delta_sqHalf}}(\sqrt{p}-1)}\right) \in [\sage{Lower(A(1,delta_sqHalf)*delta_sqHalf_l,digits)},\sage{Upper(A(1,delta_sqHalf )*delta_sqHalf_u,digits)}].
\end{align*}
Therefore, when $q=v\in\{1,2\}$ and $X\geq C=10^7$, we derive from \eqref{eq:sqhalf} that $\frac{1}{\log^2(X)}\times\sum_{\substack{\ell\leq X\\(\ell,v)=1}}\frac{\mu^2(\ell)}{\sqrt{\ell}\ \varphi_{\frac{1}{2}}(\ell)}\log\left(\frac{X}{\ell}\right)$ can be expressed as
\begin{align}\label{analytic_sqHalf} 
\int_1^X\left(\mathrm{k}_v\mathbf{F}\left(\log(t)+\mathfrak{f}_v\right)+O^*\left(\frac{\mathrm{l}_v\ \mathbf{f}}{t^{\sage{delta_sqHalf}}}\right)\right)
\frac{dt}{t\log^2(X)}\phantom{xxxxxxxxxxxxxxxxxxx}&\nonumber\\
\leq\sage{Upper(zeta(3/2)/zeta(3),digits)}\ \mathrm{k}_v\left(\frac{1}{2}+\frac{\mathfrak{f}_v}{\log(C)}\right)+\frac{\sage{Upper(1/delta_sqHalf * err_sqHalf,digits)}\ \mathrm{l}_v}{\log^2(C)}=
\begin{cases}
\sage{Upper(analytic_sqHalf_v1,digits)}&\text{ if }v=1,\\
\sage{Upper(analytic_sqHalf_v2,digits)}&\text{ if }v=2.
\end{cases}&
\end{align} 
where we have used that $\mathfrak{f}_q>0$ for all $q\in\mathbb{Z}_{>0}$. 
On the other hand, for all $10\leq X\leq 10^7$,  
 \begin{align} 
  \frac{1}{\log^2(X)}\times\sum_{\substack{\ell\leq X\\(\ell,v)=1}}\frac{\mu^2(\ell)}{\sqrt{\ell}\ \varphi_{\frac{1}{2}}(\ell)}\log\left(\frac{X}{\ell}\right)\leq
\begin{cases}\label{program_sqHalf}
\sage{Upper(program_sqHalf_v1,digits)}&\text{ if }v=1,\\
\sage{Upper(program_sqHalf_v2,digits)}&\text{ if }v=2.\\
\end{cases}
\end{align} 
The result is concluded by defining $\psi_v$ as the maximum between the bounds given in \eqref{analytic_sqHalf} and \eqref{program_sqHalf}.
 \end{proof}  

\begin{proofbold}{Lemma~\ref{SI:estimation}}
Conditions $\ell r_i\leq U$ and $Z<\ell$ imply that $r_i\leq\frac{U}{Z}$ for $i=1,2$. Therefore, from definition \eqref{def_S_I}, we derive
\begin{align} 
\mathit{S}_{\mathbf{I}}=\sum_{\substack{r_1,r_2\leq\frac{U}{Z}\\ (r_1,r_2)=1\\(r_1r_2,v)=1}}\frac{\mu(r_1)\mu(r_2)}{r_1r_2}\sum_{\substack{Z<\ell\leq U\\(\ell,r_1r_2v)=1}}\frac{\mu^2(\ell)}{\ell}\mathbf{L}_{\ell r_1}\mathbf{L}_{\ell r_2}.\label{expri}
\end{align} 
On the other hand, with the help of Lemma \ref{sum2}, for any $t>Z$, we have that
\begin{equation*}\label{Ram}
A_q(t)=\sum_{\substack{Z<\ell\leq t\\ (\ell,q)=1}}\frac{\mu^2(\ell)}{\ell}=\frac{6}{\pi^2}\frac{q}{\kappa(q)}\log\left(\frac{t}{Z}\right)+O^*\left(\frac{\sqrt{q}}{\varphi_{\frac{1}{2}}(q)}\frac{\sage{Trunc(2*Cc1,digits)}\prod_{2|q}\sage{Upper(Cc2/Cc1,digits)}}{\sqrt{Z}}\right).
\end{equation*} 
Moreover, by considering a monotone continuous function $\mathbf{L^*}$ on $[1,U]$, such that $\mathbf{L}^*(U)=0$, as $\mathbf{L^*}$ is of bounded variation, we can apply summation by parts and derive 
\begin{align}\label{boundvar}
\sum_{\substack{Z<\ell\leq U\\(\ell,q)=1}}\frac{\mu^2(\ell)}{\ell}\mathbf{L^*}(\ell)&=-\int_{Z}^UA_q(t)d\mathbf{L^*}(t)\nonumber\\
=&\frac{6}{\pi^2}\frac{q}{\kappa(q)}\int_{Z}^{U}\frac{\mathbf{L^*}(t)}{t}dt+O^*\left(\frac{\sage{Trunc(2*Cc1,digits)}\prod_{2|q}\sage{Upper(Cc2/Cc1,digits)}\ \sqrt{q}}{\varphi_{\frac{1}{2}}(q)}\frac{|\mathbf{L^*}(Z)|}{\sqrt{Z}}\right), 
\end{align}
since $\left[\phantom{\frac{}{}}A_q(t)\mathbf{L^*}(t)\right|_{Z}^{U}=\left[\log\left(\frac{t}{Z}\right)\mathbf{L^*}(t)\right|_{Z}^{U}=0$ and $\int_{Z}^Ud\mathbf{L^*}(t)=-\mathbf{L^*}(Z)$.

In particular, by taking $\mathbf{L^*}(t)=\mathbf{L}_{tr_1}\mathbf{L}_{tr_2}=\log^+\left(\frac{U}{tr_1}\right)\log^+\left(\frac{U}{tr_2}\right)$, with $r_1,r_2\leq\frac{U}{Z}$, we have a monotone decreasing function on $(0,\infty)$, thus of bounded variation, such that $\mathbf{L}_{Ur_1}\mathbf{L}_{Ur_2}=0$ and $\mathbf{L^*}(Z)=\mathbf{L}_{Zr_1}\mathbf{L}_{Zr_2}=\log\left(\frac{U}{Zr_1}\right)\log\left(\frac{U}{Zr_2}\right)$. Further, with this choice of $\mathbf{L}$, by taking $q=r_1r_2v$ and replacing \eqref{boundvar} into the innermost summation of \eqref{expri}, $\mathit{S}_{\mathbf{I}}$ equals
\begin{align}\label{mainterm:S2}   
\sum_{\substack{r_1,r_2\leq\frac{U}{Z}\\ (r_1,r_2)=1\\(r_1r_2,v)=1}}\frac{\mu(r_1)\mu(r_2)}{r_1r_2}\left(\frac{6}{\pi^2}\frac{r_1r_2v}{\kappa(r_1r_2v)}\int_{Z}^{U}\frac{\mathbf{L}_{tr_1}\mathbf{L}_{tr_2}}{t}dt\right.\phantom{xxxxxxxxxxxxxxxxxxxxx}&\nonumber\\
\left.+O^*\left(\frac{\sage{Trunc(2*Cc1,digits)}\prod_{2|r_1r_2v}\sage{Upper(Cc2/Cc1,digits)}\ \sqrt{r_1r_2v}}{\varphi_{\frac{1}{2}}(r_1r_2v)}\frac{\log\left(\frac{U}{Zr_1}\right)\log\left(\frac{U}{Zr_2}\right)}{\sqrt{Z}}\right)\right)&\nonumber\\  
=\frac{6}{\pi^2}\frac{v}{\kappa(v)}\int_{Z}^{U}\sum_{\substack{r_1,r_2\\ (r_1,r_2)=1\\(r_1r_2,v)=1}}\frac{\mu(r_1)\mu(r_2)}{\kappa(r_1)\kappa(r_2)}\frac{\mathbf{L}_{tr_1}\mathbf{L}_{tr_2}}{t}dt+\mathbf{r}_v(U),  \phantom{xxxx}& 
\end{align}  
where condition $r_i\leq\frac{U}{Z}$ above is encoded by the definition of $\mathbf{L}_{tr_i}$, $i\in\{1,2\}$, and by the range of $t$. The remainder term $\mathbf{r}_v$ can be estimated by defining $Q_v:X\mapsto\sum_{\substack{r\leq X\\ (r,v)=1}}\frac{\mu^2(r)}{\sqrt{r}\varphi_{\frac{1}{2}}(r)}\log\left(\frac{X}{r}\right)$ and using Proposition \ref{sq1.half} as follows
\begin{align}   
|\mathbf{r}_1(U)|\leq \frac{\sage{Trunc(2*Cc1,digits)}}{\sqrt{Z}}\left(\frac{2\times\sage{Trunc(Cc2/Cc1,digits)}}{\sqrt{2}(\sqrt{2}-1)}Q_2\left(\frac{U}{2Z}\right)Q_1\left(\frac{U}{Z}\right)+Q_1^2\left(\frac{U}{Z}\right)\right)\phantom{x}\label{source}&\\
\leq \frac{\sage{Trunc(2*Cc1,digits)}\left(\sage{Upper(coeff_r_v1,digits)}\ \psi_1\psi_2+\psi_1^2\right) \log^4\left(\frac{U}{Z}\right)}{\sqrt{Z}}\leq\frac{\sage{Upper(Tau1_v1,digits)}\log^4\left(\frac{U}{Z}\right)}{\sqrt{Z}}=\frac{T^{(1)}_1\log^4\left(\frac{U}{Z}\right)}{\sqrt{Z}},\nonumber&\\
|\mathbf{r}_2(U)|\leq \frac{\sage{Trunc(2*Cc2,digits)}\ \sqrt{2}}{\varphi_{\frac{1}{2}}(2)\ \sqrt{Z}}\ Q_2^2\left(\frac{U}{Z}\right)\leq\frac{\sage{Upper(Tau1_v2,digits)}\log^4\left(\frac{U}{Z}\right)}{\sqrt{Z}}=\frac{T^{(1)}_2\log^4\left(\frac{U}{Z}\right)}{\sqrt{Z}},\label{source2}&
\end{align} 
where we have used that $\frac{U}{2Z}\geq 10$. 
  
With respect to the main term of $\mathit{S}_{\mathbf{I}}$ given in \eqref{mainterm:S2}, recall the function $\tilde{m}$ defined in \eqref{tildem} and observe that  
\begin{align}\label{idddd}
\int_{Z}^{U}\sum_{\substack{r_1,r_2\\ (r_1,r_2)=1\\(r_1r_2,v)=1}}\frac{\mu(r_1)\mu(r_2)}{\kappa(r_1)\kappa(r_2)}\frac{\mathbf{L}_{tr_1}\mathbf{L}_{tr_2}}{t}dt&=\int_1^{\frac{U}{Z}}\sum_{\substack{d\\(d,v)=1}}\frac{\mu(d)}{\kappa(d)^2}\ \tilde{m}_{dv}^2\left(\frac{s}{d}\right)\frac{ds}{s},
\end{align} 
where the change of variables $s=\frac{U}{t}$ has been performed, which is valid by Lemma \ref{integral}. Hence, by combining Lemma \ref{integral} and equations \eqref{h_q}, \eqref{idddd}, $\int_{Z}^{U}\sum_{\substack{r_1,r_2\\ (r_1,r_2)=1\\(r_1r_2,v)=1}}\frac{\mu(r_1)\mu(r_2)}{\kappa(r_1)\kappa(r_2)}\frac{\mathbf{L}_{tr_1}\mathbf{L}_{tr_2}}{t}dt$ equals
\begin{align*}
\int_1^{\frac{U}{Z}}\frac{h_v(s)}{s}ds+\frac{\pi^2\kappa(v)}{6v}\check{\check{m}}_v\left(\frac{U}{Z}\right)-\frac{\pi^2\kappa(v)}{6\varphi(v)}\log\left(\frac{U}{Z}\right).
\end{align*} 
Thus, by recalling \eqref{mainterm:S2}, $\mathit{S}_{\mathbf{I}}$ may be expressed as
\begin{align}
\label{SI:estimation1} 
\frac{6}{\pi^2}\frac{v}{\kappa(v)}\int_1^{\frac{U}{Z}}\frac{h_v(s)}{s}ds+\check{\check{m}}_v\left(\frac{U}{Z}\right)-\frac{v}{\varphi(v)}\log\left(\frac{U}{Z}\right)+O^*\left(\frac{T_v^{(1)}\log^4\left(\frac{U}{Z}\right)}{\sqrt{Z}}\right). 
\end{align}
Furthermore, by Proposition \ref{tail}, $\int_1^\infty \frac{h_v(s)}{s}ds$ converges. Therefore, from Equation \eqref{SI:estimation1}, we may write 
\begin{align} 
\label{SI:estimation2} 
\mathit{S}_{\mathbf{I}}=&\frac{6v}{\pi^2\kappa(v)}\int_1^{\infty}\frac{h_v(s)}{s}ds+\check{\check{m}}_v\left(\frac{U}{Z}\right)-\frac{v}{\varphi(v)}\log\left(\frac{U}{Z}\right)\nonumber\\
&\phantom{xxxxxxxxx}+O^*\left(\frac{T_v^{(1)}\log^4\left(\frac{U}{Z}\right)}{\sqrt{Z}}+\frac{6v}{\pi^2\kappa(v)}\int_{\frac{U}{Z}}^{\infty}\frac{|h_v(s)|}{s}ds\right).
\end{align} 
 
Now, by Proposition \ref{TT}, whenever $10^{12}\leq\frac{U}{Z}$, we have
\begin{align}\label{hv11}
\left|\int_{\frac{U}{Z}}^\infty \frac{h_v(s)}{s}ds\right|&\leq\int_{\frac{U}{Z}}^\infty\frac{T_v^{(4)}}{s\log^2(s)}ds=  \frac{T_v^{(4)}}{\log\left(\frac{U}{Z}\right)},
\end{align}
whereas, if $20\leq\frac{U}{Z}\leq 10^{12}$, we have
\begin{align}\label{hv22} 
\left|\int_{\frac{U}{Z}}^\infty \frac{h_v(s)}{s}ds\right|\leq \Psi_v\left(\frac{Z\log\left(\frac{U}{Z}\right)}{U}+\frac{Z}{U}\right)+\frac{T_v^{(4)}}{\log\left(\frac{U}{Z}\right)},
\end{align} 
 where $\Psi_v=\left(T_v^{(2)}+\frac{T_v^{(3)}}{\log(20)}\right)$, by using that $\int\frac{\log(s)}{s^2}ds=-\frac{\log(s)+1}{s}$,. 

Finally, from \eqref{hv11}, \eqref{hv22} and the definitions of $T_v^{(1)},T_v^{(2)},T_v^{(3)}$, $T_v^{(4)}$ and $\Psi_v$, we derive the result.
\end{proofbold}

\subsection{The sum $\mathit{S}_{\mathbf{II}}^{(1)}$}\label{S21}

By \cite[Lemma 4.7]{1ZA20}, we have that 
 \begin{lemma}\label{sum1}
 Let $X>0$. The following estimation holds
\begin{equation}\label{sum1:eq}    
\sum_{\substack{\ell\leq X\\(\ell,q)=1}}\frac{\mu^2(\ell)}{\varphi(\ell)}=\frac{\varphi(q)}{q}\left(\log\left(X\right)+\mathfrak{a}_q\right)+O^*\left(\frac{\mathrm{A}_q\ \sage{Trunc(constant_ram,1)}\prod_{2|q}\sage{Upper(constant_ram_v2,digits)}}{\sqrt{X}}\right), 
\end{equation}
 where $\mathfrak{a}_q$ is defined as in Lemma \ref{tildetildem}, and $\mathrm{A}_q=\prod_{p|q}\left(1+\frac{p-2}{p^{\frac{3}{2}}-p-\sqrt{p}+2}\right)$.
 \end{lemma}  

\begin{proposition}\label{sumvar1log}     
Let $X\geq 10$ and $v\in\{1,2\}$. Then 
 \begin{align*} 
  \frac{1}{\log^2(X)}\times\sum_{\substack{\ell\leq X\\ (\ell,v)=1}}\frac{\mu^2(\ell)\mathrm{A}_\ell}{\varphi(\ell)}\log\left(\frac{X}{\ell}\right)\leq
\eta_v = 
\begin{cases}
\sage{Upper(constant_eta_v1,digits)}\quad&\text{ if }v=1,\\
\sage{Upper(constant_eta_v2,digits)}\quad&\text{ if }v=2.
\end{cases}
  \end{align*}  
 \end{proposition}

  \begin{proof}
  As $\sum_{\substack{\ell\leq X\\ (\ell,q)=1}}\frac{\mu^2(\ell)\mathrm{A}_\ell}{\varphi(\ell)}\log\left(\frac{X}{\ell}\right)=\int_{1}^X\left(\sum_{\substack{\ell\leq t\\ (\ell,q)=1}}\frac{\mu^2(\ell)\mathrm{A}_\ell}{\varphi(\ell)}\right)\frac{dt}{t}$, it suffices to analyze the sum inside the integral. By \cite[Thm. 3.3]{1ZA20}, with $f(p)=\frac{\mathrm{A}_p}{p-1}$, $\alpha=1$, $\beta=\frac{3}{2}$ and $\delta=\sage{delta_sumvar1log}$, we derive that
\begin{align*}
\sum_{\substack{\ell\leq X\\ (\ell,q)=1}}\frac{\mu^2(\ell)\mathrm{A}_\ell}{\varphi(\ell)}=j_q\mathbf{G}\left(\log(X)+\mathfrak{g}_q\right)+O^*\left(\frac{k_q\ \mathbf{g}}{X^{\sage{delta_sumvar1log}}}\right)
\end{align*} 
where
\begin{align*}
 &j_q=\prod_{p|q}\left(1-\frac{\mathrm{A}_p}{p-1+\mathrm{A}_p}\right),\quad\mathbf{G}=\prod_p\left(1+\frac{\mathrm{A}_p-1}{p}\right)\in[\sage{Lower(Prod_sumvar1log_l,digits)},\sage{Upper(Prod_sumvar1log_u,digits)}],\\
& \mathfrak{g}_q=\sum_{p}{\frac{\log(p)(p-1-(p-2)\mathrm{A}_p)}{(\mathrm{A}_p+p-1)(p-1)}}+\gamma+\sum_{p|q}{\frac{\log(p)\mathrm{A}_p}{\mathrm{A}_p+p-1}},\\
& k_q=\prod_{p|q}\left(1+\frac{2(p-1)-\mathrm{A}_p(p+p^{\sage{delta_sumvar1log}})}{(p-1)p^{\sage{1-delta_sumvar1log}}+\mathrm{A}_p(p+p^{\sage{delta_sumvar1log}})-p+1}\right),\\
&\sum_{p}{\frac{\log(p)(p-1-(p-2)\mathrm{A}_p)}{(\mathrm{A}_p+p-1)(p-1)}}+\gamma\in[\sage{Lower(Sum_sumvar1log_l+gamma,digits)},\sage{Upper(Sum_sumvar1log_u+gamma,digits)}],\\ 
&\mathbf{g}=\Delta_1^{\sage{delta_sumvar1log}}\prod_{p}\left(1+\frac{p(\mathrm{A}_p-1)+\mathrm{A}_pp^{\sage{delta_sumvar1log}}+1}{(p-1)p^{\sage{1-delta_sumvar1log}}}\right)\in[\sage{Lower(A(1,delta_sumvar1log)*deltaProd_sumvar1log_l,digits)},\sage{Upper(A(1,delta_sumvar1log)*deltaProd_sumvar1log_u,digits)}].
  \end{align*}
Therefore, as $\int_{1}^X\frac{\log(t)}{t}dt=\frac{\log^2(X)}{2}$ and $\mathfrak{g}_q>0$ for all $q\in\mathbb{Z}_{>0}$, we derive for all $X\geq C=10^7$ and $q=v\in\{1,2\}$ that $ \frac{1}{\log^2(X)}\times\sum_{\substack{\ell\leq X\\ (\ell,v)=1}}\frac{\mu^2(\ell)\mathrm{A}_\ell}{\ell}\log\left(\frac{X}{\ell}\right)$ may be estimated as
 \begin{align} 
 \int_{1}^X\left(j_v\mathbf{G}\left(\log(t)+\mathfrak{g}_v\right)+O^*\left(\frac{k_v\ \mathbf{g}}{t^{\sage{delta_sumvar1log}}}\right)\right)\frac{dt}{t\log^2(X)}\phantom{xxxxxxxxxxxxxxxxx}&\nonumber\\ 
\leq\ \sage{Upper(Prod_sumvar1log_u,digits)}\ j_v\left(\frac{1}{2}+\frac{\mathfrak{g}_v}{\log(C)}\right)+\frac{\sage{Upper(1/delta_sumvar1log * errProd_sumvar1log ,digits)}\ k_v}{\log^2(C)}=
\begin{cases}\label{analytic_sumvar1log}
\sage{Upper(analytic_sumvar1log_v1,digits)}&\text{ if }v=1,\\
\sage{Upper(analytic_sumvar1log_v2,digits)}&\text{ if }v=2.
\end{cases}&
\end{align} 
On the other hand, for all $10\leq X\leq 10^7$,
 \begin{align} 
  \frac{1}{\log^2(X)}\times\sum_{\substack{\ell\leq X\\ (\ell,v)=1}}\frac{\mu^2(\ell)\mathrm{A}_\ell}{\ell}\log\left(\frac{X}{\ell}\right)\leq
\begin{cases}\label{program_sumvar1log}
\sage{Upper(program_sumvar1log_v1,digits)}&\text{ if }v=1,\\ 
\sage{Upper(program_sumvar1log_v2,digits)}&\text{ if }v=2.
\end{cases}
\end{align}  
The result is concluded by defining $\eta_v$ as the maximum between the bounds given in \eqref{analytic_sumvar1log} and \eqref{program_sumvar1log}.
\end{proof} 
 
\begin{proofbold}{Lemma~\ref{S1:total}}
From \eqref{S1:eq} and the definition of $\check{\check{m}}_q$ given in \eqref{checks}, observe that 
\begin{align}\label{compress} 
\mathit{S}_{\mathbf{II}}^{(1)}&=\frac{v}{\varphi(v)}\sum_{\substack{\ell\leq Z\\(\ell,v)=1}}\frac{\mu^2(\ell)}{\varphi(\ell)}\sum_{\substack{n\leq\frac{U}{\ell}\\(n,\ell v)=1}}\frac{\mu(n)}{\varphi(n)}\log\left(\frac{U}{\ell n}\right),
\end{align}
where we have used that for any square-free $n$, $\sum_{d|n}\frac{1}{\varphi(d)}=\frac{n}{\varphi(n)}$. Moreover, from \eqref{compress}, we obtain that $\frac{\varphi(v)}{v}\mathit{S}_{\mathbf{II}}^{(1)}$ equals 
\begin{align}\label{dd2}
\sum_{\substack{\ell\leq U\\(\ell,v)=1}}\frac{\mu^2(\ell)}{\varphi(\ell)}\sum_{\substack{n\leq\frac{U}{\ell}\\(n,\ell v)=1}}\frac{\mu(n)}{\varphi(n)}\log\left(\frac{U}{\ell n}\right)-\sum_{\substack{Z<\ell\leq U\\(\ell,v)=1}}\frac{\mu^2(\ell)}{\varphi(\ell)}\sum_{\substack{n\leq\frac{U}{\ell}\\(n,\ell v)=1}}\frac{\mu(n)}{\varphi(n)}\log\left(\frac{U}{\ell n}\right)&\nonumber\\
= \log(U)-\sum_{\substack{n\leq\frac{U}{Z}\\(n,v)=1}}\frac{\mu(n)}{\varphi(n)}\sum_{\substack{Z<\ell\leq \frac{U}{n}\\(\ell,nv)=1}}\frac{\mu^2(\ell)}{\varphi(\ell)}\log\left(\frac{U}{\ell n}\right),&
\end{align}
where in the above first summation, we have used M\"obius inversion.

Now, on using Lemma \ref{sum1} and summation by parts, we deduce that $\sum_{\substack{Z<\ell\leq \frac{U}{n}\\(\ell,nv)=1}}\frac{\mu^2(\ell)}{\varphi(\ell)}\log\left(\frac{U}{\ell n}\right)$ may be estimated as
\begin{align} 
\int_{Z}^{\frac{U}{n}}\left(\frac{\varphi(nv)}{nv}\log\left(\frac{t}{Z}\right)+O^*\left(\frac{\mathrm{A}_{nv}\ \sage{Upper(2*constant_ram,digits)}\prod_{2|nv}\sage{Upper(constant_ram_v2,digits)}}{\sqrt{Z}}\right)\right)\frac{dt}{t}\phantom{xxxxxxxxxxxx}&\nonumber\\
=\frac{\varphi(nv)}{2nv}\log^2\left(\frac{\frac{U}{Z}}{n}\right)+O^*\left(\frac{\mathrm{A}_{nv}\ \sage{Upper(2*constant_ram,digits)}\prod_{2|nv}\sage{Upper(constant_ram_v2,digits)}\log\left(\frac{\frac{U}{Z}}{n}\right)}{\sqrt{Z}}\right),&\label{followup}
\end{align}
Replacing \eqref{followup} into the second term of \eqref{dd2} gives further  
\begin{align*}
\sum_{\substack{n\leq\frac{U}{Z}\\(n,v)=1}}\frac{\mu(n)}{\varphi(n)}\sum_{\substack{Z<\ell\leq \frac{U}{n}\\(\ell,nv)=1}}\frac{\mu^2(\ell)}{\varphi(\ell)}\log\left(\frac{U}{\ell n}\right)=\frac{\varphi(v)}{2v}\sum_{\substack{n\leq\frac{U}{Z}\\(n,v)=1}}\frac{\mu(n)}{n}\log^2\left(\frac{\frac{U}{Z}}{n}\right)\phantom{xxxxxx}&\nonumber\\
+O^*\left(\frac{\sage{Upper(2*constant_ram,digits)}\ \mathrm{A}_{v}}{\sqrt{Z}}\sum_{\substack{n\leq\frac{U}{Z}\\(n,v)=1}}\frac{\prod_{2|nv}\sage{Upper(constant_ram_v2,digits)}\ \mu^2(n)\mathrm{A}_{n}}{\varphi(n)}\log\left(\frac{\frac{U}{Z}}{n}\right)\right).&\nonumber
\end{align*}
The above main term corresponds to $\frac{\varphi(v)}{2v}\ \check{\check{m}}_v\left(\frac{U}{Z}\right)$. As for the error term, it can be estimated by means of Proposition \ref{sumvar1log}: if $v=2$, the factor $\prod_{2|nv}\sage{Upper(constant_ram_v2,digits)}$ is always present, whereas, if $v=1$, we have
\begin{align}\label{ethos}
\sum_{\substack{n\leq\frac{U}{Z}\\(n,v)=1}}\frac{\prod_{2|nv}\sage{Upper(constant_ram_v2,digits)}\ \mu^2(n)\mathrm{A}_{n}}{\varphi(n)}\log\left(\frac{\frac{U}{Z}}{n}\right)=\sum_{\substack{n\leq\frac{U}{Z}}}\frac{\prod_{2|n}\sage{Upper(constant_ram_v2,digits)}\ \mu^2(n)\mathrm{A}_{n}}{\varphi(n)}\log\left(\frac{\frac{U}{Z}}{n}\right)&\nonumber\\
=\sage{Upper(constant_ram_v2,digits)}\ \mathrm{A}_2\sum_{\substack{n\leq\frac{U}{2Z}\\(n,2)=1}}\frac{\mu^2(n)\mathrm{A}_{n}}{\varphi(n)}\log\left(\frac{\frac{U}{2Z}}{n}\right)+\sum_{\substack{n\leq\frac{U}{Z}\\(n,2)=1}}\frac{\mu^2(n)\mathrm{A}_{n}}{\varphi(n)}\log\left(\frac{\frac{U}{Z}}{n}\right)&\nonumber\\
\leq\sage{Upper(constant_ram_v2,digits)}\ \eta_2\log^2\left(\frac{U}{2Z}\right)+\eta_2\log^2\left(\frac{U}{Z}\right)\leq\left(\sage{Upper(constant_ram_v2,digits)}\ \eta_2+\eta_2\right)\log^2\left(\frac{U}{Z}\right),&
\end{align} 
where we have used that $\mathrm{A}_2=1$ and, since $\frac{U}{2Z}\geq 10$, Proposition \ref{sumvar1log}. We conclude the result by defining 
\begin{align*}
\Upsilon^{(1)}_v=
\begin{cases}
\sage{Upper(2*constant_ram*constant_ram_v2,digits)}\ \eta_2+\sage{Upper(2*constant_ram,digits)}\ \eta_2 &\text{ if }v=1,\\ 
\sage{Upper(2*constant_ram*constant_ram_v2,digits)}\ \eta_2&\text{ if }v=2.
\end{cases}
\end{align*}
\end{proofbold}

\begin{remark}\label{remark2}
The error term given in \eqref{followup} has arisen due to the fact that we are studying a sum whose range starts sufficiently away from $1$; if the summation $\sum_{\substack{Z<\ell\leq \frac{U}{n}\\(\ell,nv)=1}}\frac{\mu^2(\ell)}{\varphi(\ell)}\log\left(\frac{U}{\ell n}\right)$ would have started from $1$ (or any admissible constant value) rather than $Z$, a remainder term of order $\log^2(U)$ would have appeared in Lemma \ref{S1:total}, thus not even providing an asymptotic estimation for $\mathit{S}_{\mathbf{II}}^{(1)}$. This fact justifies why we have split the expression \eqref{split} into two summations.
\end{remark}

\subsection{The sum $\mathit{S}_{\mathbf{II}}^{(2)}$}\label{S22}

We need a series of lemmas that rely on an interval arithmetic computations within a range, using specifically that $v\in\{1,2\}$. As those calculations may be performed for any $q\in\mathbb{Z}_{>0}$, Theorem \ref{CONCLUSION} holds true.

\begin{proposition}\label{sumvarp}        
Let $X\geq 20$ and $v\in\{1,2\}$. Then   
\begin{align*}  
X\times\sum_{\substack{\ell> X\\(\ell,v)=1}}\frac{\mu^2(\ell)}{\varphi(\ell)^2}&\leq\varphi_v^{(1)}=
\begin{cases}
\sage{Upper(constant1_phi_v1,digits)}&\text{ if }v=1,\\
\sage{Upper(constant1_phi_v2,digits)}&\text{ if }v=2.
\end{cases}
\end{align*}
\end{proposition} 

 \begin{proof}
 By applying \cite[Thm. 4.3.1]{1ZA20} with $f(p)=\frac{1}{\varphi(p)^2}=\frac{1}{(p-1)^2}$, $\alpha=2$ and $\beta=3$, we have
 \begin{align}\label{ready}
 \sum_{\substack{\ell\leq X\\(\ell,q)=1}}\frac{\mu^2(\ell)}{\varphi(\ell)^2}&=\sum_{\substack{\ell\\(\ell,q)=1}}\frac{\mu^2(\ell)}{\varphi(\ell)^2}-\frac{\mathrm{u}_q\ \mathbf{I}}{X}+O^*\left(\frac{\mathrm{v}_q\ \mathbf{i}^{(q)}}{X^{\frac{3}{2}}}\right)
\end{align} 
where
 \begin{align} 
 \mathrm{u}_q&=\prod_{p|q}\left(1-\frac{p}{p^2-p+1}\right),\ \mathbf{I}=\prod_{p}\left(1+\frac{1}{p(p-1)}\right)\in[\sage{Lower(Prod_sumvarp_l,digits)},\sage{Upper(Prod_sumvarp_u,digits)}],\nonumber\\
 \mathrm{v}_q &=\prod_{p|q}\left(1+\frac{p^2-4p+2}{(\sqrt{p}-1)(p-1)^2+2p-1}\right),\nonumber\\
 \mathbf{i}^{(q)}&=\begin{cases}
\sage{Upper(walue__sumvarp_1,digits)}\prod_{p}\left(1+\frac{2p-1}{(\sqrt{p}-1)(p-1)^2}\right)\phantom{.}\in[\sage{Lower(i_v1_l,digits)},\sage{Upper(i_v1_u,digits)}],\text{ if }2\nmid q,\\
\sage{Upper(walue__sumvarp_2,digits)}\prod_{p}\left(1+\frac{2p-1}{(\sqrt{p}-1)(p-1)^2}\right)\in[\sage{Lower(i_v2_l,digits)},\sage{Upper(i_v2_u,digits)}],\text{ if }2|q,
\end{cases} 
  \end{align}
Therefore, for all $X\geq C=10^6$ and $q=v\in\{1,2\}$, we deduce from \eqref{ready} that
\begin{align}
X\times\sum_{\substack{\ell> X\\(\ell,v)=1}}\frac{\mu^2(\ell)}{\varphi(\ell)^2}&\leq\mathrm{u}_v\ \mathbf{I}+\frac{\mathrm{v}_v\ \mathbf{i}^{(v)}}{\sqrt{C}}\leq
\begin{cases}\label{valuesphi}
\sage{Upper(a_constant_phi_v1,digits)}&\text{ if }v=1,\\
\sage{Upper(a_constant_phi_v2,digits)}&\text{ if }v=2. 
\end{cases}
\end{align}
On the other hand, we have that for all $20 \leq X\leq 10^6$,
\begin{align}
X\times\sum_{\substack{\ell> X\\(\ell,v)=1}}\frac{\mu^2(\ell)}{\varphi(\ell)^2}\leq
\begin{cases}\label{intvaluesphi}
\sage{Upper(constant1_threshold_sumvarp,digits)}&\text{ if }v=1,\\
\sage{Upper(constant2_threshold_sumvarp,digits)}&\text{ if }v=2. 
\end{cases}
\end{align}
Finally, we define $\varphi_v^{(1)}$ by taking the maximum between the bounds \eqref{valuesphi} and \eqref{intvaluesphi}, respectively.
 \end{proof}

\begin{proposition}\label{Ss1}  
Let $X\geq 4\times 10^5$ and $v\in{1,2}$. Then 
\begin{equation*}
\frac{1}{X}\times\sum_{\substack{\ell\leq X\\(\ell,v)=1}}\frac{\mu^2(\ell)\ell^2}{\varphi(\ell)^2}\leq \varphi_v^{(2)}=
\begin{cases}
\sage{Upper(constant2_phi_v1,digits)}&\text{ if }v=1,\\
\sage{Upper(constant2_phi_v2,digits)}&\text{ if }v=2.
\end{cases}
\end{equation*} 
\end{proposition} 

  \begin{proof}
  By using summation by parts in Equation \eqref{ready}, we obtain
 \begin{align*}
\sum_{\substack{\ell\leq X\\(\ell,q)=1}}\frac{\mu^2(\ell)\ell^2}{\varphi(\ell)^2}=\mathrm{u}_q\mathbf{I}\ X+O^*\left(5\mathrm{v}_q\mathbf{i}^{(q)}\ \sqrt{X}\right),
\end{align*}
 where $\mathrm{u}_q,\mathbf{I},\mathrm{v}_q$ and $\mathbf{i}^{(q)}$ are defined in Proposition \ref{sumvarp}. Note that it was not necessary to calculate the main term of \eqref{ready}.
 Hence, when $X\geq C=10^8$ and $q=v\in\{1,2\}$, we have
\begin{equation}
\frac{1}{X}\times\sum_{\substack{\ell\leq X\\(\ell,v)=1}}\frac{\mu^2(\ell)\ell^2}{\varphi(\ell)^2}\leq \mathrm{u}_v\mathbf{I}+\frac{5\mathrm{v}_v\mathbf{i}^{(v)}}{\sqrt{C}}\leq
\begin{cases}\label{valuesphi2}
\sage{Upper(constant2_phi_v1,digits)}&\text{ if }v=1,\\
\sage{Upper(constant2_phi_v2,digits)}&\text{ if }v=2. 
\end{cases} 
\end{equation}
On the other hand, for all $X$ such that $4\times 10^5\leq X\leq 10^8$,
\begin{equation}\label{intvaluesphi2}
\frac{1}{X}\times\sum_{\substack{\ell\leq X\\(\ell,v)=1}}\frac{\mu^2(\ell)\ell^2}{\varphi(\ell)^2}\leq \begin{cases}
\sage{Upper(constant1_threshold_Ss1,digits)}&\text{ if }v=1,\\  
\sage{Upper(constant2_threshold_Ss1,digits)}&\text{ if }v=2.
\end{cases}
\end{equation} 
The result is concluded by taking the maximum between the bounds \eqref{valuesphi2} and \eqref{intvaluesphi2}, which we define as $\varphi_v^{(2)}$, $v\in\{1,2\}$.
 \end{proof}

Consider now the arithmetic function $\nu$ defined on prime numbers as $\nu(2)=1$, $\nu(p)=\frac{p}{p-2}$, if $p>3$. The following result is interesting since it describes the function whose average has an asymptotic expression with constant term equal to $\gamma+\sum_{p|q}\frac{\log(p)}{p-1}$, so that the infinite summation $T_f^q$ considered in \cite[Thms. 3.3, 4.6]{1ZA20} vanishes.

\begin{lemma}\label{Parity} 
Let $X>0$. Then
\begin{align*}
\sum_{\substack{\ell\leq X}}\frac{\mu^2(\ell)\nu(\ell)}{\ell}=\mathbf{H}\left(\log(X)+\gamma+\frac{\log(2)}{2}\right)+O^*\left(\frac{\sage{Upper(h_v1_u,digits)}}{\sqrt{X}}\right),\\
\sum_{\substack{\ell\leq X\\(\ell,2)=1}}\frac{\mu^2(\ell)\nu(\ell)}{\ell}=\frac{\mathbf{H}}{2}\left(\log(X)+\gamma+\log(2)\right)+O^*\left(\frac{\sage{Upper(h_v2_u,digits)}}{\sqrt{X}}\right),
\end{align*}
where $\mathbf{H}=\sage{Numb(ParityProd_l,ParityProd_u,digits)}\ldots$.
\end{lemma}

\begin{proof}
Let $q\in\mathbb{Z}_{>0}$. By \cite[Thm. 4.6]{1ZA20} with $f(\ell)=\frac{\nu(\ell)}{\ell}$, $\alpha=1$ and $\beta=2$, we derive that
\begin{align*}
\sum_{\substack{\ell\leq X\\(\ell,q)=1}}\frac{\mu^2(\ell)\nu(\ell)}{\ell}=\mathrm{m}_q\mathbf{H}\left(\log(X)+\mathfrak{h}_q\right)+O^*\left(\frac{\mathrm{n}_q\ \mathbf{h}^{(q)}}{\sqrt{X}}\right),
\end{align*} 
\begin{flalign} 
\text{where}&&\mathrm{m}_q&=\prod_{2|q}\frac{1}{2}\prod_{\substack{p|q\\p\geq 3}}\left(1-\frac{1}{p-1}\right),\ \mathrm{n}_q=\prod_{\substack{p|q\\p\geq 3}}\left(1-\frac{p-4}{(\sqrt{p}-1)(p-2)+2}\right),&\nonumber\\
&&\mathbf{H}&=\prod_{\substack{p\geq 3}}\left(1+\frac{1}{p(p-2)}\right)\in[\sage{Lower(ParityProd_l,digits)},\sage{Upper(ParityProd_u,digits)}),\label{valueH}&\\ 
&&\mathbf{h}^{(q)}&=
\begin{cases}
\sage{Upper(walue__Parity_1,digits)}\prod_{\substack{p\geq 3}}\left(1+\frac{2}{(\sqrt{p}-1)(p-2)}\right)\in[\sage{Lower(h_v1_l,digits)},\sage{Upper(h_v1_u,digits)}],&\text{ if }2\nmid q,\nonumber\\
\sage{Upper(walue__Parity_2,digits)}\prod_{\substack{p\geq 3}}\left(1+\frac{2}{(\sqrt{p}-1)(p-2)}\right)\in[\sage{Lower(h_v2_l,digits)},\sage{Upper(h_v2_u,digits)}],&\text{ if }2|q,
\end{cases}&
\end{flalign}
and, as $1-f(p)p+2f(p)=1-\frac{p}{p-2}+\frac{2}{p-2}=0$ for $p\geq 3$, $\mathfrak{h}_q=\sum_{2\nmid q}\frac{\log(2)}{2}+\gamma+\sum_{p|q}\frac{\log(p)}{p-1}$.
The result is concluded by considering $q=v\in\{1,2\}$. 
\end{proof} 

\begin{proofbold}{Lemma~\ref{S2:total}}
By recalling \eqref{S2:eq}, we have that $\mathit{S}_{\mathbf{II}}^{(2)}(U)$ is equal to
\begin{align}\label{first:S2}
\frac{v^2}{\varphi(v)^2}\sum_{\substack{\ell\leq Z\\(\ell,v)=1}}\frac{\mu^2(\ell)\ell}{\varphi(\ell)^2}\sum_{\substack{d\\(d,\ell v)=1}}\frac{\mu(d)}{\varphi(d)^2}-\frac{v^2}{\varphi(v)^2}\sum_{\substack{\ell\leq Z\\(\ell,v)=1}}\frac{\mu^2(\ell)\ell}{\varphi(\ell)^2}\sum_{\substack{d>\frac{U}{\ell}\\(d,\ell v)=1}}\frac{\mu(d)}{\varphi(d)^2}=&\\
\frac{v^2}{\varphi(v)^2}\sum_{\substack{\ell\leq Z\\(\ell,v)=1}}\frac{\mu^2(\ell)\ell}{\varphi(\ell)^2}\prod_{p\nmid \ell v}\left(1-\frac{1}{(p-1)^2}\right)-\frac{v^2}{\varphi(v)^2}\sum_{\substack{\ell\leq Z\\(\ell,v)=1}}\frac{\mu^2(\ell)\ell}{\varphi(\ell)^2}\sum_{\substack{d>\frac{U}{\ell}\\(d,\ell v)=1}}\frac{\mu(d)}{\varphi(d)^2}.&\nonumber
\end{align}
The inner sum of the second right hand term of \eqref{first:S2} above can be estimated with the help of Proposition \ref{sumvarp}; indeed, for any $\ell\leq Z$, we have $\frac{U}{\ell}\geq\frac{U}{Z}\geq 20$, so that, by recalling the definition of $\varphi_v^{(1)}$, 
\begin{align}\label{vali1}
\left|\sum_{\substack{d>\frac{U}{\ell}\\(d,\ell v)=1}}\frac{\mu(d)}{\varphi(d)^2}\right|\leq\sum_{\substack{d>\frac{U}{\ell}\\(d,\ell v)=1}}\frac{\mu^2(d)}{\varphi(d)^2}\leq\sum_{\substack{d>\frac{U}{\ell}\\(d,v)=1}}\frac{\mu^2(d)}{\varphi(d)^2}\leq\frac{\varphi^{(1)}_v\ell}{U}.
\end{align}   
Hence,  
\begin{align}
\left|\frac{v^2}{\varphi(v)^2}\sum_{\substack{\ell\leq Z\\(\ell,v)=1}}\frac{\mu^2(\ell)\ell}{\varphi(\ell)^2}\sum_{\substack{d>\frac{U}{\ell}\\(d,\ell v)=1}}\frac{\mu(d)}{\varphi(d)^2}\right|\leq
\frac{\varphi^{(1)}_vv^2}{\varphi(v)^2}\frac{1}{U}\sum_{\substack{\ell\leq Z\\(\ell,v)=1}}\frac{\mu^2(\ell)\ell^2}{\varphi(\ell)^2}\leq\frac{\Upsilon^{(3)}_vZ}{U},\label{vali2} 
\end{align} 
where $\Upsilon^{(3)}_v=\frac{\varphi^{(1)}_v\varphi^{(2)}_vv^2}{\varphi(v)^2}$, by using that $Z\geq 4\times 10^5$ and the definition of $\varphi_v^{(2)}$.

On the other hand, as $\prod_{p\nmid \ell v}\left(1-\frac{1}{(p-1)^2}\right)=0$ if $2\nmid \ell v$ and $\frac{\nu(\ell)}{\ell}=$ $\frac{\mu^2(\ell)\ell}{\varphi(\ell)^2}\prod_{p|\ell}\left(1-\frac{1}{(p-1)^2}\right)^{-1}$  for all square-free numbers $\ell$ such that $(\ell,2)=1$, we have that $\sum_{\substack{\ell\leq Z\\(\ell,v)=1}}\frac{\mu^2(\ell)\ell}{\varphi(\ell)^2}\prod_{p\nmid \ell v}\left(1-\frac{1}{(p-1)^2}\right)$ is equal to 
\begin{align}\label{firstly}
\frac{2}{v}\prod_{p\geq 3}\left(1-\frac{1}{(p-1)^2}\right)\sum_{\substack{\ell\leq \frac{Zv}{2}\\(\ell,2)=1}}\frac{\mu^2(\ell)\nu(\ell)}{\ell}.
\end{align}
Recall now Lemma \ref{Parity}; by using the definition of $\mathbf{H}$ given in \eqref{valueH}, we have 
\begin{align}\label{secondly}
\prod_{p\geq 3}\left(1-\frac{1}{(p-1)^2}\right)\sum_{\substack{\ell\leq \frac{Zv}{2}\\(\ell,2)=1}}\frac{\mu^2(\ell)\nu(\ell)}{\ell}=\phantom{xxxxxxxxx}&\\
\frac{1}{2}\left(\log(Z)+\gamma+\sum_{p|v}\frac{\log(p)}{p-1}\right)+O^*\left(\prod_{p\geq 3}\left(1-\frac{1}{(p-1)^2}\right)\frac{\sqrt{2}\ \sage{Upper(h_v2_u,digits)}}{\sqrt{Zv}}\right).\nonumber&
\end{align}
Therefore, by putting everything together, we derive from estimations \eqref{firstly} and \eqref{secondly} that 
\begin{align}\label{err2_S2}
\frac{v^2}{\varphi(v)^2}\sum_{\substack{\ell\leq Z\\(\ell,v)=1}}\frac{\mu^2(\ell)\ell}{\varphi(\ell)^2}\prod_{p\nmid \ell v}\left(1-\frac{1}{(p-1)^2}\right)=\frac{v}{\varphi(v)}\left(\log(Z)+\gamma+\sum_{p|v}\frac{\log(p)}{p-1}\right)&\\
+O^*\left(\frac{v^2}{\varphi(v)^2}\prod_{p\geq 3}\left(1-\frac{1}{(p-1)^2}\right)\frac{\sqrt{2}\ \sage{Upper(h_v2_u,digits)}}{\sqrt{Zv}}\right).&\nonumber
\end{align}
The result is concluded by defining $\Upsilon^{(2)}_v$ as the resulting bound on the error term given in \eqref{err2_S2}, upon replacing either $v=1$ or $v=2$ and observing that 
\begin{align*} 
\prod_{p\geq 3}\left(1-\frac{1}{(p-1)^2}\right)=\frac{1}{\mathbf{H}}&\in[\sage{Lower(Prod_Parity_l,digits)},\sage{Upper(Prod_Parity_u,digits)}].
\end{align*}
\end{proofbold}

\subsection{The sum $\mathit{S}_{\mathbf{II}}^{(3)}$ and choice of parameter}\label{S23}

As for section \S\ref{S22}, we need a series of results that rely on an interval arithmetic computations that may be carried out for any $q\in\mathbb{Z}_{>0}$.

\begin{proposition}\label{sum.1/2}         
Let $X>0$. Then  
\begin{align} \label{sum1/2general}  
\sum_{\substack{\ell\leq X\\(\ell,q)=1}}\frac{\mu^2(\ell)}{\varphi_{\frac{1}{2}}(\ell)^2}&=\mathrm{f}_q\mathbf{D}\left(\log(X)+\mathfrak{d}_q\right)+O^*\left(\frac{\sage{Upper(Error_sumHalf,digits)}\ \mathrm{h}_q}{X^{\sage{delta_sumHalf}}}\right),\\  
\label{sum2.1/2general} 
\sum_{\substack{\ell\leq X\\(\ell,q)=1}}\frac{\mu^2(\ell)\ell}{\varphi_{\frac{1}{2}}(\ell)^2}&=\mathrm{f}_q\mathbf{D}X+O^*\left(\sage{Upper(Error_sum2Half,digits)}\ \mathrm{h}_q X^{\sage{1-delta_sumHalf}}\right),
\end{align} 
where  
\begin{align*}   
\mathrm{f}_q&=\prod_{p|q}\left(1-\frac{1}{p-2\sqrt{p}+2}\right),\\
\mathfrak{d}_q&=-\sum_{p}\frac{(2\sqrt{p}-3)\log(p)}{(p-2\sqrt{p}+2)(p-1)}+\gamma+\sum_{p|q}\frac{\log(p)}{p-2\sqrt{p}+2},\nonumber\\
\mathrm{h}_q&=\prod_{p|q}\left(1+\frac{p-4\sqrt{p}-p^{\sage{delta_sumHalf}}+2}{(\sqrt{p}-1)^2p^{\sage{1-delta_sumHalf}}+p^{\sage{delta_sumHalf}}+2\sqrt{p}-1}\right),\quad\mathbf{D}=\sage{Numb(Prod_sumHalf_l,Prod_sumHalf_u,3)}\ldots.
\end{align*}
\end{proposition}
 
 \begin{proof} 
  By applying \cite[Thm. 3.3]{1ZA20} with $f(p)=\frac{1}{\varphi_{\frac{1}{2}}(p)^2}=\frac{1}{(\sqrt{p}-1)^2}$, $\alpha=1$, $\beta=\frac{3}{2}$ and $0\leq\delta=\sage{delta_sumHalf}<\frac{1}{2}$, we obtain 
  \begin{align*}
  \frac{\varphi(q)H_{f}^{q}(0)}{q}&=\mathbf{D}\ \prod_{p|q}\left(1-\frac{1}{p-2\sqrt{p}+2}\right),\\
  \frac{\kappa_{1-\delta}(q)\overline{H}_{f}^{\phantom{.}q}(-\delta)}{q^{1-\delta}}&= \mathbf{d}\ \prod_{p|q}\left(1+\frac{p-4\sqrt{p}-p^{\sage{delta_sumHalf }}+2}{(\sqrt{p}-1)^2p^{\sage{1-delta_sumHalf}}+2\sqrt{p}+p^{\sage{delta_sumHalf }}-1}\right).
  \end{align*}  
\begin{flalign*}
 \text{where} &&\mathbf{D}&=\prod_p\left(1+\frac{2}{p(\sqrt{p}-1)}\right) \in [\sage{Lower(Prod_sumHalf_l,digits)},\sage{Upper(Prod_sumHalf_u,digits)}),\\
 && \mathbf{d}&=\Delta_1^{\sage{delta_sumHalf}}\prod_{p}\left(1+\frac{2\sqrt{p}+p^{\sage{delta_sumHalf }}-1}{p^{\sage{1-delta_sumHalf }}(\sqrt{p}-1)^2}\right) \in [\sage{Lower(A(1,delta_sumHalf)*deltaProd3_sumHalf_l,digits)},\sage{Upper(A(1,delta_sumHalf )*deltaProd3_sumHalf_u,digits)}] .&
\end{flalign*}
 On the other hand, but always according to \cite[Thm. 3.2.1]{1ZA20}, we have that
  \begin{equation*}
  T_f^q=-\sum_{p\nmid q}\frac{(2\sqrt{p}-3)\log(p)}{(p-2\sqrt{p}+2)(p-1)},
  \end{equation*}
  so that, by defining $\mathrm{g}_q=\sum_{p|q}\frac{\log(p)}{p-2\sqrt{p}+2}$, we have
  \begin{flalign*}
   &&\mathfrak{d}_q=-\sum_{p}\frac{(2\sqrt{p}-3)\log(p)}{(p-2\sqrt{p}+2)(p-1)}+\gamma+\mathrm{g}_q,\phantom{xxxxxxxxxxxxxxx}\\
 \text{where}&& -\sum_{p}\frac{(2\sqrt{p}-3)\log(p)}{(p-2\sqrt{p}+2)(p-1)}+\gamma\in [\sage{Lower(Sum_sumHalf_l,digits)},\sage{Upper(Sum_sumHalf_u,digits)}],\phantom{xxxxxxxxxxxxxxx}&
  \end{flalign*}
whence Equation \eqref{sum1/2general}. Finally, a summation by parts allows us to derive expression \eqref{sum2.1/2general} from \eqref{sum1/2general}.
 \end{proof} 

The shape of the above error term becomes impractical when one wants to provide an overall estimation. It is there when one can take advantage of computer calculations under interval arithmetic. 

\begin{proposition}\label{sum.1/2threshold}       
Let $X\geq 20$ and $v\in\{1,2\}$. Then  
\begin{align*} 
\frac{1}{\log(X)}\times\sum_{\substack{\ell\leq X\\(\ell,v)=1}}\frac{\mu^2(\ell)}{\varphi_{\frac{1}{2}}(\ell)^2}&\leq\chi_v^{(1)}=
\begin{cases}
\sage{Upper(constant1_chi_v1,digits)}&\text{ if }v=1,\\
\sage{Upper(constant1_chi_v2,digits)}&\text{ if }v=2.
\end{cases}
\end{align*} 
\end{proposition} 
\begin{proof}
Observe that, for all $X$ such that $20\leq X\leq 5\times 10^{8}$, 
\begin{align} 
\frac{1}{\log(X)}\times\sum_{\ell\leq X}\frac{\mu^2(\ell)}{\varphi_{\frac{1}{2}}(\ell)^2}\leq\sage{Upper(constant1_threshold_1half,digits)},\quad\frac{1}{\log(X)}\times\sum_{\substack{\ell\leq X\\(\ell,2)=1}}\frac{\mu^2(\ell)}{\varphi_{\frac{1}{2}}(\ell)^2}\leq\sage{Upper(constant2_threshold_1half,digits)}.\label{value1.1}
\end{align} 
On the other hand, by Proposition \ref{sum.1/2}, When $q=v\in\{1,2\}$ and $X\geq C=5\times 10^8$, we conclude from \eqref{sum1/2general} that
\begin{align}  
\sum_{\substack{\ell\leq X}}\frac{\mu^2(\ell)}{\varphi_{\frac{1}{2}}(\ell)^2}&\leq\sage{Upper(Prod_sumHalf_u,digits)}\left(\log(X)\sage{Upper(Sum_sumHalf_u,digits)}\right)+\frac{\sage{Upper(Error_sumHalf,digits)}} {X^{\sage{delta_sumHalf}}}\nonumber\\
&\leq\left(\sage{Upper(Prod_sumHalf_u,digits)}+\frac{\sage{Upper(Error_sumHalf,digits)}}{C^{\sage{delta_sumHalf}}\log(C)}\right)\log(X)\leq\sage{Upper(a_constant1_chi_v1,digits)}\log(X),\nonumber\\
\sum_{\substack{\ell\leq X\\(\ell,2)=1}}\frac{\mu^2(\ell)}{\varphi_{\frac{1}{2}}(\ell)^2}&\leq\sage{Upper(F2 * Prod_sumHalf_u,digits)}\left(\log(X)\sage{Upper(Sum_sumHalf_u+G2,digits)}\right)+\frac{\sage{Upper(H2 * Error_sumHalf,digits)}}{X^{\sage{delta_sumHalf}}}\label{value1.2}\\
&\leq\left(\sage{Upper(F2 * Prod_sumHalf_u,digits)}+\frac{\sage{Upper(H2 * Error_sumHalf,digits)}}{C^{\sage{delta_sumHalf}}\log(C)}\right)\log(X)\leq\sage{Upper(a_constant1_chi_v2,digits)}\log(X).\nonumber
\end{align} 
The result is concluded by defining $\chi_v^{(1)}$ as the maximum between the bounds \eqref{value1.1} and \eqref{value1.2}.
\end{proof}

\begin{proposition}\label{sum2.1/2threshold}     
Let $X\geq 1$ and $v\in\{1,2\}$. Then 
\begin{align*}
\frac{1}{X}\times\sum_{\substack{\ell\leq X\\(\ell,v)=1}}\frac{\mu^2(\ell)\ell}{\varphi_{\frac{1}{2}}(\ell)^2}&\leq 
\chi_v^{(2)}=\begin{cases}
\sage{Upper(constant2_chi_v1,digits)}&\text{ if }v=1,\\
\sage{Upper(constant2_chi_v2,digits)}&\text{ if }v=2.
\end{cases}
\end{align*}
\end{proposition}

\begin{proof}
Observe that for all $X$ such that $1\leq X\leq 5\times 10^8$, 
\begin{align}
\frac{1}{X}\times\sum_{\ell\leq X}\frac{\mu^2(\ell)\ell}{\varphi_{\frac{1}{2}}(\ell)^2}\leq\sage{Upper(constant1_threshold2_1half,digits)},\quad\frac{1}{X}\times\sum_{\substack{\ell\leq X\\(\ell,2)=1}}\frac{\mu^2(\ell)\ell}{\varphi_{\frac{1}{2}}(\ell)^2}&\leq\sage{Upper(constant2_threshold2_1half,digits)}\label{value1.2}. 
\end{align} 
On the other hand, by Proposition \ref{sum.1/2}, when $q=v\in\{1,2\}$ and $X\geq C=5\times 10^8$, \eqref{sum2.1/2general} tells us that
\begin{align}  
\sum_{\substack{\ell\leq X}}\frac{\mu^2(\ell)\ell}{\varphi_{\frac{1}{2}}(\ell)^2}&\leq\left(\sage{Upper(Prod_sumHalf_u,digits)}+\frac{\sage{Upper(Error_sum2Half,digits)}}{C^{\sage{delta_sumHalf}}} \right)X\leq \sage{Upper(a_constant2_chi_v1,digits)}\times X,\nonumber\\ 
\sum_{\substack{\ell\leq X\\(\ell,2)=1}}\frac{\mu^2(\ell)\ell}{\varphi_{\frac{1}{2}}(\ell)^2}&\leq\left(\sage{Upper(Prod_sumHalf_u,digits)}\ \mathrm{f}_2+\frac{\sage{Upper(Error_sum2Half,digits)}\times\mathrm{h}_2 }{C^{\sage{delta_sumHalf}}}\right)X\leq \sage{Upper(a_constant2_chi_v2,digits)}\times X\label{value2.2}.
\end{align}   
The result is thus obtained by defining $\chi_v^{(2)}$ as the maximum between the bounds \eqref{value1.2} and \eqref{value2.2}.
\end{proof}

\begin{proposition}\label{sum2.1/2Logthreshold}    
Let $X$, $Z$ such that $1\leq Z<X$ and $v\in\{1,2\}$. Then
\begin{align*}
\sum_{\substack{\ell\leq Z\\(\ell,v)=1}}\frac{\mu^2(\ell)\ell}{\varphi_{\frac{1}{2}}(\ell)^2}\log\left(\frac{X}{\ell}\right)&\leq\chi_v^{(2)}\ Z\left(\log\left(\frac{X}{Z}\right) + 1\right),
\end{align*} 
where $\chi_v^{(2)}$ is defined as in Proposition \ref{sum2.1/2threshold}.
\end{proposition}

\begin{proof} 
By summation by parts and Proposition \ref{sum2.1/2threshold}, we derive 
\begin{align*}
\sum_{\substack{\ell\leq Z\\(\ell,v)=1}}\frac{\mu^2(\ell)\ell}{\varphi_{\frac{1}{2}}(\ell)^2}\log\left(\frac{X}{\ell}\right)&\leq \chi_v^{(2)}\ Z\log\left(\frac{X}{Z}\right) + \chi_v^{(2)}\ (Z-1).
\end{align*} 
\end{proof}
 
By \cite[Thm. 3.3]{1ZA20} and summation by parts, we know how to the detect the order of the summation given in Lemma \ref{Ss1Log}: by applying Proposition \ref{logint}, it is of order $\log^{-2}(X)$; by  using \cite[Thm. 4.6]{1ZA20}, summation by parts and Proposition \ref{logint}, we could derive an estimation for it. Nonetheless, given that this sum is involved in a small value, presented in the second term of the bounds \eqref{errorbound1},  we have chosen to proceed faster by observing that the non-weighted sum in the statement below is convergent. 

\begin{proposition}\label{Ss1Log}  
Let $X\geq 10^{12}$, $\theta=1-\frac{1}{\log(10^{12})}$ and $v\in\{1,2\}$. Then
\begin{align*}
\log^2(X)\times\sum_{\substack{d\leq\frac{X}{10^{12}}\\(d,v)=1}}\frac{\mu^2(d)}{d^{2-2\theta}\varphi_\theta(d)^2}\frac{1}{\log^2\left(\frac{X}{d}\right)}\leq\tau_v=
\begin{cases}
\sage{Upper(constant_tau_v1,digits)}\quad&\text{ if }v=1,\\ 
\sage{Upper(constant_tau_v2,digits)}\quad&\text{ if }v=2.\\
\end{cases}
\end{align*}  
\end{proposition}

\begin{proof} 
Define $f(p)=\frac{p^{2\theta}}{(p^\theta-1)^2}$ on prime numbers and extend it to a multiplicative function. Consider $W$ such that $1<W<X$;
by using the bound $\log^{-1}\left(\frac{X}{d}\right)\leq\log^{-1}(W)$, for $1\leq d\leq\frac{X}{W}$, and writing $\log(X)=\log\left(\frac{X}{d}\right)+\log(d)$, we derive that $\sum_{\substack{d\leq\frac{X}{W}\\(d,q)=1}}\frac{\mu^2(d)f(d)\log^2(X)}{d^2\log^2\left(\frac{X}{d}\right)}$ is bounded from above by 
\begin{align}\label{ineq:5}
\sum_{\substack{d\leq\frac{X}{W}\\(d,q)=1}}\frac{\mu^2(d)f(d)}{d^2}&+\sum_{\substack{d\leq\frac{X}{W}\\(d,q)=1}}\frac{\mu^2(d)f(d)\log(d)}{d^2}\left(\frac{2}{\log(W)}+\frac{\log(d)}{\log^2(W)}\right).
\end{align} 
Observe now that the functions $t\mapsto\frac{\log(t)}{\sqrt{t}}$ and $t\mapsto\frac{\log^2(t)}{\sqrt{t}}$ have a global maximum at $t=e^2$, with value $\frac{\log(e^2)}{\sqrt{e^2}}=\frac{2}{e}$, and at $t=e^4$, with value $\frac{16}{e^2}$, respectively. Hence, from \eqref{ineq:5}, $\sum_{\substack{d\leq\frac{X}{W}\\(d,q)=1}}\frac{\mu^2(d)f(d)\log^2(X)}{d^2\log^2\left(\frac{X}{d}\right)}$ is at most
\begin{align}
\prod_{p\nmid q}\left(1+\frac{f(p)}{p^2}\right)+\left(\frac{4}{e\log(W)}+\frac{16}{e^2\log^2(W)}\right)\prod_{p\nmid q}\left(1+\frac{f(p)}{p^{\frac{3}{2}}}\right).\label{ineq:6}
\end{align}
On the other hand, on using the definition of $f$, we obtain  
\begin{align*}
\prod_{p\nmid q}\left(1+\frac{f(p)}{p^2}\right)&=\prod_{p|q}\left(1-\frac{1}{p^{2-2\theta}(p^{\theta}-1)^2+1}\right)\prod_{p}\left(1+\frac{1}{p^{2-2\theta}(p^{\theta}-1)^2}\right),\\
\prod_{p\nmid q}\left(1+\frac{f(p)}{p^{\frac{3}{2}}}\right)&=\prod_{p|q}\left(1-\frac{1}{p^{\frac{3}{2}-2\theta}(p^{\theta}-1)^2+1}\right)\prod_{p}\left(1+\frac{1}{p^{\frac{3}{2}-2\theta}(p^{\theta}-1)^2}\right),
\end{align*}
where
\begin{align*}
\prod_{p}\left(1+\frac{1}{p^{2-2\theta}(p^{\theta}-1)^2}\right)&\in[\sage{Lower(Prod2_Ss1Log_l,digits)},\sage{Upper(Prod2_Ss1Log_u,digits)}],\\
\prod_{p}\left(1+\frac{1}{p^{\frac{3}{2}-2\theta}(p^{\theta}-1)^2}\right)&\in[\sage{Lower(Prod3_Ss1Log_l,digits)},\sage{Upper(Prod3_Ss1Log_u,digits)}] .
\end{align*}
Thus, by replacing $q=v\in\{1,2\}$ and $W=10^{12}$ into \eqref{ineq:6} and by using the above infinite product estimations, we obtain the value of $\tau_v$.
\end{proof} 

\begin{proposition}\label{Ss2Log}    
Let $X\geq 10^{12}$, $\theta=1-\frac{1}{\log(10^{12})}$ and $q\in\mathbb{Z}_{>0}$. Let $\mathrm{c}$ and $\varepsilon$ be two real numbers such that $1<Z=\mathrm{c} X^{\varepsilon}<X$ and $0<\varepsilon<1-\frac{\log(\mathrm{c})}{\log(10^{12})}$. Then
\begin{align*}
\log(X)\times\sum_{\substack{\ell\leq Z\\(\ell,q)=1}}\frac{\mu^2(\ell)}{\ell\log^2\left(\frac{X}{\ell}\right)}\left(\frac{\ell^\theta}{\varphi_\theta(\ell)}\right)^2\leq\xi_{\{\mathrm{c},\varepsilon,q\}},
\end{align*} 
for some explicit constant $\xi_{\{\mathrm{c},\varepsilon,q\}}>0$. In particular, we may define $\xi_1^{(\sage{ccc})}=\xi_{\{\sage{ccc},\sage{choice},1\}}=\sage{Upper(constant_xi_v1,digits)}$, $\xi_2^{(\sage{ccc})}=\xi_{\{\sage{ccc},\sage{choice},2\}}=\sage{Upper(constant_xi_v2,digits)}$.
, $\xi_2^{(\sage{ccc2})}=\xi_{\{\sage{ccc2},\sage{choice},2\}}=
\sage{Upper(constant2_xi_v2,digits)}$ and $\xi_1^{(\sage{ccc3})}=\xi_{\{\sage{ccc3},\sage{choice},1\}}=
\sage{Upper(constant2_xi_v1,digits)}$.
\end{proposition}

\begin{proof} 
Define $f(p)=\frac{1}{p^{1-2\theta}(p^\theta-1)^2}$. As $\theta>\frac{1}{2}$, we can use \cite[Thm. 4.6]{1ZA20} with $\alpha=1$ and $\beta=1+\theta$ to estimate the above sum without the weight $\ell\mapsto\log^{-1}\left(\frac{X}{\ell}\right)$. We derive that $\sum_{\substack{\ell\leq Z\\(\ell,q)=1}}\frac{\mu^2(\ell)}{\ell}\left(\frac{\ell^\theta}{\varphi_\theta(\ell)}\right)^2$ may be estimated as
\begin{align}\label{ineq:8}
\mathrm{M}_q(Z)+O^*\left(\frac{\mathrm{y}_q\ \mathbf{j}^{(q)}}{\sqrt{Z}}\right)=\mathrm{x}_q\ \mathbf{J}\left(\log(Z)+\mathfrak{j}_q\right)+O^*\left(\frac{\mathrm{y}_q\ \mathbf{j}^{(q)}}{\sqrt{Z}}\right),
\end{align}
where
 \begin{flalign*}
\mathrm{x}_q&=\prod_{p|q}\left(1-\frac{p-1}{p^{2-2\theta}(p^\theta-1)^2+2p^{1-\theta}-p^{1-2\theta}-1}\right),\\
\mathfrak{j}_q&=-\sum_{p}\frac{\log(p)(2p^{1-\theta}-p^{1-2\theta}-2)}{(p^{1-2\theta}(p^\theta-1)^2+1)(p-1)}+\gamma+\sum_{p|q}\frac{\log(p)}{p^{1-2\theta}(p^\theta-1)^2+1},\phantom{xxx}\\
\mathrm{y}_q&=\prod_{p|q}\left(1+\frac{p^{2\theta}-4p^\theta+2}{(\sqrt{p}-1)(p^\theta-1)^2+2p^\theta-1}\right),\\    
\mathbf{J}&=\prod_{p}\left(1+\frac{2p^{1-\theta}-p^{1-2\theta}-1}{p^{2-2\theta}(p^\theta-1)^2}\right)\in[\sage{Lower(Prod_Ss2Log_l,digits)},\sage{Upper(Prod_Ss2Log_u,digits)}],\\
-&\sum_{p}\frac{\log(p)(2p^{1-\theta}-p^{1-2\theta}-2)}{(p^{1-2\theta}(p^\theta-1)^2+1)(p-1)}+\gamma\in[\sage{Lower(Sum_Ss2Log_l,digits)},\sage{Upper(Sum_Ss2Log_u,digits)}],\\
\mathbf{j}^{(q)}&=\begin{cases}
\sage{Upper(walue_Ss2Log_1,digits)}\prod_{p}\left(1+\frac{2p^\theta-1}{(\sqrt{p}-1)(p^\theta-1)^2}\right)&\in[\sage{Lower(j_v1_l,digits)},\sage{Upper(j_v1_u,digits)}],\quad\text{ if }2\nmid q,\\
\sage{Upper(walue_Ss2Log_2,digits)}\prod_{p}\left(1+\frac{2p^\theta-1}{(\sqrt{p}-1)(p^\theta-1)^2}\right)&\in[\sage{Lower(j_v2_l,digits)},\sage{Upper(j_v2_u,digits)}],\quad\text{ if }2|q.
\end{cases} 
 \end{flalign*}
Now, by integration by parts, we derive from \eqref{ineq:8} that the summation $\sum_{\substack{\ell\leq Z\\(\ell,q)=1}}\frac{\mu^2(\ell)}{\ell\log^2\left(\frac{X}{\ell}\right)}\left(\frac{\ell^\theta}{\varphi_\theta(\ell)}\right)^2$ equals
\begin{align}  
\int_1^Z\frac{\mathrm{M}_q\phantom{}'(t)}{\log^2\left(\frac{X}{t}\right)}dt+\frac{\mathrm{M}_q(1)}{\log^2(X)}+\mathrm{y}_q\ \mathbf{j}^{(q)}O^*\left(\frac{1}{\sqrt{Z}\log^2\left(\frac{X}{Z}\right)}+\int_1^Z\frac{2dt}{t^{\frac{3}{2}}\log^3\left(\frac{X}{t}\right)}\right),\label{ineq:7}
\end{align}
Notice that the integral in right hand side of \eqref{ineq:7} can be bounded with the help of Proposition \ref{logint} $\mathbf{b)}$, giving an upper bound for the above error term of the form 
\begin{align} 
\mathrm{y}_q\ \mathbf{j}^{(q)}\left(\frac{1}{\sqrt{Z}\log^2\left(\frac{X}{Z}\right)}+\frac{4}{\log^3\left(\frac{X}{\sqrt{Z}}\right)}+\frac{4}{Z^{\sage{kk*1/2}}\log^3\left(\frac{X}{Z}\right)}\right)\label{trick1}\phantom{xxxxxxxxx}&\\ 
\leq\frac{\mathrm{y}_q\ \mathbf{j}^{(q)}\mathrm{W}_{\{\mathrm{c},\varepsilon\}}}{\log\left(\frac{10^{12(1-\varepsilon)}}{\mathrm{c}}\right)\left(1-\varepsilon-\frac{\log(\mathrm{c})}{\log(10^{12})}\right)\log(X)}=\frac{\mathrm{Y}_{\{\mathrm{c},\varepsilon,q\}}}{\log(X)},\nonumber &
\end{align} 
where we have used that $0<\varepsilon<1-\frac{\log(\mathrm{c})}{\log(10^{12})}$, so that $1<\frac{10^{12(1-\varepsilon)}}{\mathrm{c}}\leq\frac{X}{Z}$ and where $\mathrm{W}_{\{\mathrm{c},\varepsilon\}}$ is defined as
\begin{align*} 
\frac{1}{\sqrt{\mathrm{c}}10^{6\varepsilon}}+\frac{4}{\sage{kk}\log\left(\frac{10^{12(\sage{1/kk}-\varepsilon)}}{\mathrm{c}}\right)}+\frac{4}{\mathrm{c}^{\sage{kk*1/2}}10^{\sage{kk*1/2*12}\varepsilon}\log\left(\frac{10^{12(1-\varepsilon)}}{\mathrm{c}}\right)}.
\end{align*}  

On the other hand, the main term of Equation \eqref{ineq:7} can be bounded with the help of Proposition \ref{logint} $\mathbf{a)}$, giving that $\int_1^Z\frac{\mathrm{M}_q\phantom{}'(t)}{\log^2\left(\frac{X}{t}\right)}dt+\frac{\mathrm{M}_q(1)}{\log^2(X)}$ is at most
\begin{align} 
\mathrm{x}_q\ \mathbf{J}\left(\frac{1}{\log\left(\frac{X}{Z}\right)}-\frac{1}{\log(X)}+\frac{\mathfrak{j}_q}{\log^2(X)}\right)\label{trick2}\leq\frac{\mathrm{X}_{\{\mathrm{c},\varepsilon,q\}}}{\log(X)},
\end{align}
where we have used that $\mathfrak{j}_q$ is positive for all $q\in\mathbb{Z}_{>0}$ and
\begin{align*}
\mathrm{X}_{\{\mathrm{c},\varepsilon,q\}}=\mathrm{x}_q\ \mathbf{J}\left(\frac{1}{1-\varepsilon-\frac{\log(\mathrm{c})}{\log(10^{12})}}+\frac{\mathfrak{j}_q}{\log(10^{12})}\right).
\end{align*} 
The result is achieved by considering $q=v\in\{1,2\}$, setting $\xi_{\{\mathrm{c},\varepsilon,v\}}$ as any upper bound of $\mathrm{X}_{\{\mathrm{c},\varepsilon,v\}}+\mathrm{Y}_{\{\mathrm{c},\varepsilon,v\}}$ and defining $\xi_{v}^{\mathrm{c}}$ as a reasonable value of $\xi_{\{\mathrm{c},\sage{choice},v\}}$, for $\mathrm{c}\in\{\sage{ccc},\sage{ccc2},\sage{ccc3}\}$.  
\end{proof}  
 
\begin{remark}\label{remark1} As long as $1<Z<X$, and regardless of the choice of $Z$, the expressions \eqref{trick1} and \eqref{trick2} tell us that $\sum_{\substack{\ell\leq Z\\(\ell,q)=1}}\frac{\mu^2(\ell)}{\ell\log^2\left(\frac{X}{\ell}\right)}\left(\frac{\ell^\theta}{\varphi_\theta(\ell)}\right)^2\ll_{q}\frac{1}{\log\left(\frac{X}{Z}\right)}$.
\end{remark}

As $\mathit{S}_{\mathbf{II}}^{(3)}$ involves a summation of small terms, given by the extraction from $\mathit{S}_{\mathbf{II}}$ of the main terms of the functions $(\ell,d)\mapsto\check{m}_{\ell dv}$, we expect it to be small. As it turns out, this is so provided that $U$ is sufficiently large. 

By using the inequality $(A_1+B_1)^2\leq (1+\omega)A_2^2+\left(1+\frac{1}{\omega}\right)B_2^2$, valid for any $0<A_1\leq A_2$, $0<B_1\leq B_2$ and $\omega>0$, we obtain that $|\mathit{S}_{\mathbf{II}}^{(3)}|$ is at most 
\begin{align} 
\label{rama}
\frac{1+\omega_v}{U}\frac{v}{\varphi_{\frac{1}{2}}(v)^2}\sum_{\substack{\ell\leq Z\\(\ell,v)=1}}\frac{\mu^2(\ell)}{\ell}\sum_{\substack{d\leq\frac{U}{\ell}\\(d,\ell v)=1}}\frac{\mu^2(d)}{d^2}\left(\frac{\ell d}{\varphi_{\frac{1}{2}}(\ell d)}\right)^2+\phantom{xxxxxxxxxx}&\\
\label{rama2}
\frac{1}{389^2}\left(1+\frac{1}{\omega_v}\right)\frac{v^{2\theta}}{\varphi_\theta(v)^2}\sum_{\substack{\ell\leq Z\\(\ell,v)=1}}\frac{\mu^2(\ell)}{\ell}\sum_{\substack{d\leq\frac{U}{\ell}\\(d,\ell v)=1}}\frac{\mu^2(d)}{d^2}\left( \frac{(\ell d)^\theta}{\varphi_{\theta}(\ell d)}\frac{\mathds{1}_{\{\frac{U}{\ell d}>10^{12}\}}^{(d,\ell)}}{\log\left(\frac{U}{\ell d}\right)}\right)^2.&
\end{align} 
On the other hand, as $20\leq\frac{U}{Z}\leq\frac{U}{\ell}$ we can apply Proposition \ref{sum.1/2threshold} and then Proposition \ref{sum2.1/2Logthreshold} to bound the double sum given in \eqref{rama} as 
\begin{align}  
\sum_{\substack{\ell\leq Z\\(\ell,v)=1}}\frac{\mu^2(\ell)\ell}{\varphi_{\frac{1}{2}}(\ell)^2}\sum_{\substack{d\leq\frac{U}{\ell}\\(d,\ell v)=1}}\frac{\mu^2(d)}{\varphi_{\frac{1}{2}}(d)^2}
&\leq\sum_{\substack{\ell\leq Z\\(\ell,v)=1}}\frac{\mu^2(\ell)\ell}{\varphi_{\frac{1}{2}}(\ell)^2}\sum_{\substack{d\leq\frac{U}{\ell}\\(d,v)=1}}\frac{\mu^2(d)}{\varphi_{\frac{1}{2}}(d)^2}\nonumber\\
&\leq \chi_v^{(3)}Z\left(\log\left(\frac{U}{Z}\right)+1\right),&\label{part:1},
\end{align}
where $\chi_v^{(3)}= \chi_v^{(1)} \chi_v^{(2)}$. Furthermore, Proposition \ref{Ss1Log} with $X=\frac{U}{\ell}\geq 10^{12}$ gives the following estimation  
\begin{align*}
\sum_{\substack{d\leq\frac{U}{\ell 10^{12}}\\(d,\ell v)=1}}\frac{\mu^2(d)}{d^{2-2\theta}\varphi_\theta(d)^2}\frac{1}{\log^2\left(\frac{U}{\ell d}\right)}\leq\sum_{\substack{d\leq\frac{U}{\ell 10^{12}}\\(d,v)=1}}\frac{\mu^2(d)}{d^{2-2\theta}\varphi_\theta(d)^2}\frac{1}{\log^2\left(\frac{U}{\ell d}\right)}\leq\frac{\tau_v\mathds{1}_{\{U\geq 10^{12}\}}(U)}{\log^2\left(\frac{U}{\ell}\right)},
\end{align*}
so that, when $U\geq 10^{12}$, $\sum_{\substack{\ell\leq Z\\(\ell,v)=1}}\frac{\mu^2(\ell)}{\ell}\sum_{\substack{d\leq\frac{U}{\ell}\\(d,\ell v)=1}}\frac{\mu^2(d)}{d^2}\left( \frac{(\ell d)^\theta}{\varphi_{\theta}(\ell d)}\frac{\mathds{1}_{\{\frac{U}{\ell d}>10^{12}\}}^{(d,\ell)}}{\log\left(\frac{U}{\ell d}\right)}\right)^2$ is at most
\begin{align}\label{part:2} 
\tau_v\sum_{\substack{\ell\leq Z\\(\ell,v)=1}}\frac{\mu^2(\ell)}{\ell\log^2\left(\frac{U}{\ell}\right)}\left(\frac{\ell^\theta}{\varphi_\theta(\ell)}\right)^2.
\end{align}

\noindent\textbf{Choice of parameter.}  
We conclude from estimations \eqref{part:1}, \eqref{part:2} and Remark \ref{remark1} that, as long as $\frac{U}{Z}\geq 20$ (or rather, $\frac{U}{Z}\gg 1$), we have
\begin{align}
\mathit{S}_{\mathbf{II}}^{(3)}\ll_v \frac{Z}{U}\log\left(\frac{U}{Z}\right)+ \frac{Z}{U}+\frac{1}{\log\left(\frac{U}{Z}\right)}. \label{S3_ll}
\end{align}
Therefore, by recalling \eqref{split}, we derive from lemmas \ref{S1:total}, \ref{S2:total} and Equation \eqref{S3_ll} that the error term of $\sum_{\substack{d,e\\(de,v)=1}}\frac{\mu(d)\mu(e)}{[d,e]}\mathbf{L}_d\mathbf{L}_e=\mathit{S}_{\mathbf{I}}+\mathit{S}_{\mathbf{II}}$ has magnitude at most 
\begin{align}
\frac{\log^4\left(\frac{U}{Z}\right)}{\sqrt{Z}}+\frac{\log^2\left(\frac{U}{Z}\right)}{\sqrt{Z}}+\frac{1}{\sqrt{Z}}+\frac{Z}{U}+\frac{Z\log\left(\frac{U}{Z}\right)}{U}+\frac{1}{\log\left(\frac{U}{Z}\right)}\label{exhibit}. 
\end{align}
From the bound $\frac{\log^4\left(\frac{U}{Z}\right)}{\sqrt{Z}}\leq\frac{\log^4\left(U\right)}{\sqrt{Z}}$, as we are aiming for an error term as small as possible, we may suppose that $\log^\mathrm{a}(U)\ll Z$ for all $\mathrm{a}>0$.  On the other hand, in order to minimize the contribution of $\log^{-1}\left(\frac{U}{Z}\right)$ to the order of \eqref{exhibit}, given that $\log^{-1}(U)\ll\log^{-1}\left(\frac{U}{Z}\right)$ and that $\log^{-1}\left(\frac{U}{\log^{\mathrm{a}}(U)U^\varepsilon}\right)$ is of strictly higher order than $\log^{-1}(U)$ for any $\mathrm{a}>0$ and $0<\varepsilon<1$, it is plausible to suppose that $Z=O(U^\varepsilon)$ for some $0<\varepsilon<1$. Hence the overall magnitude of the expression given in \eqref{exhibit}, regardless of coprimality conditions, is $\log^{-1}(U)$, as predicted in Theorem \ref{MAI}; furthermore, the magnitude of the secondary terms therein is at most $\log^4\left(\frac{U}{Z}\right)\left(\frac{1}{\sqrt{Z}}+\frac{Z}{U}\right)=O(\log^4(U)\left(U^{-\frac{\varepsilon}{2}}+U^{-(1-\varepsilon)}\right))$, and it achieves its minimal order when $\varepsilon=\sage{choice}$. We can set then $Z=\mathrm{c}U^{\sage{choice}}$ for some constant $\mathrm{c}>0$.   

\begin{proofbold}{Lemma~\ref{S3:total}} 
Let $\mathrm{c}\in\{\sage{ccc},\sage{ccc2},\sage{ccc3}\}$.
As $0<\varepsilon=\sage{choice}<1-\frac{\log(\sage{max(ccc,ccc2,ccc3)})}{\log(10^{12})}$, we can apply the particular case of Proposition \ref{Ss2Log} to estimate \eqref{part:2}, considered only when $U\geq 10^{12}$. Moreover, by equations \eqref{rama}, \eqref{rama2} and estimation \eqref{part:1}, we have that $|\mathit{S}_{\mathbf{II}}^{(3)}|$ is bounded from above by
\begin{align}\label{OMEGA} 
(1+\omega_v)\frac{\chi_v^{(1)} \chi_v^{(2)}v}{\varphi_{\frac{1}{2}}(v)^2}\left(\frac{\log(U)}{\sage{1/(1-choice)}U^{\sage{1-choice}}}+\frac{1}{U^{\sage{1-choice}}}\right)+\left(1+\frac{1}{\omega_v}\right)\frac{\tau_v\xi_v^{(\mathrm{c})}v^{2\theta}}{\varphi_\theta(v)^2}\frac{\mathds{1}_{\{U\geq 10^{12}\}}(U)}{389^2\log(U)}.
\end{align}
On the other hand, observe that  
\begin{align*}  
(1+\omega_v)\frac{\chi_v^{(1)} \chi_v^{(2)}v}{\varphi_{\frac{1}{2}}(v)^2}\geq\frac{\chi_v^{(1)} \chi_v^{(2)}v}{\varphi_{\frac{1}{2}}(v)^2}\geq
\begin{cases} 
\sage{Lower(boundomega_v1_1,digits)}&\text{ if }v=1,\\
\sage{Lower(boundomega_v2_1,digits)}&\text{ if }v=2, 
\end{cases}\\
\frac{1}{389}\frac{\tau_v\xi_v^{(\mathrm{c})}v^{2\theta}}{\varphi_\theta(v)^2}\leq\max_{\mathrm{c}\in\{\sage{ccc},\sage{ccc2},\sage{ccc3}\}}\left\{\frac{1}{389}\frac{\tau_v\xi_v^{(\mathrm{c})}v^{2\theta}}{\varphi_\theta(v)^2}\right\}\leq
\begin{cases} 
\sage{Upper(max(boundomega_v1_2,boundomega2_v1_2,boundomega3_v1_2),digits)}&\text{ if }v=1,\\
\sage{Upper(max(boundomega_v2_2,boundomega2_v2_2,boundomega3_v2_2),digits)}&\text{ if }v=2.
\end{cases}
\end{align*}
Hence, the estimation for the first term of  \eqref{OMEGA} is numerically very big and we seek for a value $\omega_v$ that will not make it much more bigger upon multiplying it by $(1+\omega_v)$. On the other hand, as long as $\frac{1}{389}\left(1+\frac{1}{\omega_v}\right)$ remains small, the second term estimation therein will be numerically small.  Thus, we may naturally suppose that $(1+\omega_v)\frac{\chi_v^{(1)} \chi_v^{(2)}v}{\varphi_{\frac{1}{2}}(v)^2}\leq \frac{1}{v}+\frac{\chi_v^{(1)} \chi_v^{(2)}v}{\varphi_{\frac{1}{2}}(v)^2}$ and $\frac{1}{389}\left(1+\frac{1}{\omega_v}\right)\leq\frac{1}{v}$ and, for the sake of simplicity, we may set $\omega_v=\frac{\varphi_{\frac{1}{2}}(v)^2}{\chi_v^{(1)} \chi_v^{(2)}v^2}$, which is independent of the value of $\mathrm{c}$, and define
\begin{align*} 
\Upsilon^{(4)}_v=\frac{(1+\omega_v)\chi_v^{(1)} \chi_v^{(2)}v}{\varphi_{\frac{1}{2}}(v)^2},\quad\Upsilon^{(5)}_v=\Upsilon^{(5)}_{\{\mathrm{c},v\}}=\frac{1}{389^2}\left(1+\frac{1}{\omega_v}\right)\frac{\tau_v\xi_v^{(\mathrm{c})}v^{2\theta}}{\varphi_\theta(v)^2}.
\end{align*}
\end{proofbold}

\section{Main term and conclusion}\label{Conclusion}

We choose three values of $\mathrm{c}$: if $U\geq 10^{\sage{exp_prog}}$, set $Z=\sage{ccc}U^{\sage{choice}}$; if $U\geq 10^{\sage{exp_prog_bound}}$ and $v=2$, set $Z=\sage{ccc2}U^{\sage{choice}}$ and, if $U\geq 10^{\sage{exp_prog_bound}}$ and $v=1$, set $Z=\sage{ccc3}U^{\sage{choice}}$. With these choices, we have that $\frac{U}{2Z}\geq 10$ and $Z\geq 4\times 10^5$. Thus, conditions on propositions \ref{sq1.half}, \ref{sumvar1log}, \ref{sumvarp} are satisfied with $X=\frac{U}{2Z}$ and estimations \eqref{source} and \eqref{source2}, \eqref{ethos} and, \eqref{vali1}, are respectively correct. Furthermore, as $Z\geq 4\times 10^5$, Proposition \ref{Ss1} can be applied to derive inequality \eqref{vali2}. Finally, Proposition \ref{sum.1/2threshold} holds with $X=\frac{U}{Z}\geq 20$, so that we can derive \eqref{part:1}.    

By recalling identity \eqref{split}, we can combine lemmas \ref{SI:estimation}, \ref{S1:total}, \ref{S2:total} and \ref{S3:total} to derive the following estimation, 
\begin{align}\label{combine} 
\sum_{\substack{d,e\\(de,v)=1}}\frac{\mu(d)\mu(e)}{[d,e]}\log^+\left(\frac{U}{d}\right)\log^+\left(\frac{U}{e}\right)=\frac{v}{\varphi(v)}\log(U)-\mathfrak{s}_v+\Xi_v(U),
\end{align} 
where $\mathfrak{s}_v=\frac{v}{\varphi(v)}\left(\gamma+\sum_{p|v}\frac{\log(p)}{p-1}-\frac{6}{\pi^2}\frac{\varphi(v)}{\kappa(v)}\int_1^{\infty}\frac{h_v(s)}{s}ds\right)$, $v\in\{1,2\}$ and
\begin{align} \label{definitive}
 |\Xi_v(U)|\leq\frac{T_v^{(1)}\log^4\left(\frac{U}{\mathrm{c}^{\sage{1/(1-choice)}}}\right)}{\sage{1/(1-choice)}^4\sqrt{\mathrm{c}}U^{\sage{choice/2}}}
+\frac{2\Upsilon^{(1)}_v\log^2\left(\frac{U}{\mathrm{c}^{\sage{1/(1-choice)}}}\right)}{\sage{1/(1-choice)}^2\sqrt{\mathrm{c}}U^{\sage{choice/2}}}
+\frac{\sage{6*ccc/(1/(1-choice))}\Psi_v\ v}{\pi^2\kappa(v)}\frac{\log\left(\frac{U}{\mathrm{c}^{\sage{1/(1-choice)}}}\right)}{U^{\sage{1-choice}}}\phantom{xx}&\nonumber\\
+\frac{\Upsilon^{(4)}_v\log(U)}{\sage{1/(1-choice)}U^{\sage{1-choice}}}
+\frac{\sage{6*ccc}\Psi_v\ v}{\pi^2\kappa(v)\ U^{\sage{1-choice}}}
+\frac{\Upsilon^{(2)}_v}{\sqrt{\mathrm{c}}U^{\sage{choice/2}}}
+\frac{\mathrm{c}\Upsilon^{(3)}_v}{U^{\sage{1-choice}}}+\frac{\Upsilon^{(4)}_v}{U^{\sage{1-choice}}}\nonumber&\\
+\frac{\sage{12/(1-choice)}T_v^{(4)}v}{\pi^2\kappa(v)\left(1-\frac{\log\left(\mathrm{c}^{\sage{1/(1-choice)}}\right)}{\log(U)}\right)}\frac{1}{\log(U)}
+\frac{\mathds{1}_{\{U\geq 10^{12}\}}(U)\Upsilon^{(5)}_{\{\mathrm{c},v\}}}{\log(U)}.&  
\end{align}  

If $U\geq 10^{\sage{exp_prog}}$, we use that $\left(1-\frac{\log\left(\mathrm{c}^{\sage{1/(1-choice)}}\right)}{\log(U)}\right)^{-1}\leq\left(1-\frac{\log\left(\sage{ccc}^{\sage{1/(1-choice)}}\right)}{\log(10^{\sage{exp_prog}})}\right)^{-1}$ to merge the last two terms of estimation \eqref{definitive} to the order $\frac{1}{\log(U)}$. Moreover, by using that  $\log\left(\frac{U}{Z}\right)\leq\sage{1-choice}\log(U)$, we can merge the remaining lower order terms to the second order $\frac{\log^4(U)}{U^{\sage{choice/2}}}$. 

If $U\geq 10^{\sage{exp_prog_bound}}$, we can obtain a numerical bound for $\Xi_v(U)$, $v\in\{1,2\}$, by observing that the last term in \eqref{definitive} appears only if $U\geq 10^{12}$, that $\left(1-\frac{\log\left(\mathrm{c}_v^{\sage{1/(1-choice)}}\right)}{\log(U)}\right)^{-1}\leq\left(1-\frac{\log\left(\mathrm{c}_v^{\sage{1/(1-choice)}}\right)}{\log(10^{\sage{exp_prog_bound}})}\right)^{-1}$, with $\mathrm{c}_1=\sage{ccc3}$, $\mathrm{c}_2=\sage{ccc2}$, and that the functions $t\geq\mathrm{c}_v^{\sage{1/(1-choice)}}\mapsto\log^4\left(\frac{t}{\mathrm{c}_v^{\sage{1/(1-choice)}}}\right)t^{-\sage{choice/2}}$, $t\geq\mathrm{c}_v^{\sage{1/(1-choice)}}\mapsto\log^2\left(\frac{t}{\mathrm{c}_v^{\sage{1/(1-choice)}}}\right)t^{-\sage{choice/2}}$, $t\geq\mathrm{c}_v^{\sage{1/(1-choice)}}\mapsto\log\left(\frac{t}{\mathrm{c}_v^{\sage{1/(1-choice)}}}\right)t^{-\sage{choice/2}}$ and $t>0\mapsto\frac{\log(t)}{t^{\sage{choice/2}}}$ 
are all decreasing for $t\geq 10^{\sage{exp_prog_bound}}$.   Subsequently, by observing that $t\geq\mathrm{c}_v^{\sage{1/(1-choice)}}\mapsto\log^5\left(\frac{t}{\mathrm{c}_v^{\sage{1/(1-choice)}}}\right)t^{-\sage{choice/2}}$ is also decreasing for $t\geq 10^{\sage{exp_prog_bound}}$, we can merge the whole expression \eqref{definitive} to the order $\frac{1}{\log(U)}$.

Furthermore, by observing that the results given in \S\ref{A}, \S\ref{S21},  \S\ref{S22} and  \S\ref{S23} can be worked out for any fixed $q\in\mathbb{Z}_{>0}$, we can derive a general equation, similar to \eqref{combine}, in which $v$ is replaced by $q$, and whose remainder term $\Xi_q(U)$ is of order $\frac{1}{\log(U)}$.  Thus, by putting everything together, we deduce our main result.

\begin{theorem}\label{CONCLUSION}
Let $U>1$. Then for all $q\in\mathbb{Z}_{>0}$, one can determine a explicit constant $K_q>0$ such that
\begin{align*}
\sum_{\substack{d,e\\(de,q)=1}}\frac{\mu(d)\mu(e)}{[d,e]}\log^+\left(\frac{U}{d}\right)\log^+\left(\frac{U}{e}\right)=&\frac{q}{\varphi(q)}\log(U)-\mathfrak{s}_q+\Xi_q(U),
\end{align*}  
and $\Xi_q(U)=O_q^*\left(\frac{K_q}{\log(U)}\right)$, where, by recalling the definition of $h_q$ given in \eqref{h_q}, we have
\begin{equation*}
\mathfrak{s}_q=\frac{q}{\varphi(q)}\left(\gamma+\sum_{p|q}\frac{\log(p)}{p-1}-\frac{6}{\pi^2}\frac{\varphi(q)}{\kappa(q)}\int_1^{\infty}\frac{h_q(s)}{s}ds\right).
\end{equation*} 
In particular, if $U\geq 10^{\sage{exp_prog}}$, we have 
\begin{align}\label{errorbound1}  
|\Xi_1(U)|\leq 
\frac{\sage{Upper(Alt_Merge1_v1+Alt_Merge2_v1,digits+1)}\log^4\left(U\right)}{U^{\sage{choice/2}}}+\frac{\sage{Upper(MergeLog_v1,digits+1)}}{\log(U)},\quad
|\Xi_2(U)|\leq
\frac{\sage{Upper(Alt_Merge1_v2+Alt_Merge2_v2,digits+1)}\log^4\left(U\right)}{U^{\sage{choice/2}}}+\frac{\sage{Upper(MergeLog_v2,digits+1)}}{\log(U)}.
\end{align}  
If $U\geq 10^{\sage{exp_prog_bound}}$, we have
\begin{align}\label{numericalbound1} 
|\Xi_1(U)|&\leq 
\sage{Upper(numerical_v1,digits+2)},
&&|\Xi_2(U)|\leq 
\sage{Upper(numerical_v2,digits+2)},\\
\label{mergebound1} 
|\Xi_1(U)|&\leq
\frac{\sage{Upper(KV1,digits+1)}}{\log(U)},
&&|\Xi_2(U)|\leq
\frac{\sage{Upper(KV2,digits+1)}}{\log(U)}.
\end{align} 
\end{theorem}

\noindent \textbf{The specific constant in the case $v\in\{1,2\}$}. An estimation of $\int_1^{10^6}\frac{h_v(s)}{s}ds$ is provided in \cite[Proposition 6.26]{Hel19}. Inspired by this calculation, we obtain
\begin{align*} 
\int_1^{10^8}\frac{h_v(s)}{s}ds\in
\begin{cases} [-0.0495100113498 -0.049510010626]\ &,\text{if\ }v=1,\\
[2.63481269161,2.63481271383]\ &,\text{if\ }v=2.
\end{cases}
\end{align*}
Moreover, similar to the estimation given in \eqref{hv22}, we have that $\int_1^{\infty}\frac{h_v(s)}{s}ds$ can be expressed as
\begin{align*} 
\int_1^{10^{\sage{bound}}}\frac{h_v(s)}{s}ds+O^*\left(\Psi_v'\Omega+\frac{2T_v^{(4)}}{\log(10^{12})}\right),
\end{align*}
where $\Psi_v'=T_v^{(2)}+\frac{T_v^{(3)}}{\log(10^{\sage{bound}})}$ and $\Omega=\frac{\log(10^{\sage{bound}})}{10^{\sage{bound}}}-\frac{\log(10^{12})}{10^{12}}+\frac{1}{10^{\sage{bound}}}-\frac{1}{10^{12}}$. Thereupon, the constant $\mathfrak{s}_v$ coming from Theorem \ref{CONCLUSION} can be estimated as follows
\begin{align}
\mathfrak{s}_1\in  \gamma+\frac{6}{\pi^2}\times[0.049510010626,0.0495100113498]+O^*(\sage{sci11})&\nonumber\\
\in[\sage{Lower(C_constant_s1_l,11)},\sage{Upper(C_constant_s1_u,11)}],\phantom{xxxxxxxxxxxxxxxxxxxxxxx.}&\nonumber\\
\mathfrak{s}_2\in 2(\gamma+\log(2))+\frac{4}{\pi^2}\times[-2.63481271383,-2.63481269161]\phantom{xxxxxxx.}&\nonumber\\
+O^*(\sage{sci22})\phantom{xxx}&\nonumber\\
\in[\sage{Lower(C_constant_s2_l,11)},\sage{Upper(C_constant_s2_u,11)}],\phantom{xxxxxxxxxxxxxxxxxxxxxxx.}&\label{svestimations} 
\end{align} 
so that \eqref{hope} is confirmed.   

Finally, since we have carefully estimated the value of $\mathfrak{s}_v$, $v\in\{1,2\}$, by assuming that $U\geq 10^{\sage{exp_prog}}$, and since $10^{\sage{exp_prog}}$ is a moderate value, we can ask about the behavior of \eqref{combine} when $1\leq U\leq 10^{\sage{exp_prog}}$. Indeed, this study was carried out by a program, inspired by \cite[\S 6.6.1]{Hel19}, that run for about $6$ weeks combined with a recurrence algorithm, whose output we write below.

\begin{proposition}\label{SUComputation}  
Let $1\leq U\leq 10^{\sage{exp_prog}}$ and $v\in\{1,2\}$. Then the expression \\ $\left|\sum_{\substack{d,e\\(de,v)=1}}\frac{\mu(d)\mu(e)}{[d,e]}\log^+\left(\frac{U}{d}\right)\log^+\left(\frac{U}{e}\right)-\frac{v}{\varphi(v)}\log(U)+\mathfrak{s}_v\right|$ is at most $\sage{Upper(barrier_threshold_1,digits+1)}$ $\times U^{-\frac{1}{3}}$, if $v=1$, or $\sage{Upper(barrier_threshold_2,digits+1)}$ $\times U^{-\frac{1}{3}}$, if $v=2$.
\end{proposition} 
\noindent Note the similarity between the error constants above and the corresponding (interval) values of $\mathfrak{s}_1$, $\mathfrak{s}_2$. Indeed, those are the values that arise when $U$ is very close to $3^+$ and $1^+$, respectively, the error constants being much smaller for larger $U$.

\section {On the Brun--Tichmarsh theorem}\label{Application}

By proceeding as in  \cite[\S 3.2]{MV07}, we finally see how we can apply Theorem \ref{CONCLUSION} to derive an explicit result about the distribution of the prime numbers. 

Let $X,Y\in\mathbb{R}_{\geq 0}$ , $P\in\mathbb{Z}_{>0}$ and define $S(X,Y;P)=\#\{n\in(X,X+Y]\cap\mathbb{Z}_{>0},(n,P)=1\}$. As per the discussion in \S\ref{Intro}, we know that $\left\{\frac{\mathbf{L}_d}{\log(U)}\right\}_{d=1}^\infty$, $U>2$, is a sequence of parameters as in the Selberg sieve. Hence, by taking in particular $P'=\prod_{2<p\leq U}p$, by \cite[Eq. (3.10)]{MV07}, we readily have that $S(X,Y;2P')$ is at most
\begin{equation*}  
\frac{Y}{2\log^2(U)}\sum_{\substack{d,e\\(de,2)=1}}\frac{\mu(d)\mu(e)}{[d,e]}\mathbf{L}_d\mathbf{L}_e+\frac{4}{\log^2(U)}\left(\sum_{\substack{d\\(d,2)=1}}|\mu^2(d)\mathbf{L}_d|\right)^2.
\end{equation*}
Observe that the above right hand term can be estimated as 
\begin{equation*}
\sum_{\substack{d\\(d,2)=1}}|\mu^2(d)\mathbf{L}_d|=\int_0^U\left(\sum_{\substack{d\leq t\\(d,2)=1}}\mu^2(d)\right)\frac{dt}{t}=\frac{4}{\pi^2}U+O^*\left(\iota\sqrt{U}\right),
\end{equation*}
where, by using \cite[Lemma 5.2]{Hel19}, the constant $\iota$ is defined as $2\left(1-\frac{4}{\pi^2}\right)$. 
 
Therefore, by using Theorem \ref{CONCLUSION} and estimation \eqref{svestimations}, we have that $S(X,Y;2P')$ is bounded from above by
\begin{equation}\label{primebound}
\frac{Y}{2\log^2(U)}\left(2\log(U)-\sage{Lower(C_constant_s2_l,digits+2)}+\Xi_2(U)\right)+\frac{4}{\log^2(U)}\left(\frac{4}{\pi^2}U+\iota\sqrt{U}\right)^2.  
\end{equation}
From \eqref{primebound}, we see immediately that in order to derive a non-trivial estimation for $S(X,Y;2P')$, and in general, to derive a continuous and monotonic version of Selberg sieve, we must have $U^2\ll Y$. Moreover, the magnitude of the bound $\eqref{primebound}$ is $\frac{Y}{\log(U)}$, which is minimized when $U$ is as large as possible. For numerical simplicity, we take $U=\sqrt{Y}$. Therefore \eqref{primebound} can be written as
\begin{equation}\label{primestar}
\frac{2Y}{\log(Y)}+\frac{4Y}{\log^2(Y)}\left(-\sage{Lower(C_constant_s2_l/2,digits+2)}+\frac{\Xi_2(\sqrt{Y})}{2}+4\left(\frac{4}{\pi^2}+\frac{\iota}{Y^{\frac{1}{4}}}\right)^2\right).  
\end{equation} 
By supposing that $Y\geq 10^{\sage{ebrun}}$, we may use \eqref{numericalbound1} to bound $\Xi_2(\sqrt{Y})$. Moreover, by using the bound $\frac{4}{\pi^2}+\iota Y^{-\frac{1}{4}}\leq\frac{4}{\pi^2}+\iota 10^{-\sage{ebrun/4}}$, we derive from \eqref{primestar} that
\begin{equation}\label{star}
S(X,Y;2P')\leq\frac{2Y}{\log(Y)}-\frac{\sage{-Upper(Brun,digits+2)}Y}{\log^2(Y)}  . 
\end{equation}
In particular, if $Y\geq 10^{\sage{ebrun}}$, by \cite[Eq. (3.3)]{MV07} and \eqref{star}, we derive that for all $X\geq 0$,
\begin{align}\label{distribution}
\pi(X+Y)-\pi(X)&\leq\frac{2Y}{\log(Y)}-\frac{\sage{-Upper(Brun,digits+2)}Y}{\log^2(Y)} +\sqrt{Y}\leq\frac{2Y}{\log(Y)}\left(1-\frac{\sage{-Upper(constantprime,digits+1)}}{\log(Y)}\right), 
\end{align}
where we have used that $\frac{\log^2(Y)}{\sqrt{Y}}\leq\frac{\log^2(10^{\sage{ebrun}})}{\sqrt{10^{\sage{ebrun}}}}$.

Now, if $P$ is any positive integer, not necessarily of the form $2P'$ considered in \eqref{star}, we can derive from \cite[Thm. 3.6]{MV07} and \eqref{star} that for any $Y\geq 10^{\sage{ebrun}}$,
\begin{equation}\label{generalstar}
S(X,Y,P)\leq\prod_{\substack{p\leq\sqrt{Y}\\p\nmid P}}\left(1-\frac{1}{p}\right)^{-1}\times\left(\frac{2Y}{\log(Y)}-\frac{\sage{-Upper(Brun,digits+2)}Y}{\log^2(Y)}\right).
\end{equation}
On the other hand, suppose that $Y\geq 10^{\sage{ebrun}}q$, where $a,q,P\in\mathbb{Z}_{>0}$ are such that $(a,q)=(P,q)=1$. We derive from \cite[Thm. 3.8]{MV07} that the number of integers $n$ such that $X<n\leq X+Y$, $n\equiv a\ (mod\ q)$ and $(n,P)=1$ can be bounded from above by
\begin{equation}\label{BT1}
\prod_{\substack{p\leq\sqrt{\frac{Y}{q}}\\p\nmid P}}\left(1-\frac{1}{p}\right)^{-1}\times\left(\frac{2Y}{q\log\left(\frac{Y}{q}\right)}-\frac{\sage{-Upper(Brun,digits+2)}Y}{q\log^2\left(\frac{Y}{q}\right)}\right).
\end{equation}
Therefore, we can generalize \eqref{distribution} as follows.

\begin{theorem}[\textbf{Brun--Titchmarsh inequality}]\label{BT-thm} 
Let $a,q\in\mathbb{Z}_{>0}$ such that $(a,q)=1$ and let $Y$ be a real number such that $Y\geq 10^{\sage{ebrun}}q$. Then for all $X\geq 0$,
\begin{equation*}
\pi(X+Y;q,a)-\pi(X;q,a)\leq\frac{2Y}{\varphi(q)\log\left(\frac{Y}{q}\right)}\left(1-\frac{\sage{-Upper(constantprime,digits+1)}}{\log\left(\frac{Y}{q}\right)}\right).
\end{equation*}  
\end{theorem}
\begin{proof} 
As $(a,q)=1$, if $p\equiv a\ (mod\ q)$, then $p\nmid q$. Proceed as in \cite[Thm. 3.9]{MV07}. By selecting $P=\prod_{\substack{p\nmid q,p\leq\sqrt{\frac{Y}{q}}}}p$, the number of primes $p$ such that $X<p\leq X+Y$, $p\equiv a\ (mod\ q)$ and $(p,P)=1$ can be bounded with the help of \eqref{BT1}, giving the upper bound
\begin{align}\label{BT2}
\prod_{\substack{p\leq\sqrt{\frac{Y}{q}}\\p|q}}\left(1-\frac{1}{p}\right)^{-1}\times\left(\frac{2Y}{q\log\left(\frac{Y}{q}\right)}-\frac{\sage{-Upper(Brun,digits+2)}Y}{q\log^2\left(\frac{Y}{q}\right)}\right)\phantom{xxxxxxx}\nonumber&\\
\leq\frac{Y}{\varphi(q)\log\left(\frac{Y}{q}\right)}\left(2-\frac{\sage{-Upper(Brun,digits+2)}}{\log\left(\frac{Y}{q}\right)}\right).&
\end{align} 
It remains to consider the primes $X<p\leq X+Y$ such that $p\equiv a\ (mod\ q)$ and $p|P$; by definition of $P$, its cardinality is at most $\sqrt{\frac{Y}{q}}$. Hence, by using that $\frac{\varphi(q)\log^2\left(\frac{Y}{q}\right)}{q\sqrt{\frac{Y}{q}}}\leq\frac{\log^2(10^{\sage{ebrun}})}{\sqrt{10^{\sage{ebrun}}}}$, we conclude from \eqref{BT2} that $\pi(X+Y;q,a)-\pi(X;q,a)$ is at most
\begin{align}\label{BT3} 
\frac{Y}{\varphi(q)\log\left(\frac{Y}{q}\right)}\left(2-\frac{\sage{-Upper(Brun,digits+2)}}{\log\left(\frac{Y}{q}\right)}+\frac{\varphi(q)\log^2\left(\frac{Y}{q}\right)}{q\sqrt{\frac{Y}{q}}\log\left(\frac{Y}{q}\right)}\right)\phantom{xxxxxxx}\nonumber&\\
\leq\frac{Y}{\varphi(q)\log\left(\frac{Y}{q}\right)}\left(2-\frac{\sage{-Upper(2*constantprime,digits+1)}}{\log\left(\frac{Y}{q}\right)}\right).&
\end{align}
\end{proof}
The reader may refer to \cite{MA13}, \cite{MV73} and \cite{SE} for further insights about the Brun--Titchmarsh theorem.

\section{Computational matters}\label{calcul}
 
The procedure for calculating converging expressions or constants in this article is standard: we set a precision value, calculate them by a recurrence algorithm until its precision value, and then, if needed, we numerically bound the remaining term, which will be small. These calculations were carried out under the interval arithmetic header int$\_$double14.2.h, implemented by Platt in C\texttt{++} (used for example in \cite{Pla16}), providing results with double precision, higher performance and faster speed, when compared to the ARB package of Sage.  

The constant $\mathrm{c}$, introduced in \S\ref{S23}, plays an important  role in deriving Theorem \ref{CONCLUSION}. As discussed in \S\ref{Conclusion}, supposing that $U\geq 10^E$, for $E\geq\sage{exp_prog}$, the parameter $Z=\mathrm{c}_EU^{\sage{choice}}$ must satisfy $U\geq 20Z$ and $Z\geq 4\times 10^5$. Also, by Proposition \ref{Ss2Log}, we must have $\mathrm{c}_E<10^{\sage{12*(1-choice)}}$.  Therefore, we take $\mathrm{c}_E\in \left[\frac{4\times 10^5}{10^{\sage{choice}E}},\min\left\{\frac{10^{\frac{E}{\sage{1/(1-choice)}}}}{20},10^{\sage{12*(1-choice)}}\right\}\right)$. In particular, if $E=\sage{exp_prog}$, $\mathrm{c}_E\in(\sage{Lower(4*10^5/10^(exp_prog*choice),digits)},\sage{Upper(10^(exp_prog*(1-choice))/20,digits)})$, 
whereas, if $E=\sage{exp_prog_bound}$, $\mathrm{c}_E\in(\sage{Lower(4*10^5/10^(exp_prog_bound*choice),digits)},\sage{Upper(10^(exp_prog_bound*(1-choice))/20,digits)})$. With the latter choice, as we want to obtain a numerical bound in \eqref{definitive}, we additionally suppose that the three functions $t\geq 10^{\sage{exp_prog_bound}}\mapsto\log^A\left(\frac{t}{c_E^3}\right)t^{-\sage{1-choice}}$, $A\in\{1,2,4\}$ are decreasing, so that $c_E\leq\frac{10^\frac{E}{\sage{1/(1-choice)}}}{e^A}$ and $c_E\in(\sage{Lower(4*10^5/10^(exp_prog_bound*choice),digits)},\sage{Upper(10^(exp_prog_bound*(1-choice))/e^4,digits)})$.  Further, in order to merge \eqref{definitive} to the order $\frac{1}{\log(U)}$, we consider $A\in\{1,2,3,5\}$, so that $c_E\in(\sage{Lower(4*10^5/10^(exp_prog_bound*choice),digits)},\sage{Upper(10^(exp_prog_bound*(1-choice))/e^5,digits)})$.

For the sake of simplicity, we choose $\mathrm{c}$ to be integer. The value $\mathrm{c}=\sage{ccc}$ considered in the range $U\geq  10^{\sage{exp_prog}}$, has been chosen to diminish the numerical contribution of the term of magnitude $\frac{\log^4(U)}{U^{\sage{1-choice}}}$ to the error \eqref{errorbound1}. The value $\mathrm{c}=\sage{ccc2}$, considered in the range $U\geq  10^{\sage{exp_prog_bound}}$, has been chosen to provide 
the optimal value $\mathpzc{B}$ in Theorem \ref{BT-thm}, given our exposition. Moreover, in this range, the value $\mathrm{c}=\sage{ccc2}$, with $v=2$, is also optimal for merging the bound \eqref{definitive} to the expression $\frac{\mathpzc{x}_{2}}{\log(U)}$, as $\mathrm{c}$ contributes to $\mathpzc{B}$ only through $\Xi_2(U)$, introduced in \eqref{numericalbound1}; by the same reasoning, when $v=1$,  we choose $\mathrm{c}=\sage{ccc3}$ to merge the bound \eqref{definitive} to the expression $\frac{\mathpzc{x}_{1}}{\log(U)}$, described in \eqref{mergebound1}.
See the table below. 

	\begin{table}[htb] 
    \centering
      \vspace{0.01mm}
		\begin{tabular}{|| M{0.35cm} | M{1.4cm}  | M{1.3cm} ||  M{0.35cm} | M{1.4cm}  | M{1.2cm} || M{0.5cm} | M{1.5cm}  | M{1.2cm} ||} 
		
		\hline
$\mathrm{c}$         & $\mathpzc{B}$     & $\mathpzc{x}_{1}$                              & $\mathrm{c}$	        & $\mathpzc{B}$  & $\mathpzc{x}_{1}$                & $\mathrm{c}$	        & $\mathpzc{B}$ & $\mathpzc{x}_{1}$			 \\ \hline 

$0.5$ & $-0.8264$ &  $\sage{Trunc(Xiv1_value_1,ds2)}\phantom{0}$ & $\sage{CCC_5}$ & $\phantom{-}\sage{Trunc(Brun_value_5,ds)}$ & $\sage{Trunc(Xiv1_value_5,ds2)}\phantom{0}$ &  $\sage{CCC_9}$ & $\sage{Trunc(Brun_value_9,ds)}$ & $\sage{Trunc(Xiv1_value_9,ds2)}$  \\ \hline 

$\sage{CCC_2}$ & $\sage{Trunc(Brun_value_2,ds)}$ &  $\phantom{1}\sage{Trunc(Xiv1_value_2,ds2)}$ & $\sage{CCC_6}$ & $\phantom{-}\sage{Trunc(Brun_value_6,ds)}$ &  $\sage{Trunc(Xiv1_value_6,ds2)}$ & $\sage{CCC_10}$ & $-0.0199$ & $\sage{Trunc(Xiv1_value_10,ds2)}$ \\ \hline 

$\sage{CCC_3}$ & $\phantom{-}\sage{Trunc(Brun_value_3,ds)}$ &  $\phantom{1}\sage{Trunc(Xiv1_value_3,ds2)}$ & $\sage{CCC_7}$ & $\phantom{-}\sage{Trunc(Brun_value_7,ds)}$ &  $\sage{Trunc(Xiv1_value_7,ds2)}$ & $\sage{CCC_11}$ & $\sage{Trunc(Brun_value_11,ds)}$ & $\sage{Trunc(Xiv1_value_11,ds2)}\phantom{0}$ \\ \hline 

$\sage{CCC_4}$ & $\phantom{-}\sage{Trunc(Brun_value_4,ds)}\phantom{0}$ &  $\phantom{1}\sage{Trunc(Xiv1_value_4,ds2)}$ & $\sage{CCC_8}$ & $\sage{Trunc(Brun_value_8,ds)}\phantom{0}$ &  $\sage{Trunc(Xiv1_value_8,ds2)}$ & $\sage{CCC_12}$ & $\sage{Trunc(Brun_value_12,ds)}$ & $-$ \\ \hline 
		
          	\end{tabular} 
		\end{table} 

\section*{Acknowledgements}

A simplified version of this work belonged to my PhD thesis. I would like to render my warmest thanks to H. Helfgott and M. Hindry for insightful discussions during my doctoral studies. 
Thanks are also due to C. Dartyge, R. de la Bret\`eche, O. Ramar\'e and T. Trudgian for providing thorough comments on my earlier work.

\end{document}